%% file: qcovering.tex
\documentclass[leqno]{amsart}
\usepackage{amsmath,amssymb,amscd,mathrsfs}
\usepackage{graphicx,color,ifthen,comment}
\usepackage{mathptmx,courier}
\usepackage[scaled=0.95]{helvet}

\usepackage[cmtip,arrow]{xy}
\usepackage{pb-diagram,pb-xy,lamsarrow,pb-lams}

\specialcomment{Notes}{\begin{quote}\footnotesize}{\end{quote}}
\excludecomment{Notes}

\usepackage{ifpdf}
\ifpdf
\usepackage[pdftex,
            a4paper,
            pdfauthor={Michael Eisermann}
            pdftitle={Quandle coverings and their Galois correspondence}
            pdfcreator={LaTeX with hyperref},
            pdfproducer={LaTeX with hyperref},
            pdfstartview=FitH,
            bookmarksopen=true,
            bookmarksnumbered,
            colorlinks=true, 
            anchorcolor=black, 
            linkcolor=black, 
            urlcolor=black,
            citecolor=black, 
            pagecolor=black, 
            filecolor=black, 
            menucolor=black, 
            ]{hyperref}
\else
\usepackage[dvips]{hyperref}
\fi

\title{Quandle coverings and their Galois correspondence}

\author{Michael Eisermann}
\address{Institut Fourier, Universit\'e Grenoble I, France}
\email{Michael.Eisermann@ujf-grenoble.fr}
\urladdr{www-fourier.ujf-grenoble.fr/{\textasciitilde}eiserm}

\thanks{This work was carried out during the winter term 2006/2007
  while the author was on a sabbatical leave funded by a research contract
  \textit{d\'el\'egation aupr\`es du CNRS}, whose support is gratefully acknowledged}

\date{first version November 2006; this version compiled \today}

\theoremstyle{plain}
  \newtheorem{theorem}{Theorem}[section]
  \newtheorem{lemma}[theorem]{Lemma}
  \newtheorem{proposition}[theorem]{Proposition}
  \newtheorem{corollary}[theorem]{Corollary}
  
\theoremstyle{definition}
  \newtheorem{definition}[theorem]{Definition}
  
  \newtheorem{remark}[theorem]{Remark}
  \newtheorem{example}[theorem]{Example}
  \newtheorem{warning}[theorem]{Warning}
  \newtheorem*{notation}{Notation}

\newcommand{\secref}[1]{\textsection\ref{#1}}

\newcommand{\N}{\mathbb{N}}
\newcommand{\Z}{\mathbb{Z}}
\newcommand{\Q}{\mathbb{Q}}
\newcommand{\R}{\mathbb{R}}
\newcommand{\C}{\mathbb{C}}

\renewcommand{\S}{\mathbb{S}}
\newcommand{\RP}{\mathbb{RP}}
\newcommand{\T}{\mathbb{T}}

\newcommand{\K}{\mathbb{K}}

\newcommand{\GL}{\operatorname{GL}}
\newcommand{\SL}{\operatorname{SL}}
\newcommand{\PSL}{\operatorname{PSL}}

\newcommand{\g}{\mathfrak{g}}
\newcommand{\bil}[2]{\langle #1 , #2 \rangle}

\newcommand{\cs}{\mathbin{\sharp}}
\newcommand{\tsa}{\mathbin{\overline{\ast}}}

\newcommand{\ad}{\operatorname{ad}}
\newcommand{\gr}{\operatorname{gr}}

\newcommand{\id}{\operatorname{id}}
\newcommand{\pr}{\operatorname{pr}}
\newcommand{\minus}{\smallsetminus}
\newcommand{\into}{\hookrightarrow}
\newcommand{\onto}{\mathrel{\makebox[1pt][l]{$\to$}{\to}}}
\newcommand{\isoto}{\mathrel{\xrightarrow{_\sim}}}
\newcommand{\gen}[1]{\langle{#1}\rangle}

\newcommand{\length}{\operatorname{length}}

\newcommand{\Ker}{\operatorname{ker}}
\renewcommand{\Im}{\operatorname{im}}
\newcommand{\Map}{\operatorname{Map}}
\newcommand{\Hom}{\operatorname{Hom}}
\newcommand{\End}{\operatorname{End}}
\newcommand{\Aut}{\operatorname{Aut}}
\newcommand{\Sym}{\operatorname{Sym}}

\newcommand{\Inn}{\operatorname{Inn}}
\newcommand{\inn}{\operatorname{inn}}
\newcommand{\Conj}{\operatorname{Conj}}
\newcommand{\Core}{\operatorname{Core}}
\newcommand{\Alex}{\operatorname{Alex}}
\newcommand{\Adj}{\operatorname{Adj}}
\newcommand{\adj}{\operatorname{adj}}
\newcommand{\Ext}{\operatorname{Ext}}

\newcommand{\category}[1]{\operatorname{\bf{#1}}}
\newcommand{\Groups}{\category{Grp}}
\newcommand{\Quandles}{\category{Qnd}}
\newcommand{\Coverings}{\category{Cov}}
\newcommand{\connected}{*}
\newcommand{\Subgroups}{\category{Sub}}
\newcommand{\Actions}{\category{Act}}
\newcommand{\lto}[1][]{\mathrel{\smash{\overset{\smash{#1}}{\longrightarrow}}}}
\newcommand{\lot}[1][]{\mathrel{\smash{\overset{\smash{#1}}{\longleftarrow}}}}
\newcommand{\lonto}[1][]{\mathrel{\smash{\overset{\smash{#1}}%
{\makebox[1pt][l]{$\longrightarrow$}{\longrightarrow}}}}}

\begin{document} 

\dedicatory{Dedicated to Professor Egbert Brieskorn,
on the occasion of his 70th birthday}

\begin{abstract}
  This article establishes the algebraic covering theory of quandles. 
  For every connected quandle $Q$ with base point $q \in Q$, we explicitly 
  construct a universal covering $p \colon (\tilde{Q},\tilde{q}) \to (Q,q)$.
  This in turn leads us to define the algebraic fundamental group 
  $\pi_1(Q,q) := \Aut(p) = \{ g \in \Adj(Q)' \mid q^g = q \}$,
  where $\Adj(Q)$ is the adjoint group of $Q$. 
  We then establish the Galois correspondence between 
  connected coverings of $(Q,q)$ and subgroups of $\pi_1(Q,q)$.
  Quandle coverings are thus formally analogous to coverings 
  of topological spaces, and resemble Kervaire's algebraic covering
  theory of perfect groups.  A detailed investigation also reveals 
  some crucial differences, which we illustrate by numerous examples.
  
  As an application we obtain a simple formula 
  for the second (co)homology group of a quandle $Q$.  
  It has long been known that $H_1(Q) \cong H^1(Q) \cong \Z[\pi_0(Q)]$,
  and we construct natural isomorphisms $H_2(Q) \cong \pi_1(Q,q)_\mathrm{ab}$ 
  and $H^2(Q,A) \cong \Ext(Q,A) \cong \Hom(\pi_1(Q,q),A)$, reminiscent of the classical 
  Hurewicz isomorphisms in degree $1$.  This means that whenever $\pi_1(Q,q)$ 
  is known, (co)homology calculations in degree $2$ become very easy.  
\end{abstract}

\keywords{classification of quandle coverings, 
  algebraic fundamental group of a quandle, 
  groupoid coverings and Galois correspondence,
  non-abelian cohomology theory}

\copyrightinfo{2006}{Michael Eisermann}

\subjclass[2000]{
57M25, 
20L05, 
18B40, 
18G50  
}



\maketitle



\vspace{-5mm}

\makeatletter
\begin{quote}
  \renewcommand{\contentsline}[4]{#2}
  \renewcommand{\tocsection}[3]{%
    \ifthenelse{\equal{#2}{}}{\unskip\ignorespaces}{\item[#2.] \bf#3\unskip\@addpunct.\rm }}
  \renewcommand{\tocsubsection}[3]{%
    \ifthenelse{\equal{#2}{}}{\unskip\ignorespaces}{#2.~#3\unskip\@addpunct. }}
  \renewcommand{\tocsubsubsection}[3]{\hspace{-3pt}\unskip\ignorespaces}
  \footnotesize 
  \begin{center}\contentsnamefont\contentsname\end{center}
  \medskip
  \setTrue{toc}
  \@input{\jobname.toc}%
  \if@filesw
    \@xp\newwrite\csname tf@toc\endcsname
    \immediate\@xp\openout\csname tf@toc\endcsname \jobname.toc\relax
  \fi
  \global\@nobreakfalse \par
  \addvspace{32\p@\@plus14\p@}%
  \let\tableofcontents\relax
\end{quote}
\makeatother


\section{Introduction and outline of results} \label{sec:Introduction}

\subsection{Motivation and background} \label{sub:IntroMotivation}

In every group $(G,\cdot)$ one can define conjugation 
on the right $a \ast b = b^{-1} \cdot a \cdot b$, 
and its inverse, conjugation on the left $a \tsa b = b \cdot a \cdot b^{-1}$.
They enjoy the following properties for all $a,b,c\in G$:
\begin{enumerate}
\item[(Q1)]
  $a \ast a = a$ \hfill (idempotency)
\item[(Q2)]
  $(a \ast b) \tsa b = a = (a \tsa b) \ast b$ \hfill (right invertibility)
\item[(Q3)]
  $(a\ast b)\ast c = (a\ast c)\ast (b\ast c)$ \hfill (self-distributivity)
\end{enumerate}

Turning these properties into axioms, D.\,Joyce \cite{Joyce:1982} defined 
a \emph{quandle} to be a set $Q$ equipped with two binary operations 
$\ast,\tsa \colon Q\times Q\to Q$ satisfying (Q1--Q3). 
Quandles thus encode the algebraic properties of conjugation;
this axiomatic approach is most natural for studying situations 
where group multiplication is absent or of a secondary nature
(see \secref{sec:Examples} for examples).
Slightly more general, a \emph{rack} is only required 
to satisfy (Q2--Q3); such structures appear naturally 
in the study of braid actions (Brieskorn \cite{Brieskorn:1988}) 
and provide set-theoretic solutions of the Yang-Baxter equation
(Drinfel{$'$}d \cite{Drinfeld:1990}).

The main motivation to study quandles comes from knot theory:
the Wirtinger presentation of the fundamental group $\pi_K = \pi_1(\S^3 \minus K)$
of a knot or link $K \subset \S^3$ involves only conjugation but not
the group multiplication itself, and can thus be seen to define a quandle $Q_K$.
The three quandle axioms then correspond precisely to the three Reidemeister moves.
These observations were first explored in 1982 by Joyce \cite{Joyce:1982}, 
who showed that the knot quandle $Q_K$ classifies knots up to orientation.
Many authors have since rediscovered and studied this notion. 
(See the historical remarks in \secref{sub:HistoricalRemarks}.)

In the 1990s emerged the concept of rack and quandle 
(co)homology \cite{FennRourkeSanderson:1995}, and 
it has since been put to work in constructing combinatorial 
knot invariants \cite{CarterEtAl:2001,CarterEtAl:2003,CarterEtAl:2005}.
Calculating quandle cohomology, however, is difficult
even in low degrees, mainly for two reasons:
\begin{itemize}
\item
  Brute force calculations are very limited in range.
  Even when they are feasible for small quandles and small degrees,
  their results are usually difficult to interpret.
\item
  Unlike group cohomology, the topological underpinnings
  are less well developed.  Geometric methods that make 
  group theory so rich are mostly absent for quandles.
\end{itemize}

For example, given a diagram of a knot $K \subset \R^3$,
it is comparatively easy to read off a homology class 
$[K] \in H_2(Q_K)$ and to verify that it is an invariant 
of the knot \cite{CarterEtAl:2001}.  Ever since the conception 
of quandle homology, however, it was an important open 
question how to interpret this fundamental class $[K]$, 
and to determine when it vanishes.

The notion of \emph{quandle covering} \cite{Eisermann:2003}
was introduced in order to geometrically interpret 
and finally determine the second (co)homology groups 
$H_2(Q_K) \cong H^2(Q_K) \cong \Z$ for every non-trivial knot $K$. 
More precisely, $H^2(Q_K)$ is freely generated by the canonical class $[E] \in H^2(Q_K)$,
corresponding to the galois covering $E \colon \Z \curvearrowright Q_L \onto Q_K$
coming from the long knot $L$ obtained by cutting $K$ open,
while its dual $H_2(Q_K)$ is freely generated by the fundamental class $[K] \in H_2(Q_K)$.
In particular, $[K]$ vanishes if and only if the knot $K$ is trivial, 
answering Question 7.3 of \cite{CarterEtAl:2001}.
As another consequence, $[K]$ encodes the orientation of the knot $K$,
and so the pair $(Q_K,[K])$ classifies oriented knots.
(The generalization to links with several components will be established 
in \secref{sub:LinkQuandles1} and \secref{sub:LinkQuandles2} below.)

\subsection{Quandle coverings} \label{sub:IntroCoverings}

Knot quandles are somewhat special, and so it was not immediately realized
that covering techniques could be useful for arbitrary quandles as well.
The aim of the present article is to fully develop the algebraic 
covering theory of quandles.  
This will lead us to the appropriate definition of the algebraic 
fundamental group $\pi_1(Q,q)$, and to the Galois correspondence
between connected coverings and subgroups of $\pi_1(Q,q)$.%
\footnote{In a more general context it will be cautious 
  to use the notation $\pi_1^\mathrm{alg}(Q,q)$ to emphasize
  that we are dealing with purely algebraic notions derived from the quandle
  structure $(Q,\ast)$; we do not consider $Q$ as a topological space.
  When $Q$ also carries a topology, $\pi_1^\mathrm{alg}(Q,q)$ 
  should not be confused with the usual topological fundamental 
  group $\pi_1^\mathrm{top}(Q,q)$. 
  While in the present article there seems to be no danger 
  of confusion, the more distinctive notation will become 
  mandatory whenever both concepts are used alongside.}
Detailed definitions and results will be given in 
the next sections, following this overview.

\begin{definition}[see \secref{sub:QuandleCoverings}]
  A quandle homomorphism $p \colon \tilde{Q} \to Q$ is called
  a \emph{covering} if it is surjective and $p(\tilde{y}) = p(\tilde{z})$ 
  implies $\tilde{x} \ast \tilde{y} = \tilde{x} \ast \tilde{z}$ 
  for all $\tilde{x},\tilde{y},\tilde{z} \in \tilde{Q}$.
  In the words of Joyce, $\tilde{y}$ and $\tilde{z}$
  are \emph{behaviourally equivalent}, that is,
  they act in the same way on $\tilde{Q}$.
\end{definition}

\begin{example}
  Consider a group extension $p \colon \tilde{G} \onto G$ 
  and let $\tilde{Q} \subset \tilde{G}$ be a conjugacy class, 
  or more generally a union of conjugacy classes in $\tilde{G}$.
  Without loss of generality we can assume that $\tilde{Q}$ generates $\tilde{G}$.
  As noted above, $\tilde{Q}$ is a quandle with respect to conjugation, 
  and the same holds for its image $Q = p(\tilde{Q}) \subset G$.
  The projection $p \colon \tilde{Q} \onto Q$ is a quandle covering
  if and only if $p$ is a central extension.
\end{example}

As a consequence, the covering theory of \emph{quandles embedded 
in groups} is essentially the theory of central group extensions.
Most quandles, however, do not embed into groups, which is why
quandle coverings have their own distinctive features.
We will see below that unlike central extensions,
the theory of quandle coverings is inherently non-abelian.

\begin{example} \label{exm:DoubleCycle}
  Consider the cyclic group $\Z_m = \Z/m\Z$ with $m \in \N$.
  We explicitly allow $m=0$, in which case $\Z_0 = \Z$.
  The disjoint union $Q_{m,n} = \Z_m \sqcup \Z_n$ becomes 
  a quandle with $a \ast b = a$ for $a,b \in \Z_m$ 
  or $a,b \in \Z_n$, and $a \ast b = a+1$ otherwise.
  This quandle has two connected components, $\Z_m$ and $\Z_n$, 
  each is trivial as a quandle, but both act non-trivially on each other.  
  This expository example will serve us for various illustrations;
  for example, we will see in Proposition \ref{prop:HeisenbergGroup}
  that $Q_{m,n}$ embeds into a group if and only if $m=n$.

  \begin{figure}[hbtp]
    \centering
    \ifpdf\input{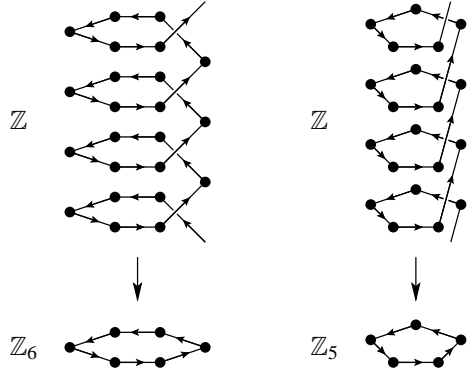}\else\input{qcircles.pstex_t}\fi
    \caption{The universal covering of the quandle $Q_{6,5}$}
    \label{fig:QuandleCovering}
  \end{figure} 

  For every factorization $m = m' m''$ and $n = n' n''$,
  the canonical projections $\Z_m \onto \Z_{m'}$ and $\Z_n \onto \Z_{n'}$
  define a map $p \colon Q_{m,n} \onto Q_{m',n'}$, which is 
  a quandle covering according to our definition.
  (See Figure \ref{fig:QuandleCovering}.)
  In this family, the trivial quandle
  $Q_{1,1} = \{0\} \sqcup \{0\}$ is the terminal object, 
  while $Q_{0,0} = \Z \sqcup \Z$ is the initial object.
  In fact, the map $Q_{0,0} \onto Q_{m,n}$ will turn out 
  to be the universal covering of $Q_{m,n}$, 
  provided that $\gcd(m,n)=1$.
\end{example}

The structure of a quandle $Q$, and in particular its coverings,
are controlled by its adjoint group $\Adj(Q)$, a notion 
introduced by Joyce \cite[\textsection 6]{Joyce:1982}
and recalled in \secref{sub:AdjointGroup} below:

\begin{definition} 
  The \emph{adjoint group} of a quandle $Q$ is the abstract group 
  generated by the elements of $Q$ subject to the 
  relations $a \ast b = b^{-1} a b$ for $a,b \in Q$.
  It comes with a natural map $\adj \colon Q \to \Adj(Q)$ 
  sending each quandle element to the corresponding group element.
  
  There exists a unique group homomorphism $\varepsilon \colon \Adj(Q) \to \Z$ 
  with $\varepsilon(\adj(Q)) = 1$.  We denote its kernel by
  $\Adj(Q)^\circ = \ker(\varepsilon)$.  If $Q$ is connected,
  then $\varepsilon$ is the abelianization of $\Adj(Q)$,
  and $\Adj(Q)^\circ$ is its commutator subgroup.
  Notice that we can reconstruct the adjoint group from $\Adj(Q)^\circ$ 
  as a semi-direct product $\Adj(Q) = \Adj(Q)^\circ \rtimes \Z$.
\end{definition}

\begin{remark}
  Even though it is easily stated, the definition of the 
  adjoint group $\Adj(Q)$ by generators and relations 
  is difficult to work with in explicit calculations.
  Little is known about such groups in general, 
  and only a few examples have been worked out. 
  Since $\Adj(Q)$ turns out to play a crucial r\^ole in determining 
  the coverings of $Q$ and the homology group $H_2(Q)$, we will 
  investigate these groups and calculate several examples 
  in \cite{Eisermann:Wirtinger}.
\end{remark}  

\begin{example} \label{exm:QmnAdjoint}
  For the quandle $Q_{m,n}$ of the previous example we will determine 
  $\Adj(Q_{m,n})$ in Proposition \ref{prop:HeisenbergGroup} below:  
  assuming $\gcd(m,n)=1$ we find $\Adj(Q_{m,n}) = \Z \times \Z$ with 
  $\adj(a) = (1,0)$ for all $a \in \Z_m$ and $\adj(b) = (0,1)$ for all $b \in \Z_n$.
  For $m=n=0$, however, $\Adj(Q_{0,0})$ is the Heisenberg group 
  $H \subset \SL_3\Z$ of upper triangular matrices.  
  Since $Q_{0,0} \onto Q_{m,n} \onto Q_{1,1}$ induces group 
  homomorphisms $H \onto \Adj(Q_{m,n}) \onto \Z \times \Z$,
  we find that $\Adj(Q_{m,n})$ is some intermediate group.
  This turns out to be $H/\gen{z^\ell}$ where $z \in H$
  generates the centre of $H$, and $\ell = \gcd(m,n)$.
\end{example}

\subsection{Galois theory for connected quandles} \label{sub:IntroGaloisConnected}

Motivated by the analogy with topological spaces,
we shall develop the covering theory of quandles 
along the usual lines:
\begin{itemize}
\item Identify the universal covering space (uniqueness, existence, explicit description).
\item Define the fundamental group $\pi_1(Q,q)$ as the group of deck transformations.
\item Establish the Galois correspondence between coverings and subgroups.
\end{itemize}

The results are most easily stated for connected quandles.
They can be suitably refined and adapted to non-connected 
quandles, as explained below and detailed in 
\secref{sec:NonConnectedBase}--\ref{sec:ExtensionCohomology}.

\begin{definition}[see \secref{sub:FundamentalGroup}]
  For a quandle $Q$ we define its \emph{fundamental group} based at $q \in Q$ 
  to be $\pi_1(Q,q) = \{ g \in \Adj(Q)^\circ \mid q^g = q \}$.
\end{definition}

Notice the judicious choice of the group $\Adj(Q)^\circ$;
the approach would not work with another group such as $\Adj(Q)$
or $\Aut(Q)$ or $\Inn(Q)$.  The right choice is not obvious,
but follows from the explicit construction of the universal
covering quandle in \secref{sub:UniversalCovering}.

\begin{proposition}[functoriality, see \secref{sub:FundamentalGroup}]
  Every quandle homomorphism $f \colon (Q,q) \to (Q',q')$ induces 
  a group homomorphism $f_* \colon \pi_1(Q,q) \to \pi_1(Q',q')$.
  We thus obtain a functor $\pi_1 \colon \Quandles_* \to \Groups$ 
  from the category of pointed quandles to the category of groups.
\end{proposition}

\begin{proposition}[lifting criterion, see \secref{sub:LiftingCriterion}]
  Let $p \colon (\tilde{Q},\tilde{q}) \to (Q,q)$ be a quandle covering
  and let $f \colon (X,x) \to (Q,q)$ be a quandle homomorphism
  from a connected quandle $X$.  Then there exists a lifting 
  $\tilde{f} \colon (X,x) \to (\tilde{Q},\tilde{q})$,
  $p \circ \tilde{f} = f$, if and only if 
  $f_* \pi_1(X,x) \subset p_* \pi_1(\tilde{Q},\tilde{q})$. 
  In this case the lifting $\tilde{f}$ is unique.
\end{proposition}


\begin{theorem}[Galois correspondence for connected coverings, 
  see \secref{sub:ConnectedGaloisCorrespondence}]
  For every connected quandle $(Q,q)$ there exists a natural 
  equivalence $\Coverings_\connected(Q,q) \cong \Subgroups(\pi_1(Q,q))$ 
  between the category of pointed connected coverings of $(Q,q)$ 
  and the category of subgroups of $\pi_1(Q,q)$.
  Moreover, a normal subgroup $K \subset \pi_1(Q,q)$ corresponds 
  to a galois covering $p \colon (\tilde{Q},\tilde{q}) \to (Q,q)$
  with deck transformation group $\Aut(p) \cong \pi_1(Q,q)/K$.  
  \qed
\end{theorem}


The Galois correspondence can be extended to non-connected coverings,
and further to principal $\Lambda$-coverings.  The latter correspond 
to extensions $\Lambda \curvearrowright \tilde{Q} \onto Q$
of the quandle $Q$ by some group $\Lambda$ as defined in
\secref{sub:QuandleExtensions}.

\begin{theorem}[Galois correspondence for general coverings,
  see \secref{sub:SemiconnectedGaloisCorrespondence}]
  For every connected quandle $(Q,q)$ there exists a natural
  equivalence $\Coverings(Q) \cong \Actions(\pi_1(Q,q))$ 
  between the category of coverings of $(Q,q)$ 
  and the category of actions of $\pi_1(Q,q)$.
  Moreover, there exists a natural bijection 
  $\Ext(Q,\Lambda) \cong \Hom(\pi_1(Q,q),\Lambda)$
  between equivalence classes of extensions 
  $\Lambda \curvearrowright \tilde{Q} \onto Q$
  and the set of group homomorphisms $\pi_1(Q,q) \to \Lambda$.
\end{theorem}

Throughout this article our guiding principle is the analogy 
between the covering theories of topological spaces and quandles.
While their overall structure is the same, the individual 
objects seem quite different.  The formal analogy may thus
come as a surprise, even more so as it pervades even 
the tiniest details.  This can in large parts be explained 
by the common feature of the fundamental groupoid,
as described in \secref{sec:FundamentalGroupoid}.
We will complete this analogy in \secref{sec:ExtensionCohomology}
by establishing the relationship with (co)homology:

\begin{theorem}[Hurewicz isomorphism for connected quandles,
  see \secref{sub:HurewiczIsomorphism}]
  For every connected quandle $Q$ we have a natural 
  isomorphism $H_2(Q) \cong \pi_1(Q,q)_\mathrm{ab}$.  
  Moreover, for every group $\Lambda$ we have natural bijections 
  $H^2(Q,\Lambda) \cong \Ext(Q,\Lambda) \cong \Hom(\pi_1(Q,q),\Lambda)$.
  If $\Lambda$ is an abelian group, or more generally a module over some ring $R$, 
  then these objects carry natural $R$-module structures and 
  the natural bijections are isomorphisms of $R$-modules.
\end{theorem}

The introduction of a cohomology $H^2(Q,\Lambda)$ 
with non-abelian coefficients $\Lambda$ is natural 
inasmuch as it allows us to treat all cases in a uniform way.  
This is analogous to the cohomology $H^1(X,\Lambda)$ of 
a topological space $X$ with non-abelian coefficients $\Lambda$,
see \cite{Olum:1958}.

\subsection{Examples and applications} \label{sub:IntroExamples}

As a general application, let us mention that
every quandle $Q$ can be obtained as a covering
of a quandle $\bar{Q} \subset G$ in some group $G$.
(Take for example the image of $Q$ in its inner automorphism group.)
This can be useful, for example, in understanding 
finite connected quandles: it suffices to consider 
conjugacy classes $\bar{Q}$ in finite groups $G$ such that $G = \gen{\bar{Q}}$,
together with their covering quandles $Q \onto \bar{Q}$;
these are parametrized by subgroups of the fundamental group $\pi_1(\bar{Q},\bar{q})$.

\begin{remark}
  For every finite connected quandle $Q$ the group $\Adj(Q)^\circ$
  is finite, whence the fundamental group $\pi_1(Q,q)$ and 
  the universal covering $\tilde{Q} \to Q$ are both finite,
  cf.\ \cite{Eisermann:Wirtinger}.
\end{remark}

\begin{example}[dihedral quandles] \label{exm:DihedralFundamentalGroup}
  The dihedral quandle $D_n$ is obtained from the cyclic group 
  $\Z_n = \Z/n\Z$ with the quandle operation $a \ast b = 2b-a$.  
  It is isomorphic to the subquandle $\Z_n \times \{1\}$ of 
  the dihedral group $\Z_n \rtimes \Z_2$, corresponding to 
  the $n$ reflections of a regular $n$-gon.  
  For $n$ odd, the quandle $D_n$ is connected, and 
  we find $\Adj(D_n) = \Z_n \rtimes \Z$ with group action 
  $(a,i) \cdot (b,j) = ( a + (-1)^i b, i+j )$, 
  and $\adj \colon D_n \to \Adj(D_n)$ is given 
  by $\adj(a) = (a,1)$, cf.\ \cite{Eisermann:Wirtinger}.  
  Since $\Adj(Q)^\circ = \Z_n \rtimes \{0\}$ 
  acts on $D_n$ by $a^{(b,0)} = a-2b$, we find 
  the fundamental group $\pi_1(D_n,0) = \{0\}$.  
  This means that evey dihedral quandle $D_n$ 
  of odd order is \emph{simply connected}.
  Equivalently, every quandle covering of $D_n$ is trivial,
  that is, equivalent to $\pr_1 \colon D_n \times F \onto D_n$,
  where $F$ is some trivial quandle.
\end{example}

\begin{Notes}
  For $n$ even, $D_n$ has two components,
  formed by even and odd elements, respectively,
  and the situation is more involved, cf.\ \cite{Eisermann:Wirtinger}.
  (Notice that $D_4 \cong Q_{2,2}$ is discussed above.)
  The adjoint group $\Adj(D_n)$ can be explicitly 
  identified as follows: $D_n$ embeds into the group 
  \[
  G = \Z_n \rtimes \Z^2 \quad\text{with}\quad
  (a,i) \cdot (b,j) = (a+(-1)^{i_1+i_2}b,i+j), 
  \]
  via the map $\phi \colon D_n \to G$ given by $a \mapsto (a,1,0)$ 
  if $a$ is even, and $a \mapsto (a,0,1)$ if $a$ is odd.
  The generated subgroup $H = \gen{\phi(D_n)}$ consists
  of those triples $(a,i_1,i_2)$ with $a+i_1$ odd and $a+i_2$ even.
  The explicit representation $\phi \colon D_n \to H$ 
  so obtained is equivalent to the adjoint representation 
  $\adj \colon D_n \to \Adj(D_n)$.  For details see \cite{Eisermann:Wirtinger}.
\end{Notes}

\begin{example}[symmetric groups] \label{exm:SymmetricGroups}
  Consider the symmetric group $S_n$ on $n \ge 3$ points, and 
  let $Q$ be the conjugacy class of the transposition $q = (12)$.
  This is a quandle with $\binom{n}{2} = \frac{n(n-1)}{2}$ elements.
  It is not difficult to see that $\Adj(Q) = A_n \rtimes \Z$, 
  where the action of $k \in \Z$ on $A_n$ is given by 
  $a \mapsto (12)^k a (12)^k$, cf.\ \cite{Eisermann:Wirtinger}.
  We thus find $\Adj(Q)^\circ = A_n$, which yields 
  the fundamental group $\pi_1(Q,q) \cong S_{n-2}$.
  The subgroups of $S_{n-2}$ thus characterize 
  the connected coverings of the quandle $Q$.  
  (For $n=3$ notice that $Q=D_3$, for which 
  we already know that $\pi_1$ is trivial;
  $\pi_1(Q,q)$ is non-trivial only for $n \ge 4$.)

  Turning to the extensions of $Q$ by some group $\Lambda$,
  we find $H^2(Q,\Lambda) \cong \Ext(Q,\Lambda) \cong \Hom(S_{n-2},\Lambda)$.
  If $\Lambda$ is abelian, we see without any further calculation that $H^2(Q,\Lambda)$ 
  is trivial for $n=3$, and isomorphic to the group of $2$-torsion elements in
  $\Lambda$ for $n\ge4$, because $(S_{n-2})_\mathrm{ab} \cong \Z_2$.
  Moreover, $H_2(Q) = 0$ for $n = 3$, and $H_2(Q) = \Z_2$ for $n \ge 4$.
\end{example}

\begin{example}[knot quandles]
  As in \cite[\textsection3]{Eisermann:2003} 
  let $L$ be a long knot and let $K$ be its corresponding closed knot.
  Both knot quandles $Q_L$ and $Q_K$ are connected, their adjoint groups 
  are $\Adj(Q_L) = \Adj(Q_K) = \pi_K$, and the natural projection 
  $p \colon Q_L \to Q_K$ is a quandle covering.  We may choose a canonical 
  base point $q_L \in Q_L$ and its image $q_K \in Q_K$.  Both map to a meridian 
  $m_L = m_K \in \pi_K$, and we denote by $\ell_K \in \pi_K$ the corresponding longitude.
  The explicit construction of universal coverings in \cite{Eisermann:2003} 
  shows $\pi_1(Q_L,q_L) = \{1\}$, and so $Q_L$ is the universal covering 
  of the quandle $Q_K$.  For the quotient $Q_K = \gen{\ell_K} \backslash Q_L$ 
  we thus find $\pi_1(Q_K,q_K) = \gen{\ell_K}$, 
  whence $\pi_1(Q_K,q_K) \cong \Z$ for every non-trivial knot $K$.  


  This observation, although not in the language of quandle coverings and fundamental groups,
  was used by Joyce \cite{Joyce:1982} in order to recover the knot group data 
  $(\pi_K,m_K,\ell_K^\pm)$ from the knot quandle $Q_K$.  According to 
  Waldhausen's result \cite{Waldhausen:1968}, the triple $(\pi_K,m_K,\ell_K)$
  classifies knots, so the knot quandle classifies knots modulo inversion.
  The remaining ambiguity can be removed by the orientation class
  $[K]\in H_2(Q_K)$, as explained in \cite[\textsection6]{Eisermann:2003}.
\end{example}

\begin{remark}[knot colouring polynomials]
  The knot quandle $Q_K$, just as the knot group $\pi_K$,
  is in general very difficult to analyze.  A standard way
  to extract information is to consider (finite) representations:
  we fix a finite quandle $Q$ with base point $q \in Q$ 
  and consider knot quandle homomorphisms 
  $\phi \colon (Q_K,q_K) \to (Q,q)$.  Each $\phi$ induces a group homomorphism
  $\phi_* \colon \pi_1(Q_K,q_K) \to \pi_1(Q,q)$, which is determined
  by the image of the canonical generator $\ell_K \in \pi_1(Q_K,q_K)$.
  We can thus define a map 
  \[
  P_Q^q \colon \{\mathrm{knots}\} \to \Z\pi_1(Q,q) \quad by \quad
  P_Q^q(K) := \sum_{\phi \colon (Q_K,q_K) \to (Q,q)} \phi_*(\ell_K) .
  \]
  This invariant is the \emph{knot colouring polynomial} associated to $(Q,q)$,
  and provides a common generalization to the invariants presented in 
  \cite{Eisermann:ColoPoly} and \cite{Niebrzydowski:2006}.  
  Colouring polynomials encode, in particular, all quandle $2$-cocycle 
  invariants, as proven in \cite{Eisermann:ColoPoly}.

  Example \ref{exm:DihedralFundamentalGroup} above shows that
  the longitude images are necessarily trivial for dihedral colourings;
  the only information extracted is the number of $n$-colourings.
  The situation is different for $Q = (12)^{S_n}$, where longitude
  images yield more refined information.
\end{remark}

\begin{example}
  We conclude with another natural and highly non-abelian 
  example, where our tools are particularly efficient.
  Consider the quandle $Q_K^\pi \subset \pi_K$ consisting of all meridians of 
  the knot $K$, that is, the conjugacy class of our preferred meridian $m_K$
  in $\pi_K$, or equivalently, the image of the natural quandle homomorphism $Q_K \to \pi_K$.
  Here we find $\Adj(Q_K^\pi) = \pi_K$, and $\pi_1(Q_K^\pi,m_K)$ is a free
  group of rank $n$ if $K = K_1 \cs \cdots \cs K_n$ is the connected
  sum of $n$ prime knots \cite[Corollary 39]{Eisermann:2003}.
  Via the Hurewicz isomorphism we obtain that $H_2(Q_K^\pi) \cong \Z^n$,
  as previously noted in \cite[Theorem 53]{Eisermann:2003}.
\end{example}

\subsection{Tournants dangereux} \label{sub:TournantsDangereux}

There are a number of subtleties where quandle coverings
do not behave as could be expected at first sight.  
First of all, they do not form a category:

\begin{example} \label{exm:DoubleDoubleCovering}
  The abelian group $Q = \Q/\Z$ becomes a connected quandle with $a \ast b = 2b-a$.
  The map $p \colon Q \to Q$, $a \mapsto 2a$, is a quandle covering.
  The composition $p \circ p \colon Q \to Q$, $a \mapsto 4a$, however, 
  is \emph{not} a covering: $0$ and $\frac{1}{4}$ do not act in the same way on $Q$.
  The same phenomenon already appears for finite quandles,
  for example $D_{4n} \lonto[_2] D_{2n} \lonto[_2] D_n$.
\end{example}

\begin{remark}
  Coverings of topological spaces suffer from the same problem, 
  see Spanier \cite{Spanier:1981}, Example 2.2.8:
  given two coverings $p \colon X \to Y$ and $q \colon Y \to Z$, 
  their composition $qp \colon X \to Z$ is not necessarily a covering.
  This phenomenon is, however, rather a pathology: the composition $qp$ is always 
  a covering if $Z$ is locally path connected and semilocally $1$-connected 
  (see \cite{Spanier:1981}, Theorems 2.2.3, 2.2.6, 2.4.10).
  These hypotheses hold, in particular, for coverings 
  of manifolds, simplicial complexes, or CW-complexes.

  When we speak of topological covering theory as our model, 
  we will neglect all topological subtleties such as 
  questions of local and semilocal connectedness. 
  The reader should think of covering theory 
  in its nicest possible form, say for CW-complexes.
\end{remark}

\begin{remark}
  There are two further aspects in which quandle coverings 
  differ significantly from the model of topological coverings:
  \begin{itemize}
  \item
    For a quandle covering $p \colon (\tilde{Q},\tilde{q}) \to (Q,q)$
    the induced map on the fundamental groups, 
    $p_* \colon \pi_1(\tilde{Q},\tilde{q}) \to \pi_1(Q,q)$,
    need not be injective.
  \item
    If $\tilde{Q}$ is simply connected, then $p$ is the universal covering of $(Q,q)$.
    The converse is not true: it may well be that $p$ is universal
    but $\tilde{Q}$ is not simply connected.
  \end{itemize}
  
  It is amusing to note that the Galois correspondence stated above 
  is salvaged because these two defects cancel each other.
\end{remark}

\begin{example}
  For $Q = \Q/\Z$ one finds $\Adj(Q) = (\Q/\Z) \rtimes \Z$ 
  with $\adj(a) = (a,1)$, cf.\ \cite{Eisermann:Wirtinger}.  
  The subgroup $\Adj(Q)^\circ = \Q/\Z \times \{0\}$ acts 
  on $Q$ via $a^{(b,0)} = a - 2b$, which implies that 
  $\pi_1(Q,0) = \{\, (0,0),\, (\frac{1}{2},0) \,\} \cong \Z_2$.
  This means that $p \colon Q \to Q$, $a \mapsto 2a$
  is the universal covering.  In particular, the universal covering quandle 
  is \emph{not} simply connected, and the induced homomorphism $p_*$ 
  between fundamental groups is \emph{not} injective.
\end{example}

\begin{Notes}

  
  \begin{proof}
    Obviously $Q$ embeds into the group $G = (\Q/\Z) \rtimes \Z$ 
    via $\phi \colon Q \to G$, $a \mapsto (a,1)$: on $G$ the group operation 
    $(a,i) \cdot (b,j) = ( a + (-1)^i b, i+j )$ induces on $Q \times \{1\} \subset G$ 
    the desired conjugation action $(a,1) \ast (b,1) =  (2b-a,1)$.
    By the universal property of $\adj_Q \colon Q \to \Adj(Q)$
    there exists a unique group homomorphism $h \colon \Adj(Q) \to G$
    such that $\phi = h \circ \adj_Q$.  We have $\gen{\phi(Q)} = G$, 
    and so $h$ is surjective.  The point is to prove that $h$ is injective. 
    
    Every finite subquandle of $Q$ is contained in 
    $D_n = (\frac{1}{n}\Z)/\Z$ for some $n \in \N$,  
    which is the dihedral quandle of order $n$.
    Every finitely generated subgroup of $\Adj(Q)$
    is thus contained in some subgroup $H_n = \gen{\adj_Q(D_n)}$.
    It suffices to show that $h$ is injective on each $H_n$.
    
    Notice that $D_n$ is one of two connected components in $D_{2n}$,
    and we have $\Adj(D_{2n}) \into (\frac{1}{2n}\Z/\Z) \rtimes \Z^2$.
    We obtain the following commutative diagram:
    \[
    \begin{diagram}
      \node{D_n}
      \arrow{e,J}
      \arrow{se,b}{f}
      \node{D_{2n}}
      \arrow{s,r}{\adj_{D_{2n}}}
      \arrow{e,J}
      \node{Q}
      \arrow{s,r}{\adj_Q}
      \\
      \node{}
      \node{\Adj(D_{2n})}
      \arrow{e,t}{k}
      \arrow{se,b}{g}
      \node{\Adj(Q)}
      \arrow{s,t}{h}
      \\
      \node{}
      \node{}
      \node{G}
    \end{diagram}
    \]
    The group homomorphism $g$ maps $(a,i_1,i_2)$ to $(a,i_1+i_2)$, and is thus 
    injective on the subgroup $K_n = \gen{f(D_n)} \subset \Adj(D_{2n})$.
    This ensures that $h$ is injective on $H_n = k(K_n) \subset \Adj(Q)$.
  \end{proof}
  
\end{Notes}

\begin{remark}
  The previous example may appear somewhat artificial, 
  because the problem essentially arises from $2$-torsion and 
  the fact that all $2$-torsion elements are $2$-divisible.
  In particular, these conditions force $Q$ to be infinite.
  Example \ref{exm:LinearGroups3} exhibits a \emph{finite} quandle 
  with a universal covering that is not simply connected.
  This is definitely not a pathological construction: 
  the phenomenon naturally occurs for certain conjugacy 
  classes in groups,  for example the conjugacy class of
  $\left[\begin{smallmatrix} 0 & 1 \\ -1 & 0 \end{smallmatrix}\right]$
  in the group $\PSL_2\K$ over a finite field $\K$.
\end{remark}

\subsection{Perfect groups} \label{sub:PerfectCoverings}

Quandle coverings resemble Kervaire's algebraic covering
theory of perfect groups \cite{Kervaire:1970}, which he applied
to algebraic K-theory in order to identify the Milnor group $K_2(A)$ 
of a ring $A$ with the Schur multiplier $H_2(\GL(A)')$.
It is illuminating to contrast the theory of quandle coverings
with Kervaire's classical results.

Recall that a group $G$ is \emph{perfect}, 
or \emph{connected} in the words of Kervaire,
if $G'=[G,G]=G$, or equivalently $H_1(G) = G_{\mathrm{ab}} = 0$.
A \emph{covering} of $G$ is a central extension $\tilde{G} \onto G$
with $\tilde{G}$ perfect.  Kervaire established a bijection between
subgroups of $H_2(G)$ and isomorphism classes of coverings $\tilde{G} \onto G$.
The theory is thus analogous to the covering theory of topological spaces, 
and consequently Kervaire defined $\pi_1(G) := H_2(G)$.

\begin{remark}
  By construction, $\pi_1(G)$ is abelian and base points play no r\^ole.
  Moreover, the covering theory of perfect groups is well-behaved 
  in the following sense:
  \begin{itemize}
  \item
    Coverings of perfect groups form a category, which means that the composition 
    of two coverings is again a covering \cite[Lemme 1]{Kervaire:1970}.
  \item
    A covering $\tilde{G} \onto G$ is universal 
    if and only if $\tilde{G}$ is simply connected, 
    that is, $\pi_1(\tilde{G}) = H_2(\tilde{G}) = 0$ \cite[Lemme 2]{Kervaire:1970}.
  \item
    For every covering $p \colon \tilde{G} \to G$ the induced map 
    $p_* \colon \pi_1(\tilde{G}) \to \pi_1(G)$ is injective
    \cite[Th\'eor\`eme de classfication]{Kervaire:1970}.
  \end{itemize}
  As we have seen above, quandle coverings do not enjoy 
  these privileges in general.  They may thus be considered 
  a ``non-standard'' covering theory that warrants a careful analysis.
\end{remark}

The analogy between coverings of quandles and perfect groups 
is not only a formal one.  Their precise relationship will 
be studied in a forthcoming article \cite{Eisermann:Wirtinger}.  
As an illustration, it can be applied to determine certain adjoint groups:

\begin{theorem}
  Let $G$ be a simply connected group, i.e.\ $H_1(G) = H_2(G)=0$.
  Consider a conjugacy class $Q = q^G$ that generates $G$,
  so that $Q$ is a connected quandle.  Then we have
  an isomorphism $\Adj(Q) \isoto G \times \Z$ given by
  $\adj(q) \mapsto (q,1)$ for all $q \in Q$.
  In particular, we obtain $\Adj(Q)' = G$
  and $\pi_1(Q,q) = C_G(q) = \{ g \in G \mid q^g = q \}$.
  \qed
\end{theorem}

This directly applies to every simple group $G$ with 
Schur multiplier $H_2(G)=0$.  Most often we have $H_2(G)\ne0$,
in which case it suffice to pass to the universal covering $\tilde{G}$.

\subsection{Generalization to non-connected quandles} \label{sub:IntroGaloisNonConnected}

One final difficulty arises when we pass from connected to non-connected quandles.
In the analogous model of topological spaces, this generalization is simple, 
because a topological space (say locally connected) is the disjoint union of its components.
For quandles, however, this is far from being true: the different components
act on each other, and this interaction is in general non-trivial.
In particular, the disjoint union is not the appropriate model.

In order to develop a covering theory for non-connected quandles
we have to treat all components  individually yet simultaneously.  
The convenient way to do this is to index the components by some fixed set $I$,
and then to deal with $I$-graded objects throughout.
(For details see Section \ref{sec:NonConnectedBase}.)
The upshot is that for a non-connected quandle $Q$
all preceding statements remain true when 
suitably interpreted in the graded sense:

\begin{definition}[grading, see \secref{sub:GradedQuandles}]
  A \emph{graded quandle} is a quandle $Q = \bigsqcup_{i \in I} Q_i$
  partitioned into subsets $(Q_i)_{i \in I}$ such that 
  $Q_i \ast Q_j = Q_i$ for all $i,j \in I$.
  A \emph{pointed quandle} $(Q,q)$ is a graded quandle 
  with a base point $q_i \in Q_i$ for each $i \in I$.
  We call $(Q,q)$ \emph{well-pointed} if $q$ specifies one base point 
  in each component, i.e.\ $Q_i$ is the component of $q_i$ in $Q$.
  In this case we define the \emph{graded fundamental group} 
  to be the product $\pi_1(Q,q) := \prod_{i \in I} \pi_1(Q,q_i)$. 
\end{definition}

\begin{theorem}[Galois correspondence, see \secref{sec:GradedGaloisCorrespondence}]
  Let $(Q,q)$ be a well-pointed quandle indexed by some set $I$.
  There exists a natural equivalence 
  $\Coverings_I(Q,q) \cong \Subgroups_I(\pi_1(Q,q))$ 
  between the category of well-pointed coverings of $(Q,q)$ 
  and the category of graded subgroups of $\pi_1(Q,q)$.
  Likewise, there exists a natural equivalence
  $\Coverings(Q) \cong \Actions(\pi_1(Q,q))$
  between the category of coverings of $(Q,q)$ 
  and the category of graded actions of $\pi_1(Q,q)$.
\end{theorem}

\begin{theorem}[Hurewicz isomorphism for general quandles,
  see \secref{sub:HurewiczIsomorphism}]
  For every well-pointed quandle $(Q,q)$ we have a natural isomorphism 
  $H_2(Q) \cong \bigoplus_{i \in I} \pi_1(Q,q_i)_\mathrm{ab}$, and
  for every graded group $\Lambda$ we have natural bijections 
  \[ 
  H^2(Q,\Lambda) \cong \Ext(Q,\Lambda) \cong \Hom(\pi_1(Q,q),\Lambda)
  = \prod_{i \in I} \Hom(\pi_1(Q,q_i),\Lambda_i) .
  \]
\end{theorem}

\begin{Notes}
\begin{example}
  Reconsider the quandle $Q = Q_{m,n}$ of Example \ref{exm:DoubleCycle}. 
  Assuming $\gcd(m,n)=1$ we have $\Adj(Q) = \Z^2$ and thus 
  $\Adj(Q)^\circ = \{ (a,-a) \mid a \in \Z \} \cong \Z$.
  For $q_1 \in \Z_m$ we find $\pi_1(Q,q_1) =  \{ (ma,-ma) \mid a \in \Z \}$,
  and $\pi_1(Q,q_2) =  \{ (-na,na) \mid a \in \Z \}$ for $q_2 \in \Z_n$.  
  The graded fundamental group is thus $\pi_1(Q,q) \cong m\Z \times n\Z$,
  and the universal covering is given by the projection 
  $Q_{0,0} \onto Q_{m,n}$ depicted above.  
  Finally we have $H_2(Q) \cong m\Z \oplus n\Z$.
  (For general parameters $m,n$ the situation is slightly
  more complicated, see Example \ref{exm:QmnCovering}.)
\end{example}
\end{Notes}

One of the motivations to study non-connected quandles
is their application to links.  Given an $n$-component 
link $K = K_1 \sqcup \dots \sqcup K_n \subset \S^3$, 
we choose a base point $q_K^i \in Q_K$ for each 
link component $K_i$, and obtain a decomposition
$Q_K = Q_K^1 \sqcup \dots \sqcup Q_K^n$ into 
components $Q_K^i = [q_K^i]$.  This establishes
a natural bijection $\pi_0(K) \isoto \pi_0(Q_K)$.
 
\begin{theorem}[see \secref{sub:LinkQuandles1}]
  For every link $K \subset \S^3$ the graded fundamental 
  group of the link quandle $Q_K$ is given by 
  $\pi_1(Q_K,q_K) = \prod_{i=1}^n \gen{\ell_K^i}$,
  where $\ell_K^i \in \Adj(Q_K) = \pi_1(\S^3 \minus K)$
  is the longitude associated to the meridian 
  $m_K^i = \adj(q_K^i) \in \Adj(Q_K)$.
\end{theorem}

We will see in \secref{sub:LinkQuandles2} that the Hurewicz 
isomorphism maps the longitude $\ell_K^i \in \pi_1(Q_K,q_K^i)$
to the orientation class $[K_i] \in H_2(Q_K)$ of the component $K_i$.
This shows that the quandle $Q_K$ is a classifying 
invariant of the link $K$ in the following sense:

\begin{theorem}[see \secref{sub:LinkQuandles2}]
  Two oriented links $K = K_1 \sqcup \dots \sqcup K_n$
  and $K' = K'_1 \sqcup \dots \sqcup K'_n$ in $\S^3$ 
  are ambient isotopic respecting orientations and numbering 
  of components if and only if there exists a quandle 
  isomorphism $\phi \colon Q_K \isoto Q_{K'}$ such that 
  $\phi_* [K_i] = [K'_i]$ for all $i$. 
\end{theorem}


\subsection{Related work} \label{sub:RelatedWork}

The present article focuses on the systematic investigation 
of quandle coverings and their Galois correspondence.  
The explicit construction of a universal covering and 
the definition of the corresponding algebraic 
fundamental group appear here for the first time.
Our construction can easily be adapted to racks: 
here $\Adj(Q)^\circ$ has to be replaced by $\Adj(Q)$,
and the definition of the fundamental group has to be adapted accordingly.
Modulo these changes, our results hold verbatim for racks instead of quandles.

As it could be expected, these notions are 
closely related to quandle extensions and cohomology, 
which have both been intensively studied in recent years.
The subject of rack cohomology originated in the work
of R.\,Fenn, C.\,Rourke, and B.\,Sanderson \cite{FennRourkeSanderson:1995},
who constructed a classifying topological space $BX$ for every rack $X$.
The corresponding quandle (co)homology theory was taken up by J.S.\,Carter 
and his collaborators, in order to construct knot invariants
(see for example \cite{CarterEtAl:2001,CarterEtAl:2003,CarterEtAl:2005}).
Quandle coverings were introduced and applied to knot quandles in \cite{Eisermann:2003}.
They have also appeared in the context of non-abelian extensions, explored 
by N.\,Andruskiewitsch and M.\,Gra\~na \cite{AndruskiewitschGrana:2003},
where a corresponding non-abelian cohomology theory was proposed.
This generalized cohomology, in turn, has been taken up and applied 
to knot invariants in \cite{CarterEtAl:2005}.

We have stated above how our approach of quandle coverings can be applied 
to complete the trilogy of cohomology $H^2(Q,\Lambda)$ and extensions 
$\Ext(Q,\Lambda)$ by the third aspect: the fundamental group $\pi_1(Q,q)$.  
The result is the natural isomorphism
\begin{equation} \label{eq:QuandleTrilogy}
  H^2(Q,\Lambda) \cong \Ext(Q,\Lambda) \cong \Hom( \pi_1(Q,q), \Lambda ).
\end{equation}
A similar isomorphism has been noted by P.\,Etingof 
and M.\,Gra\~na \cite[Cor.\,5.4]{EtingofGrana:2003}: 
for every rack $X$ and every abelian group $A$ they prove that 
$H^2(X,A) \cong H^1( \Adj(X), \Map(X,A) )$, where $\Map(X,A)$ is 
the module of maps $X \to A$ with the action of the adjoint group $\Adj(X)$. 
The formulation \eqref{eq:QuandleTrilogy} takes this
one step further and highlights the geometric meaning.
For practical calculations it is as explicit 
and direct as one could possibly wish.

\subsection{Acknowledgements} \label{sub:Acknowledgments}

The concept of quandle covering, algebraic fundamental group,
and Galois correspondence developed in 2001 
when I was working on the article \cite{Eisermann:2003}. 
At that time, however, I saw no utility of this theory 
other than its application to knot quandles.  
In the intervening years, non-abelian extensions 
have gained interest, and in November 2006
the conference \textsl{Knots in Washington XXIII} 
on ``Quandles, their homology and applications''
convinced me that covering theory would be a welcome complement.
I thank J\'ozef Przytycki and the organizers for bringing together this meeting.

\subsection{How this article is organized}

The article follows the outline given in the introduction.
Section \ref{sec:QuandleDefinitions} reviews the basic definitions
of quandle theory leading up to quandle coverings, while
Section \ref{sec:Examples} displays some detailed examples.
Section \ref{sec:CoveringQuandles} records 
elementary properties of quandle coverings.
Section \ref{sec:ConnectedCoverings} constructs the universal connected covering, 
defines the fundamental group, and establishes the Galois correspondence
for connected coverings.  Section \ref{sec:NonConnectedCoverings} 
explains how to extend these results to non-connected coverings over 
a connected base quandle, while Section \ref{sec:NonConnectedBase}
discusses the technicalities necessary for non-connected base quandles.
Section \ref{sec:FundamentalGroupoid} expounds the concept of
fundamental groupoid in order to explain the striking 
similarity between quandles and topological spaces.
Section \ref{sec:ExtensionCohomology}, finally, elucidates the correspondence
between quandle extensions and quandle cohomology in the non-abelian
and graded setting, and thus completes the trilogy $H^2$, $\Ext$, $\pi_1$.


\section{Definitions and elementary properties} \label{sec:QuandleDefinitions}

The following definitions serve to fix our notation and to make 
the presentation self-contained.  They are mainly taken from Joyce \cite{Joyce:1982},
suitably extended and tailored to our application.  
Some immediate examples are stated alongside the definitions,
more elaborate examples will be postponed until the next section.

We also seize the opportunity to record some elementary 
but useful observations, which have been somewhat 
neglected or dispersed in the published literature.
In particular, we emphasize the r\^ole played by central 
group extensions, which come to light at several places.
While on the level of groups only central extensions are visible,
quandle coverings turn out to be essentially non-abelian 
(see Example \ref{exm:SymmetricGroups} above).

\subsection{The category of quandles}

A \emph{quandle} is a set $Q$ equipped with two binary operations 
$\ast,\tsa \colon Q\times Q\to Q$ satisfying the three axioms
stated in the introduction.  These axioms are symmetric in $\ast$
and $\tsa$: if $(Q,\ast,\tsa)$ is a quandle, then so is $(Q,\tsa,\ast)$.
Moreover, each of the operations $\ast$ and $\tsa$ determines the other,
so we can simply write $(Q,\ast)$ instead of $(Q,\ast,\tsa)$.
If both operations coincide we have $(a \ast b) \ast b = a$ 
for all $a,b \in Q$, which is called an \emph{involutory} quandle.
We will use the same symbol ``$\ast$'' for different quandles,
and we will frequently denote a quandle by $Q$ instead of $(Q,\ast)$,
unless there is danger of confusion.

\begin{definition}
  A \emph{quandle homomorphism} between two quandles $Q$ and $Q'$
  is a map $\phi \colon Q \to Q'$ satisfying 
  $\phi(a \ast b) = \phi(a) \ast \phi(b)$, and hence 
  $\phi(a \tsa b) = \phi(a) \tsa \phi(b)$, for all $a,b\in Q$.
  Quandles and their homomorphisms form a category, denoted $\Quandles$.

  The automorphism group $\Aut(Q)$ 
  consists of all bijective homomorphisms $\phi \colon Q \to Q$.  
  We will adopt the convention that automorphisms of $Q$ act on the right, 
  written $a^\phi$, which means that their composition $\phi\psi$ 
  is defined by $a^{(\phi\psi)} = (a^\phi)^\psi$ for all $a \in Q$.
\end{definition}

\begin{example} \label{exm:ConjugationQuandle}
  Every group $(G,\cdot)$ defines a quandle $(G,\ast)$ with $a \ast b = b^{-1} a b$.
  This is called the \emph{conjugation quandle} of $G$ and denoted $\Conj(G)$.
  Every group homomorphism $(G,\cdot) \to (H,\cdot)$
  is also a quandle homomorphism $(G,\ast) \to (H,\ast)$.
  We thus obtain a functor $\Conj \colon \Groups \to \Quandles$
  from the category of groups to the category of quandles.
\end{example}

\begin{example} \label{exm:CoreQuandle}
  Every group $(G,\cdot)$ defines an involutory
  quandle $(G,\ast)$ with $a \ast b = b a^{-1} b$.
  This is called the \emph{core quandle} of $G$ and denoted $\Core(G)$.
  Every group homomorphism $(G,\cdot) \to (H,\cdot)$
  is also a quandle homomorphism $(G,\ast) \to (H,\ast)$.
  We thus obtain another functor $\Core \colon \Groups \to \Quandles$
  from the category of groups to the category of quandles.
\end{example}

\begin{example} \label{exm:AlexanderQuandle}
  If $A$ is a group and $T \colon A \isoto A$ an automorphism,
  then $A$ becomes a quandle with $a \ast b = T(ab^{-1})b$. 
  This is called the \emph{Alexander quandle} of $(A,T)$, 
  denoted $\Alex(A,T)$.  
  Every group homomorphism $\phi \colon (A,T) \to (B,S)$ 
  with $\phi \circ T = S \circ \phi$ is also a homomorphism
  of Alexander quandles $(A,\ast) \to (B,\ast)$.
  We thus obtain a functor $\Alex \colon \category{GrpAut} \to \Quandles$
  from the category of group automorphisms to the category of quandles.

  If $A$ is abelian, then the pair $(A,T)$ is equivalent 
  to a $\Z[t^\pm]$-module $A$ with $ta = T(a)$ for all $a \in A$.  
  Restricting to this case, we obtain a functor 
  $\Alex \colon \category{Mod}_{\Z[t^\pm]} \to \Quandles$
  from the category of $\Z[t^\pm]$-modules to the category of quandles.
\end{example}

\begin{remark}
  Our definition of Alexander quandles is more inclusive 
  than usual, in order to embrace also non-abelian groups.
  Joyce \cite[\textsection7]{Joyce:1982} used the general construction,
  but reserved the name \emph{Alexander quandle} for abelian groups $A$.  
  In this case the quandle $\Alex(A,T)$ is an \emph{abelian} quandle in the sense that 
  $(a \ast b) \ast (c \ast d)  = (a \ast c) \ast (b \ast d)$ for all $a,b,c,d \in Q$. 
  Notice the special case $\Alex(A,-\id) = \Core(A,+)$.
\end{remark}

\begin{Notes}
  \begin{remark}
    Recall that a group $(G,\cdot)$ is abelian if and only if
    the set $\End(G,\cdot)$ of endomorphisms is a group with respect 
    to pointwise multiplication,  $(f \cdot g)(x) = f(x) \cdot g(x)$.
    Likewise, if a quandle $(Q,\ast)$ is abelian, then
    the set $\End(Q,\ast)$ of endomorphisms is a quandle 
    with respect to the pointwise operation defined by
    $a^{\phi \ast \psi} = a^\phi \ast a^\psi$.
  \end{remark}
\end{Notes}

\subsection{Inner automorphisms}

The quandle axioms (Q2) and (Q3) are equivalent to saying 
that for every $a \in Q$ the right translation 
$\rho_a : x \mapsto x\ast a$ is an automorphism of $Q$.  
Such structures were studied by E.\,Brieskorn \cite{Brieskorn:1988} 
under the name ``automorphic sets'' and by C.\,Rourke and R.\,Fenn
\cite{FennRourke:1992} under the name ``rack''.

\begin{definition}
  The group $\Inn(Q)$ of \emph{inner automorphisms} is the subgroup 
  of $\Aut(Q)$ generated by all $\rho_a$ with $a\in Q$. 
  We define the map $\inn \colon Q \to \Inn(Q)$ by $a \mapsto \rho_a$.
\end{definition}

\begin{remark}
  For every $\phi \in \Aut(Q)$ and $a \in Q$ we have 
  $\inn(a^\phi) = \phi^{-1} \circ \inn(a) \circ \phi = \inn(a)^\phi$.
  In particular, the subgroup $\Inn(Q)$ is normal in $\Aut(Q)$.
\end{remark}

\begin{notation}
  In view of the map $\inn \colon Q \to \Inn(Q)$, we also write $a^b$ 
  for the operation $a \ast b = a^{\inn(b)}$ in a quandle.
  Conversely, it will sometimes be convenient to write 
  $a \ast b$ for the conjugation $b^{-1}ab$ in a group.
  In neither case will there be any danger of confusion.
\end{notation}

\begin{definition}
  A right action of a group $G$ by quandle automorphisms on $Q$ 
  is a group action $Q \times G \to Q$, $(a,g) \mapsto a^g$ such that
  $(a \ast b)^g \mapsto a^g \ast b^g$ for all $a \in Q$ and $g \in G$.  
  This is the same as a group homomorphism $h \colon G \to \Aut(Q)$
  with $h(g) \colon Q \isoto Q$, $a \mapsto a^g$.
  We say that $G$ acts by inner automorphisms if $h(G) \subset \Inn(Q)$. 
\end{definition}


\subsection{Representations and augmentations}

The following terminology has proved useful in describing
the interplay between quandles and groups.

\begin{definition}
  A \emph{representation} of a quandle $Q$ 
  in a group $G$ is a map $\phi \colon Q \to G$ such that 
  $\phi(a \ast b) = \phi(a) \ast \phi(b)$ for all $a,b\in Q$.
  In other words, a representation $Q \to G$ is 
  a quandle homomorphism $Q \to \Conj(G)$.
  An \emph{augmentation} consists of a representation $\phi \colon Q \to G$ 
  together with a group homomorphism $\alpha \colon G \to \Aut(Q)$ 
  such that $\alpha \circ \phi = \inn$.  Interpreting $\alpha$ as a right action,
  an augmentation thus satisfies $\phi(a \ast b) = \phi(a) \ast \phi(b)$ 
  as well as $\phi(a^g) = \phi(a)^g$ for all $a,b \in Q$ and $g \in G$.
\end{definition}

\begin{example} \label{exm:InnerRepresentation}
  We have $\inn(a \ast b) = \inn(a) \ast \inn(b)$,
  in other words, $\inn$ is a representation of $Q$ in $\Inn(Q)$,
  called the \emph{inner representation}.  Together with the natural 
  action of $\Inn(Q)$ on $Q$ we obtain the \emph{inner augmentation}
  $Q \lto[\inn] \Inn(Q) \lto[\mathrm{inc}] \Aut(Q)$.
\end{example}

\begin{remark}
  Since the structure of $(Q,\ast)$ can be recovered from 
  the augmentation maps $(\phi,\alpha)$ via $a \ast b = a^{\phi(b)}$, 
  one can equivalently consider a set $Q$ equipped with a group action 
  $\alpha \colon Q \times G \to Q$, $(a,g) \mapsto a^g$ and 
  a $G$-equivariant map $\phi \colon Q \to G$, $\phi(a^g) = \phi(a)^g$.
  The operation $a \ast b = a^{\phi(b)}$ defines a rack,
  and a quandle if additionally $a^{\phi(a)} = a$ for all $a \in Q$.
  The data $(G,Q,\phi,\alpha)$ has been called a \emph{crossed $G$-set} 
  by Freyd and Yetter \cite{FreydYetter:1989} and an \emph{augmented rack}
  by Fenn and Rourke \cite{FennRourke:1992}.
\end{remark}

\begin{remark}
  For an augmentation $Q \lto[\phi] G \lto[\alpha] \Aut(Q)$
  we do not require that the image quandle $\phi(Q)$ generates 
  the entire group $G$.  We can always achieve this by restricting 
  to the subgroup $H = \gen{\phi(Q)}$.  This also entails 
  $\alpha(H) = \Inn(Q)$, so that we obtain: 
  \[
  \begin{diagram}
    \node{Q}
    \arrow{s,=}
    \arrow{e,t}{\phi}
    \node{H}
    \arrow{s,J}
    \arrow{e,t,A}{\alpha|_H}
    \node{\Inn(Q)}
    \arrow{s,J}
    \\
    \node{Q}
    \arrow{e,t}{\phi}
    \node{G}
    \arrow{e,t}{\alpha}
    \node{\Aut(Q)}
  \end{diagram}
  \]
\end{remark}

\begin{example}
  Consider a quandle $Q$ that can be faithfully represented in a group $G$,
  so that we can assume $Q \subset G$ and $G = \gen{Q}$, with
  the quandle operation given by conjugation. 
  The inner representation $\inn \colon Q \to \Inn(Q)$ 
  extends to an augmentation $Q \into G \lonto[\rho] \Inn(Q)$, 
  with $\ker(\rho) = Z(G)$ and $\Inn(Q) \cong \Inn(G)$.
  In particular, $\rho \colon G \to \Inn(Q)$ is a central 
  group extension.  This observation will be 
  generalized to every augmentation below. 
\end{example}

\subsection{The adjoint group} \label{sub:AdjointGroup}

The universal representation can be constructed as follows:

\begin{definition}
  Given a quandle $Q$ we define its \emph{adjoint group} 
  $\Adj(Q) = \gen{ Q \mid R }$ to be the quotient group 
  of the group $F(Q)$ freely generated by the set $Q$
  modulo the relations induced by the quandle operation,
  $R = \{ a \ast b = b^{-1} \cdot a \cdot b \mid a,b \in Q \}$.
  By construction we obtain a canonical map 
  $\adj \colon Q \into F(Q) \onto \Adj(Q)$
  with $\adj(a \ast b) = \adj(a) \ast \adj(b)$.
\end{definition}

The group $\Adj(Q)$ can be interpreted as the ``enveloping group'' of $Q$.
Notice, however, that the map $\adj$ is in general not injective,
see Proposition \ref{prop:HeisenbergGroup} below.

\begin{remark}[universal property] \label{rem:UniversalProperty}
  The map $\adj \colon Q \to \Adj(Q)$ is the universal group representation
  of the quandle $Q$:  for every group representation $\phi \colon Q \to G$ 
  there exists a unique group homomorphism $h \colon \Adj(Q) \to G$ 
  such that $\phi = h \circ \adj$.  
\end{remark}

\begin{remark}[functoriality] \label{rem:FunctorialProperty}
  Every quandle homomorphism $\phi \colon Q \to Q'$ induces a unique
  group homomorphism $\Adj(\phi) \colon \Adj(Q) \to \Adj(Q')$
  such that $\Adj(\phi) \circ \adj_Q = \adj_{Q'} \circ \phi$.
  We thus obtain a functor $\Adj \colon \Quandles \to \Groups$.
\end{remark}

\begin{remark}[adjointness] \label{rem:AdjointProperty}
  Its name is justified by the fact that $\Adj$ is the left adjoint functor
  of $\Conj \colon \Groups \to \Quandles$, already discussed above.
  More explicitly this means that we have a natural bijection 
  $\Hom_{\Quandles} (Q,\Conj(G)) \cong \Hom_{\Groups} (\Adj(Q),G)$,
  see \cite[chap.\,IV]{MacLane:1995}.
\end{remark}

\begin{example}[adjoint action] \label{exm:AdjointAction}
  The inner representation $\inn \colon Q \to \Inn(Q)$ induces 
  a unique group homomorphism $\rho \colon \Adj(Q) \onto \Inn(Q)$ 
  such that $\inn = \rho \circ \adj$.  In this way the adjoint group 
  $\Adj(Q)$ acts on the quandle $Q$, again denoted by
  $Q \times \Adj(Q) \to Q$, $(a,g) \mapsto a^g$.
\end{example}

\begin{remark}[adjoint augmentation] \label{rem:AdjointAugmentation}
  The pair $Q \lto[\adj] \Adj(Q) \lto[\rho] \Inn(Q)$ is 
  an augmentation of the quandle $Q$ on its adjoint group $\Adj(Q)$, 
  called the \emph{adjoint augmentation}. 
  In particular, this means that $\adj(a^g) = \adj(a)^g$
  for all $a \in Q$ and $g \in \Adj(G)$. 
  By construction it is the universal augmentation, in the obvious sense.
\end{remark}


\subsection{Connected components} \label{sub:ConnectedComponents}

As is the case for many other mathematical structures, a quandle $Q$ 
is called \emph{homogeneous} if  $\Aut(Q)$ acts transitively on $Q$.  
The following definition is more specific for quandles, and 
essentially goes back to Joyce \cite[\textsection 8]{Joyce:1982}:

\begin{definition}
  A quandle $Q$ is called \emph{connected} if $\Inn(Q)$ acts transitively on $Q$.
  A \emph{connected component} of $Q$ is an orbit under the action of $\Inn(Q)$.  
  Given an element $q \in Q$ we denote by $[q]$ its connected component, 
  that is, the orbit of $q$ under the action of $\Inn(Q)$.
  Finally, we denote by $\pi_0(Q) = \{ [q] \mid q \in Q \}$ 
  the set of connected components of $Q$.
\end{definition}

\begin{proposition}[universal property]
  The set $\pi_0(Q)$ of connected components can be considered 
  as a trivial quandle, in which case the canonical projection 
  $\phi \colon Q \onto \pi_0(Q)$, $q \mapsto [q]$ becomes a quandle homomorphism.
  It is universal in the sense that every quandle homomorphism
  $Q \to X$ to a trivial quandle $X$ factors uniquely through $\phi$.
  \qed
\end{proposition}

\begin{corollary}[functoriality]
  Every quandle homomorphism $\phi \colon Q \to Q'$ induces a map 
  $\phi_* \colon \pi_0(Q) \to \pi_0(Q')$ defined by $[x] \mapsto [\phi(x)]$.
  If $\phi$ is surjective then so is $\phi_*$.  In particular, 
  the homomorphic image of a connected quandle is again connected.
  \qed
\end{corollary}

\begin{remark}
  For every quandle $Q$, the elements of a given component become 
  conjugate in $\Adj(Q)$.  Its abelianization is thus given by 
  $\alpha \colon \Adj(Q) \to \Z\pi_0(Q)$, $q \mapsto [q]$, and
  its kernel is the commutator subgroup $\Adj(Q)' = \Ker(\alpha)$.
\end{remark}

\begin{definition} \label{def:Index}
  For every quandle $Q$ there exists a unique group homomorphism 
  $\varepsilon \colon \Adj(Q) \to \Z$ with $\adj(Q) \to \{1\}$.  
  Its kernel $\Adj(Q)^\circ := \Ker(\varepsilon)$ is generated by 
  all products of the form $\adj(a)^{-1} \adj(b)$ with $a,b \in Q$.
  The image of $\Adj(Q)^\circ$ under the natural group homomorphism 
  $\Adj(Q) \to \Inn(Q)$ will be denoted by $\Inn(Q)^\circ$.   
  It is generated by products of the form $\inn(a)^{-1} \inn(b)$,
  called \emph{transvections} by Joyce \cite[\textsection 5]{Joyce:1982}.
  In his analysis of symmetric spaces \'E.\,Cartan called this
  the \emph{group of deplacements} (see Loos \cite[\textsection II.1.1]{Loos:1969}).
\end{definition}

\begin{remark} \label{rem:CommutatorOrbits}
  If $Q$ is connected, then $\varepsilon \colon \Adj(Q) \to \Z$
  is the abelianization of the adjoint group, and in this case 
  $\Adj(Q)^\circ = \Adj(Q)'$ and $\Inn(Q)^\circ = \Inn(Q)'$.

  We have $\Adj(Q) = \Adj(Q)^\circ \rtimes \Z$: choosing a base point
  $q \in Q$, every element $g \in \Adj(Q)$ can be uniquely written
  as $g = \adj(q)^{\varepsilon(g)} h$ with $h \in \Adj(Q)^\circ$.

  The components of $Q$ are the orbits under 
  the adjoint action of $\Adj(Q)$.  We obtain the same 
  orbits with respect to the subgroup $\Adj(Q)^\circ$.
  Indeed, for $a \in Q$ and $g \in \Adj(Q)$ we have $a^g = a^h$ 
  with $h = \adj(a)^{-\varepsilon(g)} g \in \Adj(Q)^\circ$.
\end{remark}

\begin{Notes}
  \begin{remark}
    If $Q$ is not connected, then the orbits under
    $\Adj(Q)^\circ$ and $\Adj(Q)'$ usually differ significantly:
    Consider the quandle $Q = Q_{m,n}$ of Example \ref{exm:DoubleCycle}
    with $\gcd(m,n)=1$, where we find $\Adj(Q) = \Z \times \Z$ and 
    $\Inn(Q) = \Z_n \times \Z_m$.  The orbits under $\Adj(Q)^\circ \cong \Z$ 
    are the two connected components, and do thus not coincide 
    with the orbits under the trivial group $\Adj(Q)' = \{\id\}$.
  \end{remark}
\end{Notes}

\subsection{Central group extensions}

Fenn and Rourke \cite{FennRourke:1992} have called the kernel 
of the natural group homomorphism $\rho \colon \Adj(Q) \onto \Inn(Q)$ 
the \emph{excess} of $Q$, but did not study $\rho$ more closely.
We will now see that $\rho$ is a central extension.

As for every group, the inner automorphism group $\Inn(\Adj(Q))$ is the image 
of the homomorphism $\gamma \colon \Adj(Q) \onto \Aut(\Adj(Q))$ defined
by conjugation, $\gamma(g) \colon x \mapsto x^g$, and its kernel is the centre 
of $\Adj(Q)$.  By definition of the adjoint group, we also have a homomorphism 
$\alpha \colon \Aut(Q) \to \Aut(\Adj(Q))$ given by $\phi \mapsto \Adj(\phi)$.

\[
\begin{diagram}
  \node{Q}
  \arrow{e,t}{\inn}
  \arrow{s,l}{\adj}
  \node{\Inn(Q)}
  \arrow{e,t,J}{}
  \arrow{s,r,A}{\beta}
  \node{\Aut(Q)}
  \arrow{s,r}{\alpha}
  \\
  \node{\Adj(Q)}
  \arrow{ne,t,A}{\rho}
  \arrow{e,t,A}{\gamma}
  \node{\Inn \Adj(Q)}
  \arrow{e,t,J}{}
  \node{\Aut \Adj(Q)}
\end{diagram}
\]

\begin{proposition} \label{prop:CentralExtension}
  We have $\alpha(\Inn(Q)) = \Inn\Adj(Q)$. The restriction of $\alpha$ 
  defines a group homomorphism $\beta \colon \Inn(Q) \onto \Inn\Adj(Q)$
  that makes the above diagram commute.  As a consequence, the group 
  homomorphism $\rho \colon \Adj(Q) \onto \Inn(Q)$ is a central extension.
\end{proposition}

\begin{proof}
  We already have $\inn = \rho \circ \adj$ by construction of $\rho$,
  so we only have to verify that $\alpha \circ \rho = \gamma$.
  Every $g \in \Adj(Q)$ acts on $Q$ by inner automorphisms,
  $\rho(g) \colon Q \isoto Q$, $a \mapsto a^g$.  
  The quandle automorphism $\rho(g)$ induces a group automorphism 
  $\Adj\rho(g) \colon \Adj(Q) \isoto \Adj(Q)$
  with $\adj(a) \mapsto \adj(a^g) = \adj(a)^g$,
  cf.\ Remark \ref{rem:AdjointAugmentation}.
  We conclude that $\Adj\rho(g) = \gamma(g)$.
  This means that the diagram is commutative
  and $\alpha(\Inn(Q)) = \Inn\Adj(Q)$.
\end{proof}

As an illustration we wish to determine the adjoint group 
of the quandle $Q_{m,n} = \Z_m \sqcup \Z_n$ from Example \ref{exm:DoubleCycle}.
Recall that it decomposes into two components, $\Z_m$ and $\Z_n$.

\begin{proposition} \label{prop:HeisenbergGroup}
  The adjoint group $\Adj(Q_{0,0})$ is isomorphic to the Heisenberg group
  \[
  H = \left\{ \left(\begin{smallmatrix} 1 & * & * \\ 0 & 1 & * \\ 0 & 0 & 1
  \end{smallmatrix}\right) \in \SL_3\Z \right\} \quad\text{generated by}\quad
  x = \left(\begin{smallmatrix} 1 & 1 & 0 \\ 0 & 1 & 0 \\ 0 & 0 & 1 \end{smallmatrix}\right),\;
  y = \left(\begin{smallmatrix} 1 & 0 & 0 \\ 0 & 1 & 1 \\ 0 & 0 & 1 \end{smallmatrix}\right),\;
  z = \left(\begin{smallmatrix} 1 & 0 & 1 \\ 0 & 1 & 0 \\ 0 & 0 & 1 \end{smallmatrix}\right).
  \]

  More generally, for parameters $m,n \in \N$ 
  the adjoint group $G = \Adj(Q_{m,n})$ is isomorphic 
  to the quotient $H_\ell = H/\gen{z^\ell}$ with $\ell = \gcd(m,n)$,
  via the isomorphism $\phi \colon G \isoto H_\ell$ 
  defined by $\adj(a) \mapsto xz^a$ for $a \in \Z_m$ 
  and $\adj(b) \mapsto yz^{-b}$ for $b \in \Z_n$.

  In particular, $\adj \colon Q_{m,n} \to \Adj(Q_{m,n})$ 
  is injective if and only if $m = n$, and 
  we have $\Adj(Q_{m,n}) \cong \Z \times \Z$ 
  if and only if  the parameters $m$ and $n$ are coprime. 
\end{proposition}

\begin{proof}
  By definition, the adjoint group $G = \Adj(Q_{m,n})$ 
  is generated by elements $s_a$ with $a \in \Z_m$ and 
  $t_b$ with $b \in \Z_n$ subject to the quandle relations
  $s_a \ast t_b = s_{a+1}$ and $t_b \ast s_a = t_{b+1}$,
  as well as $s_a \ast s_{a'} = s_a$ and $t_b \ast t_{b'} = t_b$
  for all $a,a' \in \Z_m$ and $b,b' \in \Z_n$.

  In $H$ we have $[x,y] = x^{-1} y^{-1} x y = z$ and $[x,z] = [y,z] = 1$, 
  which entails the desired relations $(xz^a) \ast (yz^{-b}) = xz^{a+1}$
  and $(yz^{-b}) \ast (xz^a) = yz^{-(b+1)}$. The quotient group 
  $H_\ell = H/\gen{z^\ell}$ thus allows a quandle representation
  $Q_{m,n} \to H_\ell$ with $a \mapsto xz^a$ for $a \in \Z_m$
  and $b \mapsto yz^{-b}$ for $b \in \Z_n$.  This induces
  a surjective group homomorphism $\phi \colon G \onto H_\ell$.


  Since $\Inn(Q_{m,n}) \cong \Z_n \times \Z_m$ is abelian, 
  the commutator group $G'$ is contained in the kernel of $G \onto \Inn(Q)$, 
  which is central according to Proposition \ref{prop:CentralExtension}.
  Consider 
  \[
  u := [s_a,t_b] = s_a^{\smash{-1}} t_b^{\smash{-1}} s_a t_b 
  = s_a^{\smash{-1}} s_{a+1} .
  \]
  Repeatedly conjugating this equation by $t_b$ yields
  \[
  u = s_{a}^{\smash{-1}} s_{a+1} = s_{a+1}^{\smash{-1}} s_{a+2} 
  = \dots = s_{a-1}^{\smash{-1}} s_{a} .
  \]
  On the other hand we find $u = t_{b+1}^{\smash{-1}} t_b$ 
  and repeatedly conjugating by $s_a$ yields
  \[
  u = t_{b+1}^{\smash{-1}} t_{b} = t_{b+2}^{\smash{-1}} t_{b+1}
  = \dots = t_{b}^{\smash{-1}} t_{b-1} .
  \]
  This shows that $u^m = u^n = 1$ and thus $u^\ell = 1$ for $\ell = \gcd(m,n)$.
  With $s := s_0$ and $t := t_0$ we finally obtain 
  $s_a = s u^a$ for all $a \in \Z_m$ and 
  $t_b = t u^{-b}$ for all $b \in \Z_n$.
  We conclude that every element of $G$ can be written 
  as $s^i t^j u^k$ with $i,j \in \Z$ and $k \in \Z_\ell$.
  The group homomorphism $\phi \colon G \onto H_\ell$ satisfies 
  $\phi(s^i t^j u^k) = x^i y^j z^k$, and is thus seen to be injective.
\end{proof}

\begin{remark}
  The natural group homomorphism $\beta \colon \Inn(Q) \onto \Inn\Adj(Q)$ 
  is surjective but in general not injective.  
  Consider for example $Q = Q_{m,n}$ with $\gcd(m,n)=1$.
  Then $\Adj(Q) \cong \Z^2$, so $\Inn\Adj(Q) = \{\id\}$, 
  whereas $\Inn(Q) \cong \Z_m \times \Z_n$.
\end{remark}

\begin{remark}
  The group homomorphism $\alpha \colon \Aut(Q) \to \Aut\Adj(Q)$
  is in general neither injective nor surjective.  The trivial quandle $Q = \{q\}$,
  for example, has trivial automorphism group $\Aut(Q) = \{\id\}$, 
  whereas the adjoint group $\Adj(Q) \cong \Z$ has $\Aut\Adj(Q) = \{\pm\id\}$.
\end{remark}

\begin{corollary} \label{cor:CentralExtension}
  For every augmentation $Q \lto[\phi] G \lto[\alpha] \Inn(Q)$
  with $G = \gen{\phi(Q)}$, the induced group homomorphism
  $h \colon \Adj(Q) \onto G$ and $\alpha \colon G \onto \Inn(Q)$ 
  are central extensions, because $\alpha \circ h = \rho$
  is a central extension according to Proposition \ref{prop:CentralExtension}.
  \qed
  \[
  \begin{diagram}
    \node[2]{\Adj(Q)}
    \arrow{se,t,A}{h}
    \arrow{s,r,3,-}{\rho}
    \\
    \node{Q}
    \arrow{ne,t}{\adj}
    \arrow[2]{e,t,1}{\phi}
    \arrow{se,b}{\inn}
    \node{}
    \arrow{s,A}
    \node{G}
    \arrow{sw,r,A}{\alpha}
    \\
    \node[2]{\Inn(Q)}
  \end{diagram}
  \]
\end{corollary}


\subsection{Functoriality}

Unlike the adjoint representation, 
the inner representation $Q \to \Inn(Q)$ is \emph{not} functorial.
To see this, it suffices to consider a quandle $Q'$ and an element 
$q' \in Q'$ that acts non-trivially, i.e.\ $\inn(q') \ne \id_{Q'}$.
The trivial quandle $Q = \{q\}$ maps into $Q'$ with $q \mapsto q'$,
but no group homomorphism $\Inn(Q) \to \Inn(Q')$ can map 
$\inn(q) = \id_Q$ to $\inn(q') \ne \id_{Q'}$.

A closer look reveals that the crucial hypothesis is surjectivity:

\begin{proposition}
  For every surjective quandle homomorphism $p \colon Q \onto \bar{Q}$ 
  there exists a unique group homomorphism 
  $h \colon \Inn(Q) \onto \Inn(\bar{Q})$
  that makes the following diagram commutative:
  \[
  \begin{diagram}
    \node{Q}
    \arrow{s,l,A}{p}
    \arrow{e,t}{\inn_{Q}}
    \node{\Inn(Q)}
    \arrow{s,r,A}{h = \Inn(p)}
    \\
    \node{\bar{Q}}
    \arrow{e,t}{\inn_{\bar{Q}}}
    \node{\Inn(\bar{Q})}
  \end{diagram}
  \]
\end{proposition} 

\begin{proof}
  Uniqueness is clear because $\Inn(Q) = \gen{\inn(Q)}$.
  In order to prove existence, first observe that for each $a \in Q$
  the inner action $x \mapsto x \ast a$ preserves the fibres of $p$.
  The same is thus true for every $g \in \Inn(Q)$, so we obtain 
  a well-defined map $\bar{g} \colon \bar{Q} \to \bar{Q}$ 
  as follows: for each $\bar{x}$ choose a preimage $x \in Q$ 
  with $p(x) = \bar{x}$ and set $\bar{x}^{\bar{g}} := p(x^g)$.
  By construction we have $\overline{f \circ g} = \bar{f} \circ \bar{g}$,
  and $g = \inn(a)$ is mapped to $\bar{g} = \inn(p(a))$.  
  This shows that the map $h \colon \Inn(Q) \to \Inn(\bar{Q})$, 
  $g \mapsto \bar{g}$, is well-defined and a surjective group homomorphism.
\end{proof}

\subsection{Quandle coverings} \label{sub:QuandleCoverings}

The following definition of quandle covering was inspired 
by \cite{Eisermann:2003}, where this approach 
was successfully used to study knot quandles.

\begin{definition} \label{def:QuandleCovering}
  A quandle homomorphism $p \colon \tilde{Q} \to Q$ is called 
  a \emph{covering} if it is surjective and $p(\tilde{x}) = p(\tilde{y})$ 
  implies $\tilde{a} \ast \tilde{x} = \tilde{a} \ast \tilde{y}$ 
  for all $\tilde{a},\tilde{x},\tilde{y} \in \tilde{Q}$.
\end{definition}

\begin{example}
  For every augmentation $Q \lto[\phi] G \lto[\alpha] \Aut(Q)$
  the homomorphism $\phi \colon Q \onto \phi(Q)$
  is a quandle covering.  In particular,
  the inner representation $\inn \colon Q \to \Inn(Q)$ defines 
  a quandle covering $Q \onto \inn(Q)$.  By definition, $\inn(Q)$ 
  is the smallest quandle covered by $Q$. In the other extreme 
  we will show in Section \ref{sec:ConnectedCoverings} below 
  how to construct the universal covering of $Q$.
\end{example}

\begin{notation}
  We shall reserve the term ``covering'' for the map $p \colon \tilde{Q} \to Q$. 
  If emphasis is desired, it is convenient to call $p \colon \tilde{Q} \to Q$
  the \emph{quandle covering} and $\tilde{Q}$ the \emph{covering quandle}.%
\end{notation}

\begin{example} \label{exm:ConjCovering}
  A surjective group homomorphism $p \colon \tilde{G} \onto G$ yields 
  a quandle covering $\Conj(\tilde{G}) \to \Conj(G)$ if and only if 
  $\ker(p) \subset \tilde{G}$ is a central subgroup.
\end{example}

\begin{Notes}
  \begin{proof}
    $(\Leftarrow)$ 
    If $\ker(p) \subset Z(\tilde{G})$, 
    then $\Conj(p)$ is a quandle covering: for $z \in \ker(p)$ 
    we find $a \ast (zb) = (zb)^{-1} a (zb) = b^{-1}ab = a \ast b$.

    $(\Rightarrow)$
    Conversely, let us assume that $a \ast b = b^{-1}ab$ 
    equals $a \ast zb$ for all $a,b \in G$ and $z \in \ker(p)$. 
    For $a=b$ we find $a = a \ast a = a \ast za = a^{-1} z^{-1} a z a$, 
    which implies that $z \in Z(\tilde{G})$.
  \end{proof}
\end{Notes}

\begin{example} \label{exm:CoreCovering}
  A surjective group homomorphism $p \colon \tilde{G} \onto G$ yields
  a quandle covering $\Core(\tilde{G}) \to \Core(G)$ if and only if 
  $\ker(p) \subset \tilde{G}$ is a central subgroup of exponent $2$.
\end{example}

\begin{Notes}
  \begin{proof}
    $(\Leftarrow)$ 
    If $\ker(p) \subset Z(\tilde{G})$ if of exponent $2$,
    then $\Core(p)$ is a quandle covering: for $z \in \ker(p)$ 
    we find $a \ast (zb) = zb \cdot a^{-1} zb = ba^{-1}b = a \ast b$.

    $(\Rightarrow)$
    Conversely, let us assume that  $a \ast b = ba^{-1}b$ 
    equals $a \ast zb$ for all $a,b \in G$ and $z \in \ker(p)$.
    For $a=b$ we find $a = a \ast a = a \ast za = z^2 a$, 
    which is equivalent to $z^2 = 1$.  
    For $a=1$ we have $b^2 = 1 \ast b = 1 \ast zb = z b z b$, 
    which is equivalent to $z \in Z(\tilde{G})$.
  \end{proof}
\end{Notes}

\begin{example} \label{exm:AlexCovering}
  A surjective group homomorphism $p \colon \tilde{A} \onto A$ 
  with $p \circ \tilde{T} = T \circ p$ yields a quandle covering 
  $\Alex(\tilde{A},\tilde{T}) \to \Alex(A,T)$ if and only if 
  $\tilde{T}$ acts trivially on $\ker(p) \subset \tilde{A}$.
\end{example}

\begin{Notes}
  \begin{proof}
    $(\Leftarrow)$ 
    If $T(z) = z$ for all $z \in \ker(p)$,
    then $a \ast zb = T(a (zb)^{-1})zb = T(ab^{-1})b = a \ast b$.

    $(\Rightarrow)$
    Conversely, let us assume that  $a \ast b = T(ab^{-1})b$ 
    equals $a \ast zb$ for all $a,b \in G$ and $z \in \ker(p)$.
    For $a=b$ we find $a = a \ast a = a \ast za = T(z^{-1}) z a$, 
    whence $T(z) = z$.  
  \end{proof}
\end{Notes}


\begin{warning}
  The composition of two central group extensions is in general
  not a central extension, and so the functor $\Conj$ shows
  that we cannot generally expect the composition 
  of two quandle coverings to be again a covering 
  (see also Example \ref{exm:DoubleDoubleCovering}).
  Similar remarks apply to the functors $\Core$ and $\Alex$.
\end{warning}

\begin{remark} \label{rem:CoveringAction}
  A surjective quandle homomorphism $p \colon \tilde{Q} \onto Q$
  is a covering if and only if the inner representation 
  $\inn \colon \tilde{Q} \to \Inn(\tilde{Q})$ factors through $p$:
  we thus obtain a representation $\tilde{\sigma} \colon Q \to \Inn(\tilde{Q})$
  by setting $\tilde{a} \ast x := \tilde{a} \ast \tilde{x}$ for all $x \in Q$ 
  and $\tilde{a},\tilde{x} \in \tilde{Q}$ with $p(\tilde{x}) = x$.
  This is well-defined because $\tilde{a} \ast \tilde{x}$
  does not depend on the choice of the preimage $\tilde{x}$.
  Moreover, $\tilde{\sigma}$ induces a group homomorphism 
  $\tilde{\rho} \colon \Adj(Q) \to \Inn(\tilde{Q})$.  
  This situation is summarized in the following commutative diagram:
  \[
  \begin{diagram}
    \node{\tilde{Q}}
    \arrow[2]{e,t}{\adj_{\tilde{Q}}}
    \arrow[2]{s,l,A}{p}
    \arrow{se,b}{\inn_{\tilde{Q}}}
    \node{}
    \node{\Adj(\tilde{Q})}
    \arrow[2]{s,r,A}{\Adj(p)}
    \arrow{sw,b,A}{\rho_{\tilde{Q}}}
    \\
    \node{}
    \node{\Inn(\tilde{Q})}
    \arrow[2]{s,r,1,A}{\Inn(p)}
    \node{}
    \\
    \node{Q}
    \arrow{e,b,3,-}{\adj_{Q}}
    \arrow{ne,t,..}{\tilde{\sigma}}
    \arrow{se,b}{\inn_Q}
    \node{}
    \arrow{e}    
    \node{\Adj(Q)}
    \arrow{nw,t,A,..}{\tilde{\rho}}
    \arrow{sw,b,A}{\rho_Q}
    \\
    \node{}
    \node{\Inn(Q)}
    \node{}
  \end{diagram}
  \]
  
  In particular, $\tilde{\rho}$ defines a natural action of the adjoint group $\Adj(Q)$ 
  on the covering quandle $\tilde{Q}$, and $p$ is equivariant with respect to this action.
\end{remark}

\begin{proposition} \label{prop:CoveringAdjInnExtention}
  For every quandle covering $p \colon \tilde{Q} \onto Q$, the induced
  group homomorphisms $\Adj(p) \colon \Adj(\tilde{Q}) \onto \Adj(Q)$ 
  and $\Inn(p) \colon \Inn(\tilde{Q}) \onto \Inn(Q)$ are central extensions.
\end{proposition}

\begin{proof}
  This follows from the commutativity of the diagram 
  and Proposition \ref{prop:CentralExtension}.
\end{proof}


\section{Examples of quandles and coverings} \label{sec:Examples}

This section recalls some classical examples where quandles
arise naturally: conjugation in groups, the adjoint action
of a Lie group on its Lie algebra, and the symmetries of a 
Riemannian symmetric space.  Our aim here is to highlight 
the notion of quandle covering and its relationship to 
central group extensions, coverings of Lie groups,
and coverings of symmetric spaces, respectively.
A more detailed analysis and the calculation of adjoint 
and fundamental groups will be the object of a forthcoming
article \cite{Eisermann:Wirtinger}.

\subsection{Trivial coverings}

Even though this is by far the least interesting case,
we shall start our tour with trivial coverings.

\begin{example}[trivial covering] \label{exm:TrivialCovering}
  Let $Q$ be a quandle and let $F$ be a non-empty set.
  We can consider $F$ as a trivial quandle, and equip
  the product $\tilde{Q} = Q \times F$ with the quandle
  operation $(a,s) \ast (b,t) = (a \ast b,s)$.
  The projection $p \colon Q \times F \to Q$ 
  given by $(q,s) \mapsto q$ is a quandle covering, 
  called \emph{trivial covering} with fibre $F$.
\end{example}

\begin{remark}
  For every quandle homomorphism $p \colon \tilde{Q} \to Q$,
  each fibre $F = p^{-1}(q)$ is a subquandle of $\tilde{Q}$.
  If $p$ is a quandle covering, then $F$ is necessarily trivial.
  The fibres over any two points of the same component are isomorphic.  
  The isomorphism is not canonical, however, and covering theory 
  studies the possible monodromy.
\end{remark}

\begin{remark}[almost trivial covering] \label{rem:AlmostTrivialCovering}
  If $Q$ decomposes into connected components $(Q_i)_{i \in I}$, 
  then we can choose a non-empty set $F_i$ for each $i \in I$
  and equip the union $\tilde{Q} = \bigsqcup_{i \in I} Q_i \times F_i$ 
  with the previous quandle operation $(a,s) \ast (b,t) = (a \ast b,s)$. 
  The result is a quandle covering $\tilde{Q} \to Q$, $(q,s) \mapsto q$ 
  that is trivial over each component, but not globally trivial
  if the fibres over different components are non-isomorphic.
\end{remark}


\subsection{Conjugation quandles}

As already noted in the introduction, every group $G$ becomes 
a quandle with respect to conjugation $a \ast b = b^{-1} a b$.
More generally, every non-empty union $Q$ of conjugacy classes in $G$ 
is a quandle with these operations, and $Q$ is a connected quandle if and only if 
$Q$ is a single conjugacy class in the generated subgroup $H = \gen{Q}$.

\begin{remark}[central extensions] \label{rem:CentralExtensions}
  Given a quandle $Q \subset G$ and a central group extension 
  $p \colon \tilde{G} \onto G$, the preimage $\tilde{Q} = p^{-1}(Q)$
  yields a quandle covering $p \colon \tilde{Q} \to Q$.
  The kernel $\Lambda = \ker(p)$ acts on the covering quandle $\tilde{Q}$ such that
  $(\lambda a) \ast b = \lambda(a \ast b)$ and $a \ast (\lambda b) = a \ast b$ 
  for all $a,b \in \tilde{Q}$ and $\lambda \in \Lambda$.
  This will be called a quandle extension,
  see Definition \ref{def:QuandleExtension}.
\end{remark}


\begin{example}[linear groups] \label{exm:LinearGroups1}
  Consider the special linear group $\SL_2\K$ over a field $\K$.
  Its centre is $Z = \{\pm\id\}$ and thus of order $2$ 
  if $\operatorname{char} \K \ne 2$.  The quotient is the projective 
  special linear group $\PSL_2\K = \SL_2\K/Z$, and by construction 
  $p \colon \SL_2\K \onto \PSL_2\K$ is a central extension.
  We will assume that $|\K| \ge 4$, so that $\SL_2\K$ is perfect 
  and $\PSL_2\K$ is simple.  (See \cite[\textsection XIII.8]{Lang:2002}.)

  The conjugacy class $\tilde{Q} = \tilde{q}^{\tilde{G}}$ 
  of $\tilde{q} = \left(\begin{smallmatrix} 
      0 & 1 \\ -1 & 0 \end{smallmatrix}\right)$ 
  defines a quandle in $\tilde{G} = \SL_2\K$.
  Its image $Q := p(\tilde{Q}) = q^{G}$ is the conjugacy class 
  of $q := p(\tilde{q}) = \pm \tilde{q}$ in $G = \PSL_2\K$.
  We have $G = \gen{Q}$ because $G$ is simple, and 
  $\tilde{G} = \gen{\tilde{Q}}$ because $\tilde{G}$ is perfect.
  (This is a general observation: $\gen{\tilde{Q}}$ is normal in $\tilde{G}$
  and maps onto $G$, so that $\tilde{G}/\gen{\tilde{Q}}$ is abelian, 
  whence $\tilde{G} = \gen{\tilde{Q}}$.)

  Suppose that there exist $a,b \in \K$ such that $a^2 + b^2 = -1$. 
  (This always holds in finite characteristic, 
  and also for $\K = \C$, but not for $\K = \R$.)
  In this case the matrix $c = \left(\begin{smallmatrix} 
      a & b \\ b & -a \end{smallmatrix}\right) \in \tilde{G}$ 
  conjugates $\tilde{q}$ to $\tilde{q}^c = -\tilde{q}$, 
  so that $Z \cdot \tilde{Q} = \tilde{Q}$.
  This means that $p \colon \tilde{Q} \onto Q$ 
  is a two-fold covering of connected quandles,
  and even an extension $Z \curvearrowright \tilde{Q} \onto Q$.

  If $a^2 + b^2 = -1$ has no solution in $\K$, as for example in $\K = \R$,
  then $\tilde{q}$ and $-\tilde{q}$ are not conjugated in $\tilde{G} = \SL_2\K$,
  so that $p^{-1}(Q) = +\tilde{Q} \sqcup -\tilde{Q}$ consists 
  of two isomorphic copies of $\tilde{Q}$.  This is again
  a two-fold quandle covering, but a trivial one.
\end{example}

\subsection{Lie groups and Lie algebras}

Every Lie group $G$ is tied to its Lie algebra $\g = T_1 G$ 
by two important maps (see for example \cite[chap.\,II]{Helgason:2001}):
\begin{enumerate}
\item
  The exponential map $\exp \colon \g \to G$, 
  characterized by its usual differential equation 
  $\frac{d}{dt} \exp(tx) = T_1 \exp(tx) \cdot x$
  with initial condition $\exp(0) = 1$.
\item
  The adjoint action $\ad \colon G \to \Aut(\g)$, 
  denoted $\ad(g) \colon x \mapsto x^g$, which can be 
  defined as the derivative at $1$ of the conjugation map
  $\gamma_g \colon G \to G$, $h \mapsto h \ast g = g^{-1} h g$.
\end{enumerate}

This naturally corresponds to a quandle structure in the following sense:

\begin{itemize}
\item
  The set $\g$ is a quandle with respect to $x \ast y = x^{\exp(y)}$.
  \\
  We recover the Lie bracket as the derivative 
  $\frac{d}{dt}\left[ x * ty \right]_{t=0} = [x,y]$.
\item
  The triple $\g \lto[\exp] G \lto[\ad] \Aut(\g)$ 
  is an augmentation of the quandle $(\g,\ast)$.
  \\
  The image $Q = \exp(\g)$ is a quandle in the group $G$,
  with respect to conjugation.
\item
  In general we have $\exp(\g) \subsetneq G$.
  If $G$ is connected and $\exp \colon \g,0 \to G,1$ 
  is a local diffeomorphism, then we have 
  $G = \gen{\exp(\g)}$ and $\ad(G) = \Inn(\g,\ast)$.
\end{itemize}

\begin{remark}
  In the finite-dimensional case, the manifold $G$ is modelled 
  on $\R^n$ or $\C^n$, and the inverse function theorem ensures
  that $\exp$ is a local diffeomorphism from an open neighbourhood 
  of $0 \in \g$ onto an open neighbourhood of $1 \in G$.  
  In the infinite-dimensional case, this result still holds 
  for Lie groups modelled on Banach spaces.  It may fail, 
  however, for complete locally convex vector spaces, 
  a setting motivated and studied by Milnor \cite{Milnor:1984}.  
  He notes that in some cases the conclusion $G = \gen{\exp(\g)}$ 
  follows from the additional property that the group $G$ is simple, 
  because $\gen{\exp(\g)}$ is a normal subgroup.
\end{remark}

\begin{remark}[central extensions again]
  If $p \colon \tilde{G} \to G$ is a connected covering of the Lie group $G$,
  then $\tilde{G}$ carries a unique Lie group structure such that $p$ is 
  a Lie group homomorphism.  The linear isomorphism $T_1 p \colon T_1 \tilde{G} \isoto T_1 G$
  provides an isomorphism of Lie algebras $\tilde{\g} \isoto \g$, and so 
  we obtain another augmentation $\g \lto[\exp] \tilde{G} \lto[\ad] \Aut(\g)$.
  This can be summarized as follows:
  \[
  \begin{diagram}
    \node{\tilde{\g}}
    \arrow{e,t,A}{\exp_{\tilde{G}}}
    \arrow{s,l}{\cong}
    \node{\tilde{Q}}
    \arrow{e,t,J}{\mathrm{inc}}
    \arrow{s,l,A}{p}
    \node{\tilde{G}}
    \arrow{e,t,A}{\ad_{\tilde{G}}}
    \arrow{s,l,A}{p}
    \node{\Inn(\tilde{\g})}
    \arrow{s,r}{\cong}
    \\
    \node{\g}
    \arrow{e,t,A}{\exp_{G}}
    \node{Q}
    \arrow{e,t,J}{\mathrm{inc}}
    \node{G}
    \arrow{e,t,A}{\ad_{G}}
    \node{\Inn(\g)}
  \end{diagram}
  \]
  Assuming $G = \gen{Q}$ and $\tilde{G} = \gen{\tilde{Q}}$,
  we recover a well-known fact of Lie group theory:
  $p \colon \tilde{G} \to G$ is a central group extension, 
  because both $\tilde{G}$ and $G$ are intermediate to 
  the central extension $\Adj(\g,\ast) \onto \Inn(\g,\ast)$,
  see Corollary \ref{cor:CentralExtension}.
  In particular, $p \colon \tilde{Q} \to Q$ 
  is a quandle covering, see Remark \ref{rem:CentralExtensions}.
\end{remark}

\subsection{Infinite-dimensional Lie algebras}

Contrary to the finite-dimensional case,
not every infinite-dimensional Lie algebra $(L,[,])$
can be realized as the tangent space of a Lie group $G$.
This fails even for Banach Lie algebras, as remarked 
by van\,Est and Korthagen \cite{VanEstKorthagen:1964}.
(See also Serre \cite{Serre:1964}, Part II, \textsection V.8.) 
It is worth noting that the construction of 
the quandle $(L,\ast)$ can still be carried out.

The obvious idea is to define $x \ast y$ 
by the initial condition $x \ast 0 = x$ and the differential 
equation $\frac{d}{dt} ( x \ast ty ) = [ x \ast ty, y ]$.
This equation has at most one analytic solution, namely
\[ 
x \ast y = \sum_{k=0}^\infty \frac{1}{k!} [\dots[[x,y],y]\dots,y] .
\]
In order to ensure convergence, it suffices to impose
some reasonable condition on the topology of $L$:
all obstacles disappear, for example, if $L$ is a Banach Lie algebra.
It is then an amusing exercise to verify that $(L,\ast)$ is indeed a quandle:
\begin{enumerate} 
\item[(Q1)]
  Antisymmetry $[x,x]=0$ translates to idempotency $x \ast x = x$.
\item[(Q2)]
  The functional equation $\exp(y) \circ \exp(-y) = \id$ 
  ensures invertibility.
\item[(Q3)]
  The Jacobi identity of the Lie bracket $[,]$ translates 
  to self-distributivity of the quandle operation $\ast$.
\end{enumerate}

We conclude that constructing the quandle $(L,\ast)$ 
is a rather benign topological problem.  The natural group 
that appears here is $G = \gen{\exp(L)} = \Inn(L,\ast)$,
but in general this need not be a Lie group.
The much deeper problem of constructing a Lie group $G$ 
realizing the Lie algebra $L$ involves the Lie algebra structure 
of $L$ in a more profound way and will in general lead to non-trivial obstructions.
The lesson to be learned from this excursion is that although 
a Lie group $G$ may be too much to ask, the less ambitious 
quandle structure $(L,\ast)$ can still be rescued.

\subsection{Reflection quandles}

Consider $\R^n$ with $a \ast b = a \tsa b = 2b-a$, 
which is the symmetry about the point $b$.
This defines a connected involutory quandle $Q = (\R^n,\ast)$,
called the $n$-dimensional \emph{reflection quandle}.
Since $b$ is the unique fix-point of $\inn(b)$, 
we see that $\inn \colon Q \to \Inn(Q)$ is injective.
More precisely, $(\R^n,\ast)$ is isomorphic to conjugacy class
of reflections in the semidirect product 
$\Inn(\R^n,\ast) \cong (\R^n,+) \rtimes \{\pm\id\}$.

\begin{example} \label{exm:TorusCovering1}
  The quandle structure passes to the quotient group $\T^n = \R^n/\Z^n$, 
  where it can again be formulated as $a \ast b = 2b-a$.
  In this way the torus $\T^n$ inherits a unique quandle 
  structure such that the projection $p \colon \R^n \to \T^n$ 
  is a quandle homomorphism.  The quotient map $p$ is \emph{not} 
  a quandle covering, because $\inn_Q$ is injective
  and does not factor through $p$.
\end{example}

\begin{example} \label{exm:TorusCovering2}
  We can produce quandle coverings $\T^n \to \T^n$ as follows.
  Consider the subgroup $\Lambda = p(\frac{1}{2}\Z^n) = \{[0],[\frac{1}{2}]\}^n$ 
  acting on $\T^n$ by translation.  For $b,b' \in \T^n$ we have
  $\inn(b) = \inn(b')$ if and only if $b-b' \in \Lambda$.
  The quotient $\Lambda\backslash\T^n$ carries a unique quandle structure 
  such that the projection $\T^n \onto \Lambda\backslash\T^n$ is a quandle covering.
  (This quotient can be identified with $\T^n \lto[_{\cdot 2}] \T^n$.)
  Similar remarks apply to the quotient by any subgroup of $\Lambda$.
\end{example}


\subsection{Spherical quandles}

We can equip the unit sphere $\S^n \subset \R^{n+1}$ 
with the operation $a \ast b = 2 \bil{a}{b} b - a$,
which is the unique involution fixing $b$ 
and mapping $x \mapsto -x$ for $x$ orthogonal to $b$.
This turns $(\S^n,\ast)$ into a connected involutory quandle,
called the $n$-dimensional \emph{spherical quandle}.

\begin{example} \label{exm:ProjectiveQuandle}
  For $\lambda = \pm1$ and $a,b \in \S^n$ we have 
  $(\lambda a) \ast b = \lambda(a \ast b)$
  and $a \ast (\lambda b) = a \ast b$.
  This means that the projective space $\RP^n = \S^n/\{\pm1\}$
  inherits a unique quandle structure $[a] \ast [b] = [a \ast b]$
  such that the projection $p \colon \S^n \to \RP^n$ is a quandle covering.  
  The map $p$ is, of course, also a covering of topological spaces.
\end{example}

\begin{remark}
  The inner action defines a representation of the quandle 
  $(\S^n,\ast)$ in the orthogonal group $\operatorname{O}(n+1)$, 
  and into $\operatorname{SO}(n+1)$ if $n$ is even.
  This representation is not faithful because 
  $\inn(b) = \inn(-b)$ for all $b \in \S^n$,
  but we obtain a faithful representation of
  the projective quandle $(\RP^n,\ast)$.
  A faithful representation of the spherical quandle $(\S^n,\ast)$ 
  is obtained by lifting to the double covering
  $\operatorname{Pin}(n+1) \onto \operatorname{O}(n+1)$,
  see \cite[\textsection XIX.4]{Lang:2002}.
  \[
  \begin{diagram}
    \node{(\S^n,\ast)}
    \arrow{s,A}
    \arrow{e,J}
    \node{\operatorname{Pin}(n+1)}
    \arrow{s,A}
    \\
    \node{(\RP^n,\ast)}
    \arrow{e,J}
    \node{\operatorname{O}(n+1)}
  \end{diagram}
  \]
\end{remark}



\subsection{Symmetric spaces}

Reflection quandles and spherical quandles have a beautiful 
common generalization: globally symmetric Riemannian manifolds.
They have been introduced and classified by \'Elie Cartan 
in the 1920s and form a classical object of Riemannian geometry.
(See Helgason \cite[\textsection IV.3]{Helgason:2001},
Loos \cite{Loos:1969}, 
Klingenberg \cite[\textsection 2.2]{Klingenberg:1995},
Lang \cite[\textsection XIII.5]{Lang:1999}.)
We briefly recall some elementary properties in order to 
characterize the quandle coverings that naturally arise in this context.%
\footnote{ In the classification of symmetric spaces one usually 
  passes to universal coverings and then concentrates on simply connected spaces.
  The observations that follow concern non-simply connected symmetric spaces,
  because we are particularly interested in the coverings themselves.  
  We will not appeal to the classification, so our remarks can be 
  considered an elementary complement to the simply connected case.}

\begin{definition}
  A \emph{symmetric space} is a smooth connected manifold $X$ equipped with 
  a Riemannian metric such that for each point $x \in X$ there exists 
  an isometry $s_x \colon X \isoto X$ that reverses every geodesic arc 
  $\gamma \colon ({\left]-\varepsilon,+\varepsilon\right[},0) \to (X,x)$,
  meaning that $s_x \circ \gamma(t) = \gamma(-t)$.  
\end{definition}

In a symmetric space every geodesic arc can be prolonged 
to a complete geodesic $\R \to X$, and the Hopf-Rinow theorem 
implies that $X$ is a complete Riemannian manifold.  
Conversely, the fact that $X$ is connected and complete
ensures that any two points $x,x' \in X$ can be joined by 
a geodesic, and so the symmetry $s_x$ is unique for each $x$.

\begin{proposition} 
  A symmetric space $X$ is an involutory quandle 
  with respect to the operation $\ast \colon X \times X \to X$
  defined by the symmetry $x \ast y = s_y(x)$.
\end{proposition}

\begin{proof}
  Axiom (Q1) follows from $s_x(x) = x$, and Axiom (Q2) from $s_x^2 = \id_X$.
  For (Q3) notice that the isometry $s_z s_y s_z$ reverses every geodesic 
  $(\R,0) \to (X,s_z(y))$, and so we conclude $s_z s_y s_z = s_{s_z(y)}$ 
  by uniqueness of the symmetry about $s_z(y)$.
\end{proof}

\begin{remark} \label{rem:ConnectedSymmetricSpaces}
  For a symmetric space $X$, topological connectedness 
  entails algebraic connectedness.  The quandle $(X,\ast)$ 
  is even \emph{strongly} connected: since any two points $x,x' \in X$ 
  can be joined by a geodesic $\gamma \colon \R \to X$ with 
  $\gamma(0) = x$ and $\gamma(1) = x'$, the symmetry about 
  $y = \gamma(\frac{1}{2})$ maps $x$ to $x'$.  In other words,
  we do not need a product of successive symmetries to go from $x$ to $x'$; 
  one step suffices.  For the quandle $(X,\ast)$ this means that $x' = x^g$ 
  for some $g \in \inn(X)$, rather than $g \in \Inn(X)$ as usual.
\end{remark}

In favourable cases a covering $p \colon \tilde{X} \onto X$ 
of symmetric spaces is also a quandle covering $(\tilde{X},\ast) \onto (X,\ast)$,
as for $\S^n \onto \RP^n$, but in general it need not be,
as illustrated by the example $\R^n \onto \T^n$ above.
For Lie groups this phenomenon is easy to understand:

\begin{example}
  Consider a Lie group $G$ with a bi-invariant Riemannian metric, for example, 
  a compact Lie group.  (See \cite[\textsection IV.6]{Helgason:2001}).
  In this case a smooth map $(\R,0) \to (G,1)$ is a geodesic
  if and only if it is a group homomorphism, and $G$ is 
  a symmetric space: the geodesic-reversing involution 
  at $1 \in G$ is just $s_1(g) = g^{-1}$, and for any other point 
  $h \in G$ we find $s_h(g) = h g^{-1} h$.  We thus recover
  the core quandle of $G$ of Example \ref{exm:CoreQuandle},
  and we deduce from Example \ref{exm:CoreCovering} that a covering 
  $p \colon \tilde{G} \onto G$ of Lie groups is a quandle
  covering if and only if $\ker(p)$ is a group of exponent $2$.
  This is actually the general condition:
\end{example}

\begin{theorem} \label{thm:SymmetricExponent2}
  Let $X$ be a symmetric space.  For every connected covering 
  $p \colon \tilde{X} \to X$ the covering space $\tilde{X}$
  carries a unique Riemannian structure such that $p$ is 
  a local isometry.  Equipped with this canonical structure, 
  $\tilde{X}$ is itself a symmetric space and $p$ is 
  a quandle homomorphism.  It is a quandle covering 
  if and only if $\Aut(p)$ is a group of exponent $2$.
\end{theorem}

The proof relies on the following observation,
which is interesting in its own right:

\begin{lemma} \label{lem:GeodesicLoop}
  Let $X$ be a homogeneous Riemannian manifold.
  Then in every homotopy class $c \in \pi_1(X,x)$ there exists a loop 
  $\gamma \colon [0,1] \to X$, with $\gamma(0) = \gamma(1) = x$,
  minimizing the arc-length of all loops in $c$.  Every such loop
  $\gamma$ is a closed geodesic, satisfying $\gamma'(0) = \gamma'(1)$, 
  so that its continuation defines a geodesic $(\R,0) \to (X,x)$ of period $1$.
  \qed
\end{lemma}

Notice that we do not consider free homotopy classes,
but homotopy classes based at $x$.  Moreover, $X$ need not 
be compact; the crucial hypothesis is homogeneity. 
For the special case of symmetric spaces, which is of interest to us here,
the conclusion $\gamma'(0) = \gamma'(1)$ can be obtained by parallel 
transport along $\gamma$, see \cite[Corollary 2.2.7]{Klingenberg:1995}.

\begin{Notes}
  \begin{proof}
    Let $p \colon (\tilde{X},\tilde{x}) \onto (X,x)$ be the universal covering 
    of the space $X$.  Since $X$ is a complete Riemannian manifold, so is $\tilde{X}$.
    Every loop $\alpha \colon [0,1] \to X$ based at $x = \alpha(0) = \alpha(1)$ lifts 
    to a path $\tilde\alpha \colon [0,1] \to \tilde{X}$ from $\tilde{x} = \tilde\alpha(0)$
    to some point $\hat{x} := \tilde\alpha(1)$.  By the Hopf-Rinow theorem there 
    exists a length-minimizing geodesic $\tilde\gamma \colon [0,1] \to \tilde{X}$
    from $\tilde\gamma(0) = \tilde{x}$ to $\tilde\gamma(1) = \hat{x}$. 
    Its projection $\gamma = p\tilde\gamma \colon [0,1] \to X$ is a geodesic loop 
    with $[\alpha] = [\gamma]$.  Moreover $\length(\alpha) = \length(\tilde\alpha) \ge
    \length(\tilde{\gamma}) = \length(\gamma)$, so that $\gamma$ minimizes 
    the arc-length of all loops in the homotopy class $c=[\alpha]$.
    
    It remains to show that $\gamma'(0) = \gamma'(1)$.
    If the tangent vectors $\gamma'(0) \ne \gamma'(1)$ formed an angle, 
    then we could obtain a shorter loop $\beta \colon [0,1] \to X$ 
    by smoothing the angle.  Of course this would move $\beta(0)$ 
    away from the base point $x$, but since $X$ is assumed homogeneous, 
    we can move it back by some isometry to ensure $\beta(0) = x$.  
    We can do this in such a way that $\beta$ is arbitrarily close 
    to $\gamma$, so that both loops are homotopic in $X$ fixing the base point.  
    But then $\length(\beta) < \length(\gamma)$ contradicts the 
    minimality of $\gamma$.  We conclude that $\gamma'(0) = \gamma'(1)$, 
    so that we obtain a closed geodesic as desired.
  \end{proof}
\end{Notes}

\begin{proof}[Proof of the theorem]
  The symmetry $s_x \colon (X,x) \to (X,x)$ acts as 
  inversion on $\pi_1^\mathrm{top}(X,x)$, which implies that
  this group is abelian.   Every connected covering 
  $p \colon (\tilde{X},\tilde{x}) \to (X,x)$ is thus galois,
  and the symmetry $s_x \colon (X,x) \to (X,x)$ lifts to a symmetry  
  $s_{\tilde{x}} \colon (\tilde{X},\tilde{x}) \to (\tilde{X},\tilde{x})$.
  This turns $\tilde{X}$ into a Riemannian symmetric space,
  and we obtain a quandle $(\tilde{X},\ast)$.
  The projection $p$ is a quandle homomorphism:
  for $a,b \in \tilde{X}$ we have $p \circ s_{b} = s_{p(b)} \circ p$,
  whence $p(a \ast b) = p \circ s_{b}(a) = s_{p(b)} \circ p(a) = p(a) \ast p(b)$.

  Any two points $a, b \in p^{-1}(x)$ are related
  by a unique deck transformation $h \in \Aut(p)$ such that $h(a) = b$, 
  and by a geodesic $\gamma \colon \R \to \tilde{X}$ with $\gamma(0) = a$ 
  and $\gamma(1) = b$ such that $\gamma|_{[0,1]}$ is length-minimizing.
  We thus have $\gamma(s) \ast \gamma(t) = \gamma(2t-s)$ for all $s,t \in \R$,
  and also $h \gamma(t) = \gamma(t+1)$ according to Lemma \ref{lem:GeodesicLoop}.
    
  If $p$ is a quandle covering, then $s_{a} = s_{b}$
  entails $h^2(a) = \gamma(2) = \gamma(0) \ast \gamma(1) = a \ast b = a$.
  This shows that the deck transformation $h^2 \colon \tilde{X} \to \tilde{X}$ 
  fixes $a$ and is thus the identity. 
  
  Conversely, if $h^2 = \id$, then
  $a \ast b = \gamma(0) \ast \gamma(1) = \gamma(2) = h^2(a) = a$.
  This implies that $s_{b} = s_{a}$, because both 
  are liftings of $s_x = p \circ s_{a} = p \circ s_{b}$ fixing $a$.
\end{proof}

\begin{remark}
  The examples of Lie groups and symmetric Riemannian manifolds
  are manifestly of a topological nature, and the quandles that emerge
  naturally are \emph{topological quandles}, analogous to topological groups.  
  It is conceivable to define the adjoint group in the topological category, 
  so that the adjoint augmentation $Q \to \Adj(Q) \to \Inn(Q)$ is continuous
  and universal in an appropriate sense.  Likewise, the theory of (algebraic 
  i.e.\ discrete) quandle coverings can be adapted to continuous quandle coverings, 
  and a topological Galois correspondence can be established.
  We postpone this generalization and consider only the algebraic aspect,
  that is, discrete quandles, in this article.
\end{remark}

\subsection{Historical remarks} \label{sub:HistoricalRemarks}

As early as 1942, M.\,Takasaki \cite{Takasaki:1942} 
introduced  the notion of ``kei'' (i.e.\ involutory quandle) 
as an abstraction of symmetric spaces, and later 
O.\,Loos \cite{Loos:1969} extensively studied symmetric spaces 
as differential manifolds with an involutory quandle structure.
Racks first appeared around 1959 under the name ``wracks'' 
in unpublished correspondence between J.H.\,Conway 
and G.C.\,Wraith (see \cite{FennRourke:1992}).
D.\,Joyce published the first comprehensive treatment 
of quandles in 1982, and also coined the name ``quandle''.
Independently, S.\,Matveev studied the equivalent notion of 
``distributive groupoid'' (which is not a groupoid in the usual sense).
Racks were rediscovered on many occasions and studied under various names:
as ``automorphic sets'' by E.\,Brieskorn \cite{Brieskorn:1988},
as ``crossed $G$-sets'' by  P.J.\,Freyd and D.N.\,Yetter \cite{FreydYetter:1989},
as ``racks'' by R.\,Fenn and C.\,Rourke \cite{FennRourke:1992},
and as ``crystals'' by L.H.\,Kauffman \cite{Kauffman:2001}.
For a detailed review see \cite{FennRourke:1992}.


\section{The category of quandle coverings} \label{sec:CoveringQuandles}

This section initiates the systematic study of quandle coverings.
They correspond vaguely to central group extensions, 
but also incorporate intrinsically non-abelian features.  
The best analogy seems to be with coverings of topological spaces.
Throughout this article we will use this analogy 
as a guiding principle wherever possible.

\subsection{The category of quandle coverings}

We have already seen that the composition of quandle coverings is 
in general not a quandle covering (see \secref{sub:TournantsDangereux}).
In order to obtain a category we have to consider coverings 
over a fixed base quandle:

\begin{definition}
  Let $p \colon \tilde{Q} \to Q$ and $\hat{p} \colon \hat{Q} \to Q$ 
  be two quandle coverings.  A \emph{covering morphism} from $p$ to $\hat{p}$ 
  (over $Q$) is a quandle homomorphism $\phi \colon \tilde{Q} \to \hat{Q}$
  such that $p = \hat{p} \circ \phi$.
  \[
  \begin{diagram}
    \node{\tilde{Q}}
    \arrow[2]{e,t}{\phi}
    \arrow{se,b}{\smash{p}}
    \node[2]{\hat{Q}}
    \arrow{sw,b}{\smash{\hat{p}}}
    \\
    \node[2]{Q}
  \end{diagram}
  \]
\end{definition}

\begin{proposition} \label{prop:CoveringMorphismEquivariance}
  A map $\phi \colon \tilde{Q} \to \hat{Q}$ with 
  $p = \hat{p} \circ \phi$ is a covering morphism if and only if 
  $\phi$ is equivariant with respect to $\Adj(Q)$, or
  equivalently, its subgroup $\Adj(Q)^\circ$.
\end{proposition}

\begin{proof}
  Consider $\tilde{a},\tilde{b} \in \tilde{Q}$
  and $b = p(\tilde{b}) = \hat{p}\phi(\tilde{b})$.
  Since both $p$ and $\hat{p}$ are coverings, we have 
  on the one hand $\phi(\tilde{a} \ast \tilde{b}) = \phi(\tilde{a}^{\adj(b)})$
  and on the other hand $\phi(\tilde{a}) \ast \phi(\tilde{b}) = \phi(\tilde{a})^{\adj(b)}$.
  This proves the desired equivalence.
  It suffices to assume equivariance under 
  the subgroup $\Adj(Q)^\circ$, by replacing $\adj(b)$ 
  with $\adj(a)^{-1}\adj(b) \in \Adj(Q)^\circ$
  where $a = p(\tilde{a})$.
\end{proof}


\begin{proposition}
  Given a quandle $Q$, the coverings  $p \colon \tilde{Q} \to Q$
  together with their covering morphisms form a category, called 
  the category of \emph{coverings over $Q$}, denoted $\Coverings(Q)$.
\end{proposition}

\begin{proof}
  The only point to verify is that, given three coverings 
  $p_i \colon \tilde{Q}_i \to Q$ with $i=1,2,3$, the composition of 
  two covering morphisms $\phi_1 \colon \tilde{Q}_1 \to \tilde{Q}_2$ 
  and $\phi_2 \colon \tilde{Q}_2 \to \tilde{Q}_3$ is again
  a covering morphism.  We already know that $\Quandles$ is a category,
  so $\phi = \phi_2 \circ \phi_1 \colon \tilde{Q}_1 \to \tilde{Q}_3$
  is a quandle homomorphism.  Moreover, $p_3 \circ \phi = p_3 \circ \phi_2 \circ \phi_1
  = p_2 \circ \phi_1 = p_1$. 
\end{proof}

\begin{remark}
  Every surjective covering morphism $\phi \colon \tilde{Q} \onto \hat{Q}$ 
  is itself a quandle covering: if $\phi(\tilde{x}) = \phi(\tilde{y})$ then
  $p(\tilde{x}) = \hat{p}\phi(\tilde{x}) = \hat{p}\phi(\tilde{y}) = p(\tilde{y})$ 
  and so $\inn(\tilde{x}) = \inn(\tilde{y})$.
\end{remark}

\begin{definition}
  For a quandle covering $p \colon \tilde{Q} \to Q$
  we define $\Aut(p)$ to be the group of covering automorphisms of $p$, 
  also called the \emph{group of deck transformations} of the covering $p$.
  We will adopt the convention that deck transformations of $p$ 
  act on the left, which means that their composition $\phi\psi$ 
  is defined by $(\phi\psi)(\tilde{q}) = \phi(\psi(\tilde{q}))$ 
  for all $\tilde{q} \in \tilde{Q}$.  
\end{definition}

We let $\Aut(p)$ act on the left because this is 
the most convenient (and traditional) way 
to denote two commuting actions:

\begin{proposition} \label{prop:CommutingActions}
  Given a quandle covering $p \colon \tilde{Q} \to Q$,
  two groups naturally act on the covering quandle $\tilde{Q}$:
  the group of deck transformations $\Aut(p)$ acts on the left 
  while the group of inner automorphisms $\Inn(\tilde{Q})$ 
  acts on the right.  Both actions commute.
\end{proposition}

\begin{proof}
  Consider $\phi \in \Aut(p)$ and $\tilde{x},\tilde{y} \in \tilde{Q}$.
  Then $\phi(\tilde{x} \ast \tilde{y}) 
  = \phi(\tilde{x}) \ast \phi(\tilde{y})
  = \phi(\tilde{x}) \ast \tilde{y}$,
  which means that $\phi$ and $\inn(\tilde{y})$ commute.  
  Since the group $\Inn(\tilde{Q})$ is generated by $\inn(\tilde{Q})$,
  this proves that the actions of $\Aut(p)$ and $\Inn(\tilde{Q})$ commute.
\end{proof}

\subsection{Pointed quandles and coverings}

As in the case of topological spaces, we have to choose base points
in order to obtain uniqueness properties of coverings.

\begin{definition} \label{def:PointedQuandles}
  A \emph{pointed quandle} $(Q,q)$ is a quandle $Q$ with 
  a specified base point $q \in Q$.  A homomorphism (resp.\ covering) 
  $\phi \colon (Q,q) \to (Q',q')$ between pointed quandles 
  is a quandle homomorphism (resp.\ covering)
  $\phi \colon Q \to Q'$ such that $\phi(q) = q'$.
  Pointed quandles and their homomorphisms form a category,
  denoted $\Quandles_*$.  Likewise, coverings
  $p \colon (\tilde{Q},\tilde{q}) \to (Q,q)$ over a fixed
  base quandle $(Q,q)$ form a category, denoted $\Coverings(Q,q)$.
\end{definition}

\begin{definition}
  Let $f \colon (X,x) \to (Q,q)$ and $p \colon (\tilde{Q},\tilde{q}) \to (Q,q)$
  be homomorphisms of pointed quandles.  A \emph{lifting} of $f$ over $p$ is 
  a quandle homomorphism $\tilde{f} \colon (X,x) \to (\tilde{Q},\tilde{q})$ 
  such that $p \circ \tilde{f} = f$.
  \[
  \begin{diagram}
    \node[3]{(\tilde{Q},\tilde{q})}
    \arrow{s,r}{p}
    \\
    \node{(X,x)}
    \arrow[2]{e,t}{f}
    \arrow{ene,l,..}{\tilde{f}}
    \node[2]{(Q,q)}
  \end{diagram}
  \]
\end{definition}

\begin{proposition}[lifting uniqueness] \label{prop:LiftingUniqueness}
  Let $f \colon (X,x) \to (Q,q)$ be a quandle homomorphism, and 
  let $p \colon (\tilde{Q},\tilde{q}) \to (Q,q)$ be a quandle covering.
  Then any two liftings $\tilde{f}_1, \tilde{f}_2 \colon (X,x) \to (\tilde{Q},\tilde{q})$ 
  of $f$ over $p$ coincide on the component of $x$ in $X$.  In particular, 
  if $X$ is connected, then $f$ admits at most one lifting over $p$. 
\end{proposition}

\begin{proof}
  The quandle homomorphism $f$ induces a group homomorphism
  $h \colon \Adj(X) \to \Adj(Q)$.  Since $p$ is a covering,
  the group $\Adj(Q)$ acts on $\tilde{Q}$, and so does $\Adj(X)$ via $h$.
  In this way, all the maps in the above triangle are equivariant 
  with respect to the action of $\Adj(X)$.  If $\tilde{f}_1$ and $\tilde{f}_2$
  coincide on one point $x$, they coincide on its entire orbit,
  which is precisely the connected component of $x$ in $X$.
\end{proof}

\begin{corollary} \label{cor:CoveringMorphismUniqueness}
  Between a connected covering $p \colon (\tilde{Q},\tilde{q}) \to (Q,q)$ 
  and an arbitrary covering $\hat{p} \colon (\hat{Q},\hat{q}) \to (Q,q)$
  there can be at most one covering morphism 
  $\phi \colon (\tilde{Q},\tilde{q}) \to (\hat{Q},\hat{q})$.
\end{corollary}

\begin{proof}
  The equation $p = \hat{p} \circ \phi$ means 
  that $\phi$ is a lifting of $p$ over $\hat{p}$.
\end{proof}

\begin{corollary} \label{cor:FreeAction}
  Let $p \colon \tilde{Q} \to Q$ be a quandle covering.
  If $\tilde{Q}$ is connected, then the group $\Aut(p)$ 
  of deck transformations acts freely on each fibre.
\end{corollary}

\begin{proof}
  Choose a base point $q \in Q$ and consider the fibre $F = p^{-1}(q)$.
  Every deck transformation $\phi \in \Aut(p)$ satisfies $\phi(F) = F$,
  and so $\Aut(p)$ acts on the set $F$.  If $\phi$ fixes a point 
  $\tilde{q} \in F$, then $\phi = \id$ by the previous corollary.
\end{proof}

\subsection{Galois coverings}

As for topological coverings, the galois case is most prominent:

\begin{definition}
  A covering $p \colon \tilde{Q} \to Q$ is said to be \emph{galois}
  if $\tilde{Q}$ is connected and $\Aut(p)$ acts transitively on each fibre.  
  (It necessarily acts freely by the previous corollary.)
\end{definition}

Numerous examples are provided by central group extensions 
(Remark \ref{rem:CentralExtensions} and Example \ref{exm:LinearGroups1}) 
and coverings of symmetric spaces (Examples \ref {exm:TorusCovering2} 
and \ref{exm:ProjectiveQuandle}, and Theorem \ref{thm:SymmetricExponent2}).

\begin{remark}
  Every galois covering $p \colon \tilde{Q} \to Q$ comes with
  the natural action $\Lambda \curvearrowright \tilde{Q}$ of the 
  deck transformation group $\Lambda = \Aut(p)$ satisfying 
  the following two axioms:
  \begin{enumerate} 
  \item[(E1)]
    $(\lambda\tilde{x}) \ast \tilde{y} = \lambda (\tilde{x}\ast\tilde{y})$
    and $\tilde{x} \ast (\lambda\tilde{y}) = \tilde{x}\ast\tilde{y}$
    for all $\tilde{x},\tilde{y}\in\tilde{Q}$ and $\lambda\in\Lambda$.
  \item[(E2)]
    $\Lambda$ acts freely and transitively on each fibre $p^{-1}(x)$.
  \end{enumerate}
  
  Axiom (E1) says that $\Lambda$ acts by automorphisms and 
  the left action of $\Lambda$ commutes with the right action of 
  $\Inn(\tilde{Q})$, cf.\ Proposition \ref{prop:CommutingActions}.
  We denote such an action simply by $\Lambda \curvearrowright \tilde{Q}$.
  In this situation the quotient $Q := \Lambda\backslash\tilde{Q}$ 
  carries a unique quandle structure that turns the projection 
  $p \colon \tilde{Q} \to Q$ into a quandle covering.
  Axiom (E2) then says that $p \colon \tilde{Q} \to Q$ 
  is a \emph{principal} $\Lambda$-covering, in the sense 
  that each fibre is a principal $\Lambda$-set.
\end{remark}

\subsection{Quandle extensions} \label{sub:QuandleExtensions}

The freeness expressed in (E2) relies on the connectedness 
of $\tilde{Q}$.  As an extreme counter-example, consider the trivial 
covering $p \colon \tilde{Q} = Q \times F \to Q$ where $Q$ is 
a connected quandle and $F$ is a set with at least three elements.  
Here the deck transformation group $\Aut(p) = \Sym(F)$ 
is too large: it acts transitively but not freely.

If the covering quandle $\tilde{Q}$ is non-connected,
we can nevertheless salvage the above properties by passing
from the group $\Aut(p)$ to a subgroup $\Lambda$ that satisfies (E2).
We are thus led to the concept of a principal $\Lambda$-covering.
Motivated by the terminology used in group theory, we will call 
this a quandle extension:

\begin{definition} \label{def:QuandleExtension}
  An \emph{extension} $E \colon \Lambda \curvearrowright \tilde{Q} \lto[p] Q$
  of a quandle $Q$ by a group $\Lambda$ consists of a surjective 
  quandle homomorphism $p \colon \tilde{Q} \to Q$ and 
  a group action $\Lambda \curvearrowright \tilde{Q}$ 
  satisfying the above axioms (E1) and (E2).
  This can also be called a \emph{principal $\Lambda$-covering} of $Q$.
\end{definition}

Quandle extensions are intermediate between galois coverings 
and general coverings: 

\begin{proposition} \label{prop:CoveringExtensions}
  In every extension $E \colon \Lambda \curvearrowright \tilde{Q} \lto[p] Q$
  the projection $p \colon \tilde{Q} \to Q$ is a quandle covering.
  It is a galois covering if and only if $\tilde{Q}$ is connected.

  Conversely, every galois covering $p \colon \tilde{Q} \to Q$
  defines an extension of $Q$, with the group $\Lambda = \Aut(p)$ 
  acting naturally on $\tilde{Q}$ by deck transformations.
  \qed
\end{proposition}

We have already seen quandle extensions in the general Examples
\ref{exm:ConjCovering}, \ref{exm:CoreCovering}, \ref{exm:AlexCovering},
and the more concrete Examples \ref{exm:TrivialCovering}, \ref{exm:LinearGroups1},
\ref{exm:TorusCovering2}, \ref{exm:ProjectiveQuandle}.
Here is another natural construction, which essentially
goes back to Joyce \cite[\textsection{7}]{Joyce:1982}
and will be proven universal in \secref{sub:UniversalCovering}.

\begin{example} \label{exm:AlexanderQuotient}
  As in Example \ref{exm:AlexanderQuandle} we consider 
  a group $G$ with automorphism $T \colon G \isoto G$
  and the associated Alexander quandle $Q = \Alex(G,T)$.
  Suppose that $H \subset G$ is a subgroup such that $T|_H = \id_H$.
  Then $H \times G \to G$, $(h,g) \mapsto hg$ defines
  a free action of $H$ on the quandle $Q$ satisfying axiom (E1) above.  
  As a consequence, the quotient set $\bar{Q} = H \backslash G$ 
  carries a unique quandle structure such that the projection 
  $p \colon Q \onto \bar{Q}$ is a quandle homomorphism,
  and $H \curvearrowright Q \onto \bar{Q}$ is a quandle extension.
\end{example}

Coverings of $Q$ form a category, which provides us with 
a natural notion of isomorphism, i.e.\ equivalence of coverings.
Here is the appropriate notion for extensions:

\begin{definition} \label{def:EquivalentExtensions}
  Let $Q$ be a quandle and let $\Lambda$ be a group.
  An \emph{equivalence}, or \emph{isomorphism}, between extensions 
  $E_1 \colon \Lambda \curvearrowright Q_1 \lto[p_1] Q$ and 
  $E_2 \colon \Lambda \curvearrowright Q_2 \lto[p_2] Q$ is 
  a quandle isomorphism $\phi \colon Q_1 \isoto Q_2$ that respects projections, 
  $p_1 = p_2 \phi$, and is equivariant, $\phi \lambda = \lambda \phi$ 
  for all $\lambda \in \Lambda$.  We denote by $\Ext(Q,\Lambda)$ 
  the set of equivalence classes of extensions of $Q$ by $\Lambda$.
\end{definition}

On could also define the seemingly weaker notion of 
\emph{homomorphism} between extensions $E_1$ and $E_2$ as a quandle 
homomorphism $\phi \colon Q_1 \to Q_2$ that respects projections 
and is $\Lambda$-equivariant.  This leads to the following observation,
which is a variant of the well-known Five Lemma for short exact sequences 
in abelian categories (see \cite[\textsection VIII.4]{MacLane:1998}).

\begin{proposition}
  Every homomorphism $\phi \colon Q_1 \to Q_2$ between two quandle 
  extensions $E_1 \colon \Lambda \curvearrowright Q_1 \lto[p_1] Q$ 
  and $E_2 \colon \Lambda \curvearrowright Q_2 \lto[p_2] Q$
  is an isomorphism of extensions.
  \qed
\end{proposition}

The proof is a straightforward diagram chase, and will be omitted.

\begin{Notes}
  \begin{proof}
    Let us first show that $\phi \colon Q_1 \to Q_2$ is surjective.
    Given $q_2 \in Q_2$ consider $q := p_2(q_2)$: there exists $q'_1$ 
    such that $p_1(q'_1) = q$ because $p_1$ is surjective.  Now both 
    $q'_2 = \phi(q'_1)$ and $q_2$ satisfy $p_2(q_2) = p_2(q'_2) = q$ 
    because $p_1 = p_2 \circ \phi$.  There exists $\lambda \in \Lambda$ 
    such that $\lambda q'_2 = q_2$, because $\Lambda$ acts transitively 
    on the fibre $p_2^{-1}(q)$.  We set $q_1 := \lambda q'_1$ and find 
    $\phi(q_1) = q_2$ by equivariance.
    
    For injectivity, consider $q_1, q'_1 \in Q_1$ such that $\phi(q_1) = \phi(q'_1)$.
    Since $p_1 = p_2 \circ \phi$, we know that $q_1$ and $q'_1$ belong
    to the same fibre $p_1^{-1}(q)$ over $q := p_1(q_1) = p_1(q'_1)$.
    There exists $\lambda \in \Lambda$ such that $\lambda q'_1 = q_1$,
    because $\Lambda$ acts transitively on $p_1^{-1}(q)$.
    By equivariance we find $\phi(q_1) = \lambda \phi(q'_1)$,
    whence $\lambda = 1$ because $\Lambda$ acts freely on $p_2^{-1}(q)$.
    We conclude that $q_1 = q'_1$, hence $\phi$ is injective.
  \end{proof}
\end{Notes}

\subsection{Pull-backs}

Given quandle homomorphisms $p \colon \tilde{Q} \to Q$ and 
$f \colon X \to Q$ we construct their \emph{pull-back},
or \emph{fibred product} $\tilde{X} = X \times_Q \tilde{Q}$ as follows:
\[
\begin{CD}
  \tilde{X} @>{\tilde{f}}>> \tilde{Q} \\
  @V{\tilde{p}}VV @VV{p}V \\
  X @>{f}>> Q
\end{CD}
\]

The set $\tilde{X} := \{ (x,\tilde{a}) \in  X \times \tilde{Q} \mid f(x) = p(\tilde{a}) \}$
can be equipped with a quandle operation 
$(x,\tilde{a}) \ast (y,\tilde{b}) := (x \ast y, \tilde{a} \ast \tilde{b})$
such that the projections 
$\tilde{p}(x,\tilde{a}) = x$ and 
$\tilde{f}(x,\tilde{a}) = \tilde{a}$ are quandle homomorphisms and make the above diagram commute.  
The triple $(\tilde{X},\tilde{p},\tilde{f})$ is universal in the usual sense 
that any other candidate uniquely factors through it, and 
this property characterizes it up to unique isomorphism.

The quandle homomorphism $f^* p := \tilde{p} \colon \tilde{X} \to X$
is called the \emph{pull-back} of $p$ along $f$.

\begin{proposition}
  If $p$ is a covering, then its pull-back $f^* p$ is again a covering.
  Thus every quandle homomorphism $f \colon X \to Q$ induces 
  a covariant functor $f^* \colon \Coverings(Q) \to \Coverings(X)$
  by sending each covering $p \colon \tilde{Q} \to Q$
  to its pull-back $f^* p \colon \tilde{X} \to X$,
  and every morphism between coverings to the induced
  morphism between their pull-backs.
\end{proposition}

\begin{proof}
  Suppose that $p \colon \tilde{Q} \to Q$ is a covering,
  that is, $p$ is surjective and $p(\tilde{a}) = p(\tilde{b})$ 
  implies $\inn(\tilde{a}) = \inn(\tilde{b})$.
  Then $\tilde{p} \colon \tilde{X} \to X$ is surjective, and for all 
  $\tilde{x} = (x,\tilde{a})$ and $\tilde{y} = (y,\tilde{b})$ 
  the equality $\tilde{p}(\tilde{x}) = \tilde{p}(\tilde{y})$ entails $x = y$ 
  as well as $p(\tilde{a}) = f(x) = f(y) = p(\tilde{b})$.
  These in turn imply that $\inn(\tilde{x}) = \inn(\tilde{y})$, as claimed.
  This construction is natural with respect
  to covering morphisms, whence $f^*$ is a functor.
\end{proof}

For extensions $\Lambda \curvearrowright \tilde{Q} \onto Q$
we record the following observations:

\begin{proposition}[functoriality in $Q$]
  The pull-back of an extension $E \colon \Lambda \curvearrowright \tilde{Q} \lto[p] Q$
  along a quandle homomorphism $f \colon X \to Q$ inherits a natural $\Lambda$-action
  and defines an extension $f^* E \colon \Lambda \curvearrowright \tilde{X} \lto[f^* p] X$.
  We thus obtain a natural map $f^* \colon \Ext(Q,\Lambda) \to \Ext(X,\Lambda)$.
\end{proposition}

\begin{proof}
  The action on $\tilde{X}$ is given by $\lambda(x,\tilde{a}) = (x,\lambda\tilde{a})$
  for $\lambda \in \Lambda$.  Axioms (E1) and (E2) carry over from $\tilde{Q}$
  to $\tilde{X}$, so that $f^* E$ is an extension, as claimed.
\end{proof}

\begin{proposition}[functoriality in $\Lambda$] 
  Every group homomorphism $h \colon \Lambda \to \Lambda'$ induces 
  a natural map on extensions, $h_* \colon \Ext(Q,\Lambda) \to \Ext(Q,\Lambda')$.
\end{proposition}

\begin{proof}
  Given an extension $E \colon \Lambda \curvearrowright \tilde{Q} \lto[p] Q$,
  the induced extension $h_* E$ is defined as the product $\Lambda' \times \tilde{Q}$ 
  modulo the relation $(\lambda',\lambda\tilde{a}) \sim (\lambda'h(\lambda),\tilde{a})$
  for $\lambda \in \Lambda$.  The quotient $\hat{Q}$ inherits the quandle structure
  $[\lambda',a] \ast [\lambda'',b] = [\lambda', a \ast b]$, and the extension 
  $h_* E \colon \Lambda' \curvearrowright \hat{Q} \lto[\hat{p}] Q$ is defined by 
  the projection $\hat{p}[\lambda',a] = p(a)$ and the action 
  $\lambda' [\lambda'',a] = [\lambda'\lambda'',a]$.
  This construction is well-defined on isomorphism classes of extensions, 
  so that we obtain $h_* \colon \Ext(Q,\Lambda) \to \Ext(Q,\Lambda')$ as desired.
\end{proof}

The preceding propositions can be restated as saying that
$\Ext(Q,\Lambda)$ is a contravariant functor in $Q$ and
a covariant functor in $\Lambda$.  In general $\Ext(Q,\Lambda)$
is only a set, with the class of the trivial extension 
as zero element.  We obtain a group structure if $\Lambda$ is abelian:

\begin{proposition}[module structure]
  If $\Lambda$ is an abelian group, or more generally a module over some ring $R$, 
  then $\Ext(Q,\Lambda)$ carries a natural $R$-module structure, and the pull-back
  $f^* \colon \Ext(Q,\Lambda) \to \Ext(X,\Lambda)$ is a homomorphism of $R$-modules.
\end{proposition}

\begin{proof}
  The group $\Lambda$ is abelian if and only if its multiplication
  $\mu \colon \Lambda \times \Lambda \to \Lambda$
  is a group homomorphism.  In this case we obtain a 
  binary operation on $\Ext(Q,\Lambda)$ as follows:
  \[
  \otimes \colon \Ext(Q,\Lambda) \times \Ext(Q,\Lambda) 
  \lto[P] \Ext(Q,\Lambda\times\Lambda)
  \lto[\mu_*] \Ext(Q,\Lambda)
  \]
  Here $P$ is the fibred product and $\mu_*$ is the induced map as above.
  More explicitly, given two extensions $E_1 \colon \Lambda \curvearrowright Q_1 \lto[p_1] Q$ 
  and $E_2 \colon \Lambda \curvearrowright Q_2 \lto[p_2] Q$, their composition
  $E_3 = E_1 \otimes E_2$ is the fibred product $Q_1 \times_Q Q_2$
  modulo the relation $(\lambda a_1, a_2) \sim (a_1,\lambda a_2)$ for $\lambda \in \Lambda$.
  The quotient $Q_3$ inherits the quandle structure
  $[a_1,a_2] \ast [b_1,b_2] = [a_1 \ast b_1, a_2 \ast b_2]$, and the extension 
  $E_3 \colon \Lambda \curvearrowright Q_3 \lto[p_3] Q$ is defined by 
  the projection $p_3[a_1,a_2] = p_1(a_1) = p_2(a_2)$ and 
  the action $\lambda [a_1,a_2] = [\lambda a_1,a_2] = [a_1,\lambda a_2]$.

  The composition is well-defined and associative on isomorphism classes 
  of extensions.  The neutral element is given by the trivial extension 
  $E_0 \colon \Lambda \curvearrowright \Lambda \times Q \lto[\pr_2] Q$.
  The inverse of $E_1$ is obtained by replacing the action of $\Lambda$ 
  with the inverse action via $\lambda \mapsto \lambda^{-1}$.
  The details are easily verified and will be omitted.
\end{proof}


\section{Classification of connected coverings} \label{sec:ConnectedCoverings}

In order to avoid clumsy notation, we will first classify 
connected coverings.  The passage to arbitrary coverings 
over a connected base quandle is then straightforward, 
and will be treated in Section \ref{sec:NonConnectedCoverings}.
Assuming that the base quandle is connected is technically easier 
and corresponds most closely to our model, the Galois 
correspondence for coverings over a connected topological space.
The non-connected case will be treated in Section \ref{sec:NonConnectedBase}.

\subsection{Explicit construction of universal covering quandles} \label{sub:UniversalCovering}

Our first task is to ensure the existence of a universal covering quandle.
As usual, universality is defined as follows:

\begin{definition} \label{def:UniversalCovering}
  A pointed quandle covering $p \colon (\tilde{Q},\tilde{q}) \to (Q,q)$
  is \emph{universal} if for each covering $\hat{p} \colon (\hat{Q},\hat{q}) \to (Q,q)$
  there exists a unique covering morphism $\phi \colon (\tilde{Q},\tilde{q}) \to(\hat{Q},\hat{q})$.
  In other words, a universal covering is an initial object in the category $\Coverings(Q,q)$.
  Two universal coverings of $(Q,q)$ are isomorphic by a unique isomorphism,
  so that we can unambiguously speak of \emph{the} universal covering of $(Q,q)$.
\end{definition}

The following explicit construction has been inspired 
by \cite[Lemma 25]{Eisermann:2003}.

\begin{lemma} \label{lem:UniversalCovering}
  Consider a connected quandle $Q$ with base point $q \in Q$.
  Recall that the commutator subgroup $\Adj(Q)'$ is the kernel 
  of the group homomorphism $\varepsilon \colon \Adj(Q) \to \Z$ 
  with $\varepsilon(\adj(Q)) = 1$.  We define
  \[
  \tilde{Q} := \bigl\{\; (a,g) \in Q \times \Adj(Q)' \;|\; a=q^g \;\bigr\}, 
  \quad \tilde{q} := (q,1) . 
  \]
  The set $\tilde{Q}$ becomes a connected quandle with the operations
  \begin{align*}
    (a,g) \ast (b,h) & := \bigl(\; a \ast b,\; g \cdot \adj(a)^{-1} \cdot \adj(b) \;\bigr),  
    \\
    (a,g) \tsa (b,h) & := \bigl(\; a \tsa b,\; g \cdot \adj(a) \cdot \adj(b)^{-1} \;\bigr).
  \end{align*}
  The quandle $\tilde{Q}$ comes with a natural augmentation
  $\tilde{Q} \lto[\rho] \Adj(Q) \lto[\alpha] \Inn(\tilde{Q})$,
  where $\rho(b,h) = \adj(b)$ and $\alpha$ is defined by the action
  \[
  \tilde{Q} \times \Adj(Q) \to \tilde{Q} \quad\text{with}\quad
  (a,g)^h := \bigl(\; a^h,\; \adj(q)^{-\varepsilon(h)} \cdot g h \;\bigr) .
  \]
  By construction, the subgroup $\Adj(Q)' = \ker(\varepsilon)$ 
  acts freely and transitively on $\tilde{Q}$.
  The canonical projection $p \colon \tilde{Q} \to Q$ 
  given by $p(a,g)=a$ is a surjective quandle homomorphism, 
  and equivariant with respect to the action of $\Adj(Q)$.
\end{lemma}

\begin{proof}
  Since $Q$ is connected, we have $\adj(a)^{-1} \adj(b) \in \Adj(Q)'$,
  which ensures that the operations $\ast$ and $\tsa$ are 
  well-defined.  The first quandle axiom (Q1) is obvious:
  \[
  (a,g) \ast (a,g) = \bigl(\; a \ast a,\; g \cdot \adj(a)^{-1} \cdot \adj(a) \;\bigr) = (a,g) .
  \]
  
  The second axiom (Q2) follows using $\adj(a \ast b) = \adj(b)^{-1} \adj(a) \adj(b)$:
  \[
  \bigl( (a,g) \ast (b,h) \bigr) \tsa (b,h) = 
  \bigl(\; a ,\; g \cdot \adj(a)^{-1} \cdot \adj(b) 
  \cdot \adj(a \ast b) \cdot \adj(b)^{-1} \;\bigr) = (a,g) .
  \]
  
  For the third axiom (Q3) notice that each $(a,g) \in \tilde{Q}$ satisfies 
  $a = q^g$, which entails $\adj(a) = g^{-1} \cdot \adj(q) \cdot g$.
  The quandle operations can thus be reformulated as
  \begin{align*}
    (a,g) \ast (b,h) & = \bigl(\; a \ast b,\; \adj(q)^{-1} \cdot g \cdot \adj(b) \;\bigr),  
    \\
    (a,g) \tsa (b,h) & = \bigl(\; a \tsa b,\; \adj(q) \cdot g \cdot \adj(b)^{-1} \;\bigr).
  \end{align*}
  This implies self-distributivity, because
  \begin{align*}
    \bigl( (a,g) \ast (b,h) \bigr) \ast (c,k) 
    & = \bigl(\; (a \ast b) \ast c ,\; \adj(q)^{-2} g \adj(b) \adj(c) \;\bigr)  \qquad\text{equals}
    \\
    \bigl( (a,g) \ast (c,k) \bigr) \ast \bigl( (b,h) \ast (c,k) \bigr) 
    & = \bigl(\; (a \ast c) \ast (b \ast c) ,\; \adj(q)^{-2} g \adj(c) \adj(b \ast c) \;\bigr) .
  \end{align*}

  The projection $p \colon \tilde{Q} \to Q$, $p(a,g) = a$,
  is a quandle homomorphism, which implies that 
  $\rho = \adj \circ p \colon \tilde{Q} \to \Adj(Q)$
  is a representation.  Moreover, the action $\alpha$ satisfies 
  $(a,g) \ast (b,h) = (a,g)^{\adj(b)}$, so that $(\rho,\alpha)$ 
  is an augmentation.  Since $\adj(Q)$ generates the group $\Adj(Q)$,
  this also shows that $\Adj(Q)$ acts on $\tilde{Q}$ by inner automorphisms, 
  and that $p$ is equivariant with respect to the action of $\Adj(Q)$.
  Under this action, the subgroup $\Adj(Q)'$ acts freely
  and transitively on $\tilde{Q}$, which shows 
  that $\tilde{Q}$ is connected. 
\end{proof}

The reader will notice a close resemblance with the construction 
of the universal covering for a connected topological space. 
In order to construct $\tilde{Q}$ from $Q$, we keep track not only
of the points $a \in Q$ but also the paths $g \in \Adj(Q)'$ 
leading from our base point $q$ to the point $a$ in question.
Forgetting the extra information projects back to $Q$,
while keeping it defines the universal covering $\tilde{Q} \onto Q$,
as we shall now prove:

\begin{theorem} \label{thm:UniversalCovering}
  Let $Q$ be a connected quandle with base point $q \in Q$ and let 
  $(\tilde{Q},\tilde{q})$ be defined as in Lemma \ref{lem:UniversalCovering} above.
  Then the canonical projection $p \colon (\tilde{Q},\tilde{q}) \to (Q,q)$ 
  is the universal quandle covering of $(Q,q)$.
\end{theorem}

\begin{proof}
  It is clear from its construction that 
  $p \colon (\tilde{Q},\tilde{q}) \to (Q,q)$ 
  is a covering.  We want to show that for every other 
  covering $\hat{p} \colon (\hat{Q},\hat{q}) \to (Q,q)$ 
  there exists a unique quandle homomorphism
  $\phi \colon (\tilde{Q},\tilde{q}) \to (\hat{Q},\hat{q})$
  with $\hat{p} \circ \phi = p$.  Uniqueness is clear 
  from Corollary \ref{cor:CoveringMorphismUniqueness},
  the crucial point is thus to show existence.

  We recall from Remark \ref{rem:CoveringAction} 
  that every covering $\hat{p} \colon \hat{Q} \to Q$
  induces an action of $\Adj(Q)$ on $\hat{Q}$ by inner automorphisms,
  and that $\hat{p}$ is equivariant with respect to this action.
  For our covering $p \colon \tilde{Q} \to Q$ this action has been
  made explicit in the preceding Lemma \ref{lem:UniversalCovering}.

  We define $\phi \colon (\tilde{Q},\tilde{q}) \to (\hat{Q},\hat{q})$
  by $\phi(a,g) = \hat{q}^g$.  This is an equivariant map with respect to $\Adj(Q)'$.
  Both maps $\hat{p} \phi$ and $p$ are thus equivariant and coincide in 
  $\tilde{q} = (q,1)$.  Since $\tilde{Q}$ is connected we conclude $\hat{p} \phi = p$.
  Proposition \ref{prop:CoveringMorphismEquivariance} now
  shows that $\phi$ is a quandle homomorphism, and hence a 
  covering morphism from $p$ to $\hat{p}$ as desired.
\end{proof}

\begin{remark}
  In Lemma \ref{lem:UniversalCovering}, all 
  the information of $(a,g) \in \tilde{Q}$ is contained 
  in the second coordinate $g$, so we could just as well dispense 
  with the first coordinate $a = q^g$.  This means that we consider 
  the group $G = \Adj(Q)'$ equipped with quandle operations 
  \[
  g \ast h = x^{-1}gh^{-1}xh 
  \quad\text{and}\quad 
  g \tsa h = xgh^{-1}x^{-1}h,
  \]
  where $x = \adj(q)$.  This is the (non-abelian) Alexander quandle 
  $\Alex(G,T)$ with automorphism $T \colon G \isoto G$ given by 
  $g \mapsto x^{-1}gx$.  These formulae already appear 
  in the work of Joyce \cite[\textsection{7}]{Joyce:1982}
  on the representation theory of homogeneous quandles.
  There the natural choice is $G = \Aut(Q)$, whereas 
  the universal covering requires $G = \Adj(Q)'$.
  
  The notation proposed in the preceding lemma emphasizes 
  the interpretation of $\tilde{Q}$ as a path fibration, 
  where $(a,g)$ designates a path $g$ from $q$ to the endpoint $a$.
  This extra information of base points will become necessary when 
  we consider quandles with more than one connected component,
  see Lemma \ref{lem:MultiCoveringQuandle} below.
\end{remark}

\subsection{Fundamental group of a quandle} \label{sub:FundamentalGroup}

As announced in the introduction, once we have understood the universal 
covering $p \colon (\tilde{Q},\tilde{q}) \to (Q,q)$ of a quandle $(Q,q)$,
we can define the fundamental group $\pi_1(Q,q)$ as the group 
$\Aut(p)$ of deck transformations: 

\begin{definition} \label{def:FundamentalGroup}
  We call $\pi_1(Q,q) = \{ g \in \Adj(Q)' \mid q^g = q \}$ 
  the \emph{fundamental group} of the quandle $Q$ based at $q \in Q$. 
\end{definition}

\begin{proposition} \label{prop:CanonicalIdentification}
  For the universal covering $p \colon (\tilde{Q},\tilde{q}) \to (Q,q)$ as above, 
  we obtain a canonical group isomorphism $\phi \colon \pi_1(Q,q) \isoto \Aut(p)$
  from the left action $\pi_1(Q,q) \times \tilde{Q} \to \tilde{Q}$
  defined by $h \cdot (a,g) = (a,hg)$.
\end{proposition}

\begin{proof}
  The action is well-defined and induces
  an injective group homomorphism $\pi_1(Q,q) \to \Aut(\tilde{Q})$.
  By construction it respects the projection $p \colon \tilde{Q} \to Q$, 
  so we obtain $\phi \colon \pi_1(Q,q) \to \Aut(p)$.
  The action of $\pi_1(Q,q)$ is free and transitive 
  on the fibre $p^{-1}(q) = \{ (q,g) \mid q^g = q \}$.
  Given a covering automorphism $\alpha \in \Aut(p)$
  there exists thus a unique element $h \in \pi_1(Q,q)$
  with $\alpha(\tilde{q}) = h \cdot \tilde{q}$.
  This means that $\alpha = \phi(h)$, because $\tilde{Q}$
  is connected (see Corollary \ref{cor:CoveringMorphismUniqueness}).
  This proves that $\phi$ is also surjective.
\end{proof}

\begin{proposition}[functoriality]
  Every quandle homomorphism $f \colon (X,x) \to (Y,y)$ induces 
  a homomorphism $f_* \colon \pi_1(X,x) \to \pi_1(Y,y)$
  of fundamental groups.  We thus obtain a functor
  $\pi_1 \colon \Quandles_* \to \Groups$ from the
  category of pointed quandles to the category of groups.
\end{proposition}

\begin{proof}
  Every quandle homomorphism $f \colon X \to Y$ induces 
  a group homomorphism $h = \Adj(f) \colon \Adj(X) \to \Adj(Y)$.
  In this way $\Adj(X)$ acts on $Y$, and $f$ becomes equivariant.
  In particular, every $g \in \Adj(X)'$ with $x^g = x$ is mapped 
  to $h(g) \in \Adj(Y)'$ with $y^{h(g)} = y$, which proves 
  the first claim.  Moreover, this construction respects composition.
\end{proof}

\begin{proposition}
  We have $\pi_1(Q,q^g) = \pi_1(Q,q)^g$ for every $g \in \Adj(Q)$, 
  or more generally for every $g \in \Aut(Q)$.  Thus, if $Q$ is connected, 
  or homogeneous, then the isomorphism class of the fundamental group 
  $\pi_1(Q,q)$ is independent of the choice of base point $q \in Q$.
  \qed
\end{proposition}

\subsection{Coverings and monodromy}

As for topological coverings, two groups naturally act on a 
quandle covering $p \colon \tilde{Q} \to Q$: the deck transformation
group $\Aut(p)$ acts on the left, while the adjoint group $\Adj(Q)$
and in particular its subgroup $\pi_1(Q,q)$ act on the right.
Both actions are connected as follows:

\begin{proposition}[monodromy action] \label{prop:MonodromyAction}
  Every galois covering $p \colon \tilde{Q} \to Q$ induces a natural 
  surjective group homomorphism $h \colon \pi_1(Q,q) \onto \Aut(p)$.

  More generally, every quandle extension 
  $E \colon \Lambda \curvearrowright \tilde{Q} \lto[p] Q$
  of a connected quandle $Q$ by a group $\Lambda$ induces 
  a natural group homomorphism $h \colon \pi_1(Q,q) \to \Lambda$.
  Moreover, $h$ is surjective if and only if $\tilde{Q}$ 
  is connected; in this case $p$ is a galois covering.

  In both settings, $h$ is an isomorphism if and only if 
  $p$ is the universal covering of $Q$.
\end{proposition}

\begin{proof}
  Every galois covering $p \colon \tilde{Q} \to Q$ 
  defines an extension, with the group $\Lambda = \Aut(p)$ 
  acting naturally on $\tilde{Q}$ by deck transformations
  (see Proposition \ref{prop:CoveringExtensions}).
  We will thus concentrate on the more general 
  formulation of extensions.

  Since the covering $p \colon \tilde{Q} \to Q$ 
  is equivariant under the natural action of $\Adj(Q)$,
  every $g \in \pi_1(Q,q)$ maps the fibre $F = p^{-1}(q)$ to itself.
  In particular, there exists a unique element $h(g) \in \Lambda$
  such that $\tilde{q}^g = h(g) \tilde{q}$.  
  For $g_1,g_2 \in \pi_1(Q,q)$ we find that 
  \[
  \tilde{q}^{g_1 g_2} = (h(g_1) \tilde{q})^{g_2} 
  = h(g_1) (\tilde{q}^{g_2}) =  h(g_1)h(g_2) \tilde{q}, 
  \]
  since both actions commute (see Proposition \ref{prop:CommutingActions}).  
  We conclude that $h(g_1 g_2) = h(g_1) h(g_2)$, whence $h$ is a group homomorphism.

  If $\tilde{Q}$ is connected, there exists for each 
  $\hat{q} \in F$ a group element $g \in \Adj(Q)'$ 
  such that $\tilde{q}^g = \hat{q}$ (see Remark \ref{rem:CommutatorOrbits}).  
  By equivariance this equation projects to $q^g = q$, and so we have 
  $g \in \pi_1(Q,q)$.  This implies that $h$ is surjective.

  Conversely, if $h$ is surjective, then $\tilde{Q}$ is connected:
  given $\hat{q} \in \tilde{Q}$, there exists $g_1 \in \Adj(Q)$
  such that $p(\hat{q})^{g_1} = q$, because $Q$ is connected.
  This implies that $\hat{q}^{g_1} = \lambda \tilde{q}$ for some 
  $\lambda \in \Lambda$.  Since $h$ is assumed to be surjective,
  there exists $g_2 \in \pi_1(Q,q)$ such that $h(g_2) = \lambda^{-1}$.
  We conclude that $\hat{q}^{g_1 g_2} = \tilde{q}$, as desired.

  Finally, if $h$ is an isomorphism, then $\Adj(Q)'$ acts freely on $\tilde{Q}$.
  We thus obtain an isomorphism between $(\tilde{Q},\tilde{q})$ and the
  universal covering constructed in Theorem \ref{thm:UniversalCovering}.
\end{proof}

\begin{proposition} \label{prop:CoveringFundamentalImage}
  For every quandle covering $p \colon (\tilde{Q},\tilde{q}) \onto (Q,q)$ 
  the induced group homomorphism $p_* \colon \pi_1(\tilde{Q},\tilde{q}) \to \pi_1(Q,q)$
  has image $\Im(p_*) = \{ g \in \Adj(Q)^\circ \mid \tilde{q}^g = \tilde{q} \}$ 
  and kernel $\ker(p_*) = \ker\left[ \; \Adj(p) \colon \Adj(\tilde{Q}) \to \Adj(Q) \; \right]$.
\end{proposition}

\begin{proof}
  We know by Proposition \ref{prop:CoveringAdjInnExtention}
  that $\phi = \Adj(p) \colon \Adj(\tilde{Q}) \onto \Adj(Q)$
  is a central extension. 
  By Definition \ref{def:Index} we have 
  $\varepsilon_{\smash{\tilde{Q}}} = \varepsilon_Q \circ \phi$, 
  so that $\phi$ maps $\Adj(\tilde{Q})^\circ$ onto $\Adj(Q)^\circ$.
  The action of $\Adj(Q)$ on $\tilde{Q}$ is such 
  that $\tilde{q}^{\tilde{g}} = \tilde{q}^{\phi(\tilde{g})}$ 
  for all $\tilde{g} \in \Adj(\tilde{Q})$, cf.\ Remark \ref{rem:CoveringAction}.  
  
  If $\tilde{g} \in \pi_1(\tilde{Q},\tilde{q})$ then 
  $g = \phi(\tilde{g})$ satisfies $g \in  \Adj(Q)^\circ$ and $\tilde{q}^g = \tilde{q}$.
  Conversely, for each $g \in \Adj(Q)^\circ$ with $\tilde{q}^g = \tilde{q}$, every
  preimage $\tilde{g} \in \phi^{-1}(g)$ satisfies $\tilde{g} \in \Adj(\tilde{Q})^\circ$
  and $\tilde{q}^{\tilde{g}} = \tilde{q}$, whence $\tilde{g} \in \pi_1(\tilde{Q},\tilde{q})$
  and $g = p_*(\tilde{g})$.  Existence of $\tilde{g}$ is ensured by the surjectivity of $\phi$.
  
  Finally, $\tilde{g} \in \ker(p_*)$ is equivalent to 
  $\tilde{g} \in \Adj(\tilde{Q})^\circ$ and $\tilde{q}^{\tilde{g}} = \tilde{q}$
  and $\phi(\tilde{g}) = 1$.  This last condition entails the
  two previous ones: if $\phi(\tilde{g}) = 1$ then $\tilde{g} \in \Adj(\tilde{Q})^\circ$
  and $\tilde{q}^{\tilde{g}} = \tilde{q}^{\phi(\tilde{g})} = \tilde{q}$,
  so that $\tilde{g} \in \ker(p_*)$.  We conclude that $\ker(p_*) = \ker(\Adj(p))$.
\end{proof}

\begin{warning}
  For a connected quandle covering $p \colon \tilde{Q} \onto Q$ 
  the adjoint group homomorphism $\Adj(\tilde{Q}) \onto \Adj(Q)$ 
  can have non-trivial kernel, and so 
  $p_* \colon \pi_1(\tilde{Q},\tilde{q}) \to \pi_1(Q,q)$ 
  is in general not injective.  In this respect the covering
  theory of quandles differs sharply from coverings of topological 
  spaces, where $p_*$ is injective for every covering.
\end{warning}

\begin{example} \label{exm:LinearGroups2}
  As in Example \ref{exm:LinearGroups1}, consider a group $\tilde{G}$ 
  and a conjugacy class $\tilde{Q} \subset \tilde{G}$ such that 
  $\tilde{G} = \gen{\tilde{Q}}$.  Assume that $\Lambda \subset Z(\tilde{G})$ 
  is a non-trivial central subgroup such that $\Lambda \cdot \tilde{Q} = \tilde{Q}$.
  The quotient map $p \colon \tilde{G} \to G := \tilde{G}/\Lambda$ sends 
  $\tilde{Q}$ to a conjugacy class $Q = p(\tilde{Q})$ in $G$ with $G = \gen{Q}$.
  We thus obtain an extension $\Lambda \curvearrowright \tilde{Q} \lto[p] Q$.

  Since $\tilde{Q}$ embeds into a group, the adjoint map
  $\tilde{Q} \to \Adj(\tilde{Q})$ is injective.  
  The group homomorphism $h = \Adj(p) \colon \Adj(\tilde{Q}) \to \Adj(Q)$
  is not injective because $\tilde{q}$ and $\lambda\tilde{q}$ 
  with $\lambda \in \Lambda \minus \{1\}$, are distinct 
  in $\tilde{Q}$ but get identified in $Q$.  The element 
  $\tilde{z} = \adj(\tilde{q})^{-1} \adj(\lambda\tilde{q})$ in $\Adj(\tilde{Q})'$
  is thus contained in $\ker(h)$, and thus in the centre of $\Adj(\tilde{Q})$.
  In particular $\tilde{q}^{\tilde{z}} = \tilde{q}$, and so
  $\tilde{z} \in \pi_1(\tilde{Q},\tilde{q})$ is a non-trivial 
  element that maps to $p_*(\tilde{z}) = 1$ in $\pi_1(Q,q)$.
\end{example}



\subsection{The lifting criterion} \label{sub:LiftingCriterion}

As for topological coverings, the fundamental group provides 
a simple criterion for the lifting over a quandle covering:

\begin{proposition}[lifting criterion] \label{prop:LiftingCriterion}
  Let $f \colon (X,x) \to (Q,q)$ be a quandle homomorphism, and let 
  $p \colon (\tilde{Q},\tilde{q}) \to (Q,q)$ be a quandle covering.
  Assume further that $(X,x)$ is connected.  Then there exists a lifting
  $\tilde{f} \colon (X,x) \to (\tilde{Q},\tilde{q})$ if and only if
  $f_* \pi_1(X,x) \subset p_* \pi_1(\tilde{Q},\tilde{q})$. 
\end{proposition}

\begin{proof}
  We already know from Corollary \ref{cor:CoveringMorphismUniqueness} 
  that $\tilde{f}$ is unique, and so we only have to consider existence.
  Let us begin with the easy case:
  If a lifting $\tilde{f}$ exists, then $f = p \tilde{f}$ 
  implies $f_* = p_* \tilde{f}_*$ and thus $f_* \pi_1(X,x) 
  = p_* \tilde{f}_* \pi_1(X,x) \subset p_* \pi_1(\tilde{Q},\tilde{q})$.

  Conversely, assume $f_* \pi_1(X,x) \subset p_* \pi_1(\tilde{Q},\tilde{q})$.
  Since $p$ is a covering, the group $\Adj(Q)$ acts on $\tilde{Q}$ by inner automorphisms.
  The quandle homomorphism $f \colon X \to Q$ induces a group homomorphism 
  $f_* = \Adj(f) \colon \Adj(X) \to \Adj(Q)$, and in this way $\Adj(X)$ also acts on $\tilde{Q}$.
  By connectedness, every element of $X$ can be written as $x^g$ with some $g \in \Adj(X)'$.
  We can thus define $\tilde{f} \colon (X,x) \to (\tilde{Q},\tilde{q})$
  by setting $\tilde{f} \colon x^g \mapsto \tilde{q}^g$,
  and our hypothesis ensures that this is well-defined.
  By construction, the map $\tilde{f}$ is $\Adj(X)'$-equivariant.  
  Both maps $p \tilde{f}$ and $f$ are $\Adj(X)'$-equivariant and 
  coincide in $x$; since $X$ is connected we obtain $p \tilde{f} = f$.
  As in Proposition \ref{prop:CoveringMorphismEquivariance}
  we conclude that $\tilde{f}$ is a quandle homomorphism.
\end{proof}

\begin{definition} \label{Def:SimplyConnected}
  A quandle $Q$ is \emph{simply connected}
  if it is connected and $\pi_1(Q,q) = \{1\}$.
\end{definition}

Notice that connectedness implies that $\pi_1(Q,q) \cong \pi_1(Q,q')$
for all $q,q' \in Q$.  It thus suffices to verify triviality
of $\pi_1(Q,q)$ for \emph{one} base point $q \in Q$; the property 
of being simply connected is independent of this choice, and hence well-defined.

\begin{proposition} \label{prop:SimplyConnected}
  For a quandle $Q$ the following properties are equivalent:
  \begin{enumerate}
  \item
    The quandle $Q$ is simply connected.
  \item
    Every covering $p \colon \tilde{Q} \to Q$ is equivalent to 
    a trivial covering $\pr_1 \colon Q \times F \to Q$. 
  \item
    Every quandle homomorphism $f \colon (Q,q) \to (\bar{Q},\bar{q})$ lifts uniquely 
    over each quandle covering $p \colon (\tilde{Q},\tilde{q}) \to (\bar{Q},\bar{q})$.
  \item
    Every covering $p \colon (Q,q) \to (\bar{Q},\bar{q})$ 
    is universal in the category $\Coverings(\bar{Q},\bar{q})$.
  \end{enumerate}
\end{proposition}

\begin{proof}
  (1) $\Rightarrow$ (2):
  We choose a base point $q \in Q$ and define $F := p^{-1}(q)$.
  According to the Lifting Criterion, for each $\tilde{q} \in F$ 
  there exists a unique quandle homomorphism 
  $\phi_{\tilde{q}} \colon (Q,q) \to (\tilde{Q},\tilde{q})$
  such that $p \circ \phi_{\tilde{q}} = \id_Q$.
  Its image is the connected component of $\tilde{q}$ in $\tilde{Q}$.
  We thus have a bijection $\psi \colon \pi_0(\tilde{Q}) \to F$
  such that $\psi([\tilde{q}]) = \tilde{q}$ for every $\tilde{q} \in F$.
  Putting this information together 
  we obtain mutually inverse quandle isomorphisms
  $\Phi \colon Q \times F \to \tilde{Q}$, 
  $\Phi(x,\tilde{q}) = \phi_{\tilde{q}}(x)$
  and $\Psi \colon \tilde{Q} \to Q \times F$, 
  $\Psi(\tilde{x}) = (p(\tilde{x}),\psi([\tilde{x}]))$.
  
  (2) $\Rightarrow$ (3): 
  By hypothesis (2) and Remark \ref{rem:AlmostTrivialCovering},
  $Q$ must be connected, which ensures uniqueness.
  Existence follows from the pull-back construction,
  because $f^* p$ is a covering over $(Q,q)$ and trivial by hypothesis.

  (3) $\Rightarrow$ (4): 
  This is clear from Definition \ref{def:UniversalCovering}.

  (4) $\Rightarrow$ (1): 
  The identity $\id_Q \colon (Q,q) \to (Q,q)$ is a covering.
  If it is universal, then $Q$ must be connected 
  by Remark \ref{rem:AlmostTrivialCovering}.
  Moreover, $(Q,q)$ must be isomorphic to the explicit model
  $(\tilde{Q},\tilde{q})$ of Theorem \ref{thm:UniversalCovering}
  via the projection map $p \colon (\tilde{Q},\tilde{q}) \to (Q,q)$.
  This implies $\pi_1(Q,q) = \{1\}$, whence $Q$ is simply connected.
\end{proof}

\begin{example}
  For a long knot $L$, the knot quandle $Q_L$ is simply connected 
  by \cite[Theorem 30]{Eisermann:2003}.  The natural quandle
  projection $Q_L \onto Q_K$ is thus the universal covering
  of the knot quandle $Q_K$ associated to the closed knot $K$.
\end{example}

\begin{warning}
  For a universal quandle covering $p \colon (\tilde{Q},\tilde{q}) \to (Q,q)$ 
  the covering quandle $\tilde{Q}$ need not be simply connected.  
  This is another aspect in which quandle coverings differ from 
  topological coverings, where every universal covering is simply connected.
\end{warning}

\begin{example} \label{exm:LinearGroups3}
  We continue Example \ref{exm:LinearGroups2} using the same notation.
  The universal covering $\hat{p} \colon (\hat{Q},\hat{q}) \to (Q,q)$ of $(Q,q)$ 
  induces a covering $\tilde{p} \colon (\hat{Q},\hat{q}) \to (\tilde{Q},\tilde{q})$.
  This means that $\Adj(\hat{p}) \colon \Adj(\hat{Q}) \onto \Adj(Q)$
  factors as $\Adj(\hat{Q}) \lonto[g] \Adj(\tilde{Q}) \lonto[h] \Adj(Q)$.
  We have already found a non-trivial element $\tilde{z} \in \pi_1(\tilde{Q},\tilde{q})$
  with $h(\tilde{z}) = 1$ in $\pi_1(Q,q)$.  Every preimage $\hat{z} \in g^{-1}(\tilde{z})$
  lies in centre of $\Adj(\hat{Q})$ and also in the commutator subgroup, 
  and thus provides a non-trivial element $\hat{z} \in \pi_1(\hat{Q},\hat{q})$.
\end{example}

\subsection{Galois correspondence} \label{sub:ConnectedGaloisCorrespondence}

Let $(Q,q)$ be a connected quandle.  We wish to establish 
a correspondence between the following two categories.  
On the one hand, we have the category $\Coverings_\connected(Q,q)$ formed by 
pointed connected coverings $p \colon (\tilde{Q},\tilde{q}) \to (Q,q)$
and their pointed covering morphisms.  On the other hand, we have 
the category $\Subgroups(\pi_1(Q,q))$ formed by subgroups of $\pi_1(Q,q)$ 
and homomorphisms given by inclusion.  The Galois correspondence establishes 
a natural equivalence $\Coverings_\connected(Q,q) \cong \Subgroups(\pi_1(Q,q))$.

\begin{remark}
  In $\Subgroups(\pi_1(Q,q))$ inclusion defines a partial order
  on the set of subgroups.  Likewise, in $\Coverings_\connected(Q,q)$ 
  each set of covering morphisms $\Hom(p,p')$ is either empty or contains exactly 
  one element (see Corollary \ref{cor:CoveringMorphismUniqueness}),
  which expresses a partial preorder. 
\end{remark}

\begin{lemma}
  There exists a unique functor 
  $\Phi \colon \Coverings_\connected(Q,q) \to \Subgroups(\pi_1(Q,q))$ 
  mapping each covering $p \colon (\hat{Q},\hat{q}) \to (Q,q)$
  to the subgroup $p_* \pi_1(\hat{Q},\hat{q}) \subset \pi_1(Q,q)$.
\end{lemma}

\begin{proof}
  Obviously $\Phi$ is well-defined on objects.
  Every covering morphism $\phi$ from $p$ to $p'$ entails that 
  $p_* \pi_1(\hat{Q},\hat{q}) = p_* \phi_* \pi_1(\hat{Q},\hat{q})
  \subset p_* \pi_1(\hat{Q}',\hat{q}')$, so that $\Phi$ is indeed a functor.
\end{proof}

\begin{lemma}
  There exists a unique functor 
  $\Psi \colon \Subgroups(\pi_1(Q,q)) \to \Coverings_\connected(Q,q)$
  mapping each subgroup $K \subset \pi_1(Q,q)$ to the quotient
  $\tilde{Q}_K := K \backslash \tilde{Q}$ of the universal
  covering $\tilde{Q}$. 
\end{lemma}

\begin{proof}
  We consider the  universal covering $p \colon (\tilde{Q},\tilde{q}) \to (Q,q)$
  constructed in Lemma \ref{lem:UniversalCovering}.  Given a subgroup 
  $K \subset \pi_1(Q,q)$, we identify $K$ with the corresponding subgroup of $\Aut(p)$, 
  via the monodromy action explained in Proposition \ref{prop:MonodromyAction}.
  This allows us to define the quotient $\tilde{Q}_K := K \backslash \tilde{Q}$ 
  with base point $\tilde{q}_K = [\tilde{q}]$ and 
  projection $p_K \colon (\tilde{Q}_K,\tilde{q}_K) \to (Q,q)$
  defined by $p_K([\tilde{x}]) = p(\tilde{x})$.  
  The result is the covering $\Psi(K) := p_K$ we wish to consider.  

  Moreover, if $K \subset L \subset \pi_1(Q,q)$, then 
  the covering $p_L$ is a quotient of the covering $p_K$.
  We thus have a covering morphism from $p_K$ to $p_L$, 
  so that $\Psi$ is indeed a functor. 
\end{proof}

\begin{theorem}[Galois correspondence] \label{thm:ConnectedCorrespondence}
  Let $(Q,q)$ be a connected quandle.  Then the functors
  $\Phi \colon \Coverings_\connected(Q,q) \to \Subgroups(\pi_1(Q,q))$ and
  $\Psi \colon \Subgroups(\pi_1(Q,q)) \to \Coverings_\connected(Q,q)$ 
  establish a natural equivalence between the category of pointed connected 
  coverings of $(Q,q)$ and the category of subgroups of $\pi_1(Q,q)$.
\end{theorem}

\begin{proof}
  We will first prove that $\Phi \Psi = \id$.
  Consider a subgroup $K \subset \pi_1(Q,q)$ and the associated 
  covering $p_K \colon (\tilde{Q}_K,\tilde{q}_K) \to (Q,q)$.
  By Proposition \ref{prop:CoveringFundamentalImage} we know 
  that the image group $(p_K)_* \pi_1(\tilde{Q}_K,\tilde{q}_K)$
  consists of all $g \in \Adj(Q)'$ such that $\tilde{q}_K^g = \tilde{q}_K$.
  Comparing this with the construction of the universal covering
  $(\tilde{Q},\tilde{q})$ and its quotient $(\tilde{Q}_K,\tilde{q}_K)$
  we obtain precisely the group $K$ with which we started out.

  Conversely, let us prove that $\Psi \Phi \cong \id$.
  For every connected covering $p \colon (\hat{Q},\hat{q}) \to (Q,q)$ 
  the associated group $K = p_* \pi_1(\hat{Q},\hat{q})$ defines 
  a covering $p_K \colon (\tilde{Q}_K,\tilde{q}_K) \to (Q,q)$ as above.
  We already know that 
  $(p_K)_* \pi_1(\tilde{Q}_K,\tilde{q}_K) = K = p_* \pi_1(\hat{Q},\hat{q})$.
  The Lifting Criterion (Proposition \ref{prop:LiftingCriterion})
  implies that there exist covering morphisms
  $f \colon (\tilde{Q}_K,\tilde{q}_K) \to (\hat{Q},\hat{q})$
  and $g \colon (\hat{Q},\hat{q}) \to (\tilde{Q}_K,\tilde{q}_K)$.
  By the usual uniqueness argument (Corollary \ref{cor:CoveringMorphismUniqueness})
  we conclude that $f \circ g = \id_{\hat{Q}}$ and $g \circ f = \id_{\tilde{Q}_K}$
\end{proof}

\begin{proposition}[monodromy and deck transformation group]
  Consider a connected covering $p \colon (\tilde{Q},\tilde{q}) \to (Q,q)$ and 
  the associated subgroup $K = p_* \pi_1(\tilde{Q},\tilde{q}) \subset \pi_1(Q,q)$.
  \begin{enumerate}
  \item
    The natural right action $F \times \pi_1(Q,q) \to F$ induces 
    a bijection between the fibre $F = p^{-1}(q)$ and 
    the quotient set $K \backslash \pi_1(Q,q)$.
    In particular, the cardinality of $F$ equals
    the index of the subgroup $K$ in $\pi_1(Q,q)$.
  \item
    Let $N = \{ g \in \pi_1(Q,q) \mid K^g = K \}$ 
    be the normalizer of $K$ in $\pi_1(Q,q)$.  
    There exists a covering transformation 
    $(\tilde{Q},\tilde{q}) \to (\tilde{Q},\hat{q})$ 
    if and only if there exists an element $g \in N$ 
    such that $\tilde{q}^g = \hat{q}$.
  \item
    We have a natural short exact sequence $K \into N \onto \Aut(p)$.  
    The covering $p$ is galois if and only if the subgroup $K$ is normal in $\pi_1(Q,q)$.
    In this case the deck transformation group is $\Aut(p) \cong \pi_1(Q,q)/K$.
  \end{enumerate}
\end{proposition}

\begin{proof}
  Since $\tilde{Q}$ is connected, $\pi_1(Q,q)$ acts transitively
  on the fibre $F = p^{-1}(q)$.  The stabilizer of $\tilde{q}$ is precisely 
  the subgroup $K$, cf.\ Proposition \ref{prop:CoveringFundamentalImage}.
  Given $g \in \pi_1(Q,q)$ there exists a covering 
  automorphism $\phi \colon (\tilde{Q},\tilde{q}) \to (\tilde{Q},\tilde{q}^g)$
  if and only if the subgroups $p_* \pi_1(\tilde{Q},\tilde{q}) = K$ 
  and $p_* \pi_1(\tilde{Q},\tilde{q}^g) = K^g$ coincide 
  (see the Lifting Criterion, Proposition \ref{prop:LiftingCriterion}).
  In this case $\phi$ is unique, and so $g \mapsto \phi$ defines 
  a surjective group homomorphism $N \onto \Aut(p)$, 
  as in the proof of Proposition \ref{prop:MonodromyAction}.
\end{proof}


\section{Classification of non-connected coverings} \label{sec:NonConnectedCoverings}

\subsection{Non-connected covering quandles}

In this section we deal with coverings $p \colon \tilde{Q} \to Q$ 
where the base quandle $Q$ is connected but the covering quandle 
$\tilde{Q}$ can be non-connected.  Non-connected base quandles
are more delicate and will be treated in the next section.


\begin{proposition}
  Consider a family of quandle coverings $p_i \colon \tilde{Q}_i \to Q$ 
  indexed by $i \in I$.  Let $\tilde{Q} = \bigsqcup_{i \in I} \tilde{Q}_i \times \{i\}$
  be their disjoint union with projection $p \colon \tilde{Q} \to Q$, $p(a,i) = p_i(a)$.
  There exists a unique quandle structure on $\tilde{Q}$ that extends 
  the one on each $\tilde{Q}_i$ and turns $p$ into a quandle covering.
  The result is called the \emph{union} of the given quandle coverings over $Q$, 
  denoted by $(\tilde{Q},p) = \bigoplus_{i \in I}(\tilde{Q}_i,p_i)$.
\end{proposition}

\begin{proof}
  The point is to define the quandle structure on $\tilde{Q}$.
  Since each $p_i$ is a covering, the base quandle $Q$ acts on $\tilde{Q}_i$ 
  such that $a \ast b = a \ast p_i(b)$ for all $a,b \in \tilde{Q}_i$.
  If there is a compatible quandle structure on $\tilde{Q}$ such that
  $p \colon \tilde{Q} \to Q$ becomes a covering, then $Q$ acts on $\tilde{Q}$ 
  and we must have $(a,i) \ast (b,j) = (a,i) \ast p_j(b,j) = (a \ast p_j(b), i)$.
  This shows that there can be at most one such structure.
  In order to prove existence, we equip $\tilde{Q}$ with
  the operation $(a,i) \ast (b,j) := (a * p_j(b), i)$.
  If $I$ is non-empty, then it is easily verified that 
  this definition turns $\tilde{Q}$ into a quandle, 
  and that $p$ becomes a quandle covering of $Q$.
\end{proof}

\begin{proposition}
  Let $p \colon \tilde{Q} \to Q$ be a covering of the connected quandle $Q$.
  We can decompose $\tilde{Q}$ into connected components
  $(\tilde{Q}_i)_{i \in I}$ and define $p_i \colon \tilde{Q}_i \to Q$
  by restriction.  Then each $p_i$ is a covering, and 
  $(\tilde{Q},p) = \bigoplus_{i \in I}(\tilde{Q}_i,p_i)$ is their union.
\end{proposition}

\begin{proof}
  Notice that each $\tilde{Q}_i$ is an orbit under the action 
  of $\Adj(Q)$ on $\tilde{Q}$, and each $p_i$ is a covering 
  because it is an $\Adj(Q)$-equivariant map.
  By construction we have the equality of sets and maps,
  $(\tilde{Q},p) = \bigoplus_{i \in I}(\tilde{Q}_i,p_i)$.
  The equality of their quandle structures follows from 
  the uniqueness part of the previous proposition.
\end{proof}

\subsection{Galois correspondence} \label{sub:SemiconnectedGaloisCorrespondence}

Theorem \ref{thm:ConnectedCorrespondence} above established the correspondence 
between connected coverings and subgroups of the fundamental group.  
In the general setting it is more convenient to classify  coverings
by actions of the fundamental group on the fibre.

\begin{definition}[the category of $G$-sets]
  Let $G$ be a group.  A \emph{$G$-set} is a pair $(X,\alpha)$ consisting of a set $X$ 
  and a right action $\alpha \colon X \times G \to X$, denoted by $\alpha(x,g) = x^g$.  
  A \emph{morphism} $\phi \colon (X,\alpha) \to (Y,\beta)$ 
  between two $G$-sets 
  is an equivariant map $\phi \colon X \to Y$, i.e.\ satisfying 
  $\phi(x^g) = \phi(x)^g$ for all $x \in X$ and $g \in G$.  
  The class of $G$-sets and their morphisms form a category, 
  denoted by $\Actions(G)$.
\end{definition}


\begin{lemma}
  There exists a canonical functor $\Phi \colon \Coverings(Q) \to \Actions(\pi_1(Q,q))$ 
  mapping each covering $p \colon \tilde{Q} \to Q$ to $(F,\alpha)$ where 
  $F = p^{-1}(q)$ is the fibre over $q$, and $\alpha \colon F \times \pi_1(Q,q) \to F$
  is the monodromy action.
\end{lemma}

\begin{proof}
  Given a covering $p \colon \tilde{Q} \to Q$,
  the natural action of $\Adj(Q)$ on $\tilde{Q}$ restricts
  to an action of $\pi_1(Q,q)$ on the fibre $F = p^{-1}(q)$.  
  This defines $\Phi$ on objects.
  
  Every covering morphism $\phi \colon \tilde{Q} \to \hat{Q}$ 
  is equivariant with respect to the action of $\Adj(Q)$.
  It maps the fibre $F = p^{-1}(q)$ to the fibre $\hat{F} = \hat{p}^{-1}(q)$, 
  and the restriction $\phi_q \colon F \to \hat{F}$ is equivariant with 
  respect to the action of $\pi_1(Q,q)$.  Hence $\Phi$ is indeed a functor.
\end{proof}

\begin{lemma}
  There exists a canonical functor $\Psi \colon \Actions(\pi_1(Q,q)) \to \Coverings(Q)$
  mapping each action $\alpha \colon F \times \pi_1(Q,q) \to F$ to the
  covering $p_\alpha \colon \tilde{Q}_\alpha \to Q$ with 
  $\tilde{Q}_\alpha = (F \times \tilde{Q})/{\pi_1(Q,q)}$,
  where $\tilde{Q}$ is the universal connected covering of $Q$.
\end{lemma}

\begin{proof}
  We start with the universal connected covering 
  $p \colon (\tilde{Q},\tilde{q}) \to (Q,q)$.
  According to Proposition \ref{prop:MonodromyAction}
  we have a group isomorphism $h \colon \pi_1(Q,q) \isoto \Aut(p)$,
  such that $h(g) \tilde{q} = \tilde{q}^g$ for all $g \in \pi_1(Q,q)$.
  Given $(F,\alpha) \in \Actions(\pi_1(Q,q))$, we quotient 
  the product $F \times \tilde{Q}$ by the equivalence relation
  $(x^g,\tilde{a}) \sim (x,h(g)\tilde{a})$ for all $x \in F$, 
  $\tilde{a} \in \tilde{Q}$, and $g \in \pi_1(Q,q)$.
  The quotient $\tilde{Q}_\alpha := (F \times \tilde{Q})/_\sim$ 
  inherits the quandle structure 
  $[x,\tilde{a}] \ast [y,\tilde{b}] := [x,\tilde{a} \ast \tilde{b}]$.
  The projection $p_\alpha \colon \tilde{Q}_\alpha \to Q$, 
  $p_\alpha([x,\tilde{a}]) := p(\tilde{a})$ is well-defined 
  and a quandle covering.  As a consequence,
  the action of $\Adj(Q)$ on $\tilde{Q}_\alpha$ is given by 
  $[x,\tilde{a}]^g = [x,\tilde{a}^g]$ for all $g \in \Adj(Q)$.

  A morphism $\phi \colon (X,\alpha) \to (Y,\beta)$ of $G$-sets
  induces a map $\phi\times\id \colon X \times \tilde{Q} \to Y \times \tilde{Q}$
  that descends to a quandle homomorphism on the quotients, 
  $\bar{\phi} \colon \tilde{Q}_\alpha \to \tilde{Q}_\beta$.
  This turns out to be a covering morphism from $p_\alpha$
  to $p_\beta$, so that $\Psi$ is indeed a functor. 
\end{proof}

\begin{theorem}[Galois correspondence] \label{thm:NonConnectedCorrespondence}
  Let $(Q,q)$ be a connected quandle.  The functors
  $\Phi \colon \Coverings(Q) \to \Actions(\pi_1(Q,q))$ and
  $\Psi \colon \Actions(\pi_1(Q,q)) \to \Coverings(Q)$ establish
  a natural equivalence between the category of coverings of $Q$ 
  and the category of sets endowed with an action of $\pi_1(Q,q)$.
\end{theorem}

\begin{proof}
  Before we begin, let us point out that strictly speaking
  the compositions $\Psi \Phi$ and $\Phi \Psi$ are \emph{not}
  the identity functors.  They are, however, naturally equivalent
  to the identity functors, in the sense of \cite[\textsection I.4]{MacLane:1998},
  and this is what we have to show.

  We will first prove that $\Phi \Psi \cong \id$.
  Consider an action $\alpha \colon X \times \pi_1(Q,q) \to X$
  and the associated covering $p_\alpha \colon \tilde{Q}_\alpha \to Q$
  with fibre $F_\alpha := p_\alpha^{-1}(q)$.  Recall that $\Aut(p)$ acts
  freely and transitively from the left on the fibre $p^{-1}(q)$ of 
  the universal covering $p \colon (\tilde{Q},\tilde{q}) \to (Q,q)$.  
  The map $\psi_\alpha \colon X \to F_\alpha$, $x \mapsto [x,\tilde{q}]$,
  is thus a bijection.  Moreover, we find
  \[
  \psi_\alpha(x^g) = [x^g,\tilde{q}] = [x,h(g)\tilde{q}] = [x,\tilde{q}^g] 
  = [x,\tilde{q}]^g = \psi_\alpha(x)^g 
  \]
  for every $g \in \pi_1(Q,q)$. 
  This shows that $\psi_\alpha \colon X \to F_\alpha$ 
  is an equivalence of $\pi_1(Q,q)$-sets, as claimed.
  Naturality in $\alpha$ is easily verified.

  Conversely, let us prove that $\Psi \Phi \cong \id$.
  Consider a quandle covering $\hat{p} \colon \hat{Q} \to Q$ 
  with fibre $F = \hat{p}^{-1}(q)$ and monodromy action
  $\alpha \colon F \times \pi_1(Q,q) \to F$.
  The universal property of the covering 
  $p \colon (\tilde{Q},\tilde{q}) \to (Q,q)$
  ensures that there exists a unique covering morphism 
  $\phi_{\hat{p}} \colon F \times \tilde{Q} \to \hat{Q}$ 
  over $Q$ such that $\phi_{\hat{p}}(x,\tilde{q}) = x$ for all $x \in F$.
  More explicitly, this map is given by $(x,(q,g)) \mapsto x^g$
  for all $x \in F$ and $g \in \Adj(Q)'$.
  By construction, this map is surjective and equivariant 
  with respect to the action of $\Adj(Q)'$.  

  For $g \in \pi_1(Q,q)$ we find $\phi_{\hat{p}}(x^g,\tilde{a}) = \phi_{\hat{p}}(x,h(g)\tilde{a})$
  for all $x \in F$ and $\tilde{a} \in \tilde{Q}$.
  This means that $\phi_{\hat{p}}$ descends to a covering morphism 
  $\bar\phi_{\hat{p}} \colon \tilde{Q}_\alpha \to \hat{Q}$.
  Conversely, if $\phi_{\hat{p}}(x,\tilde{a}) = \phi_{\hat{p}}(y,\tilde{b})$,
  then both maps $\phi_{\hat{p}}(x,-)$ and $\phi_{\hat{p}}(y,-)$
  have as image the same component of $\hat{Q}$, which 
  takes us back to the case of connected coverings.
  We thus see that $(x,\tilde{a})$ and $(y,\tilde{b})$ get identified 
  in $\tilde{Q}_\alpha$, which proves that $\bar\phi_{\hat{p}}$ 
  is a covering isomorphism.  Naturality in $\hat{p}$ is easily verified.
\end{proof}

\begin{theorem} \label{thm:ExtensionFundamentalCorrespondence}
  Let $Q$ be a connected quandle with base point $q \in Q$
  and let $\Lambda$ be a group.  There exists a natural 
  bijection $\Ext(Q,\Lambda) \cong \Hom(\pi_1(Q,q),\Lambda)$.
  If $\Lambda$ is an abelian group, or more generally a module over some ring $R$,
  then both objects carry natural $R$-module structures and 
  the bijection is an $R$-module isomorphism.
\end{theorem}

\begin{proof}
  Every extension $E \colon \Lambda \curvearrowright \tilde{Q} \lto[p] Q$ induces 
  a group homomorphism $h \colon \pi_1(Q,q) \to \Lambda$ as in Proposition 
  \ref{prop:MonodromyAction}.  Choosing a base point $\tilde{q}$ in the
  fibre $F = p^{-1}(q)$, we can identify $\Lambda$ with $F$ via
  the bijection $\Lambda \isoto F$, $\lambda \mapsto \lambda \tilde{q}$.
  The monodromy action of $\pi_1(Q,q)$ then translates to right multiplication
  $\alpha \colon \Lambda \times \pi_1(Q,q) \to \Lambda$
  with $(\lambda,g) \mapsto \lambda \cdot h(g)$.
  
  Conversely, every group homomorphism $h$ defines a right action 
  $\alpha \colon \Lambda \times \pi_1(Q,q) \to \Lambda$ 
  by $(\lambda,g) \mapsto \lambda \cdot h(g)$.  
  Via Theorem \ref{thm:NonConnectedCorrespondence} the action $\alpha$
  corresponds to a covering $p_\alpha \colon \tilde{Q}_\alpha \onto Q$.  
  Multiplication on the left defines an action of $\Lambda$ on $\Lambda \times \tilde{Q}$, 
  which descends to the quotient $\tilde{Q}_\alpha$ and defines an extension 
  $E \colon \Lambda \curvearrowright \tilde{Q}_\alpha \lto[p] Q$.

  These constructions are easily seen to establish a natural bijection, as desired. 
\end{proof}


\section{Non-connected base quandles} \label{sec:NonConnectedBase}

\subsection{Graded quandles} \label{sub:GradedQuandles}

So far we have concentrated on connected base quandles.
In order to develop a covering theory over non-connected quandles
we have to treat all components  individually yet simultaneously.  
The convenient way to do this is to index the components by some fixed set $I$,
and then to deal with $I$-graded objects throughout.  The following
example illustrates the notions that will appear: 

\begin{example}
  Consider a quandle $Q$ and its decomposition 
  $Q = \bigsqcup_{i \in I} Q_i$ into connected components.
  For every covering $p \colon \tilde{Q} \to Q$ the quandle $\tilde{Q}$ 
  is graded, with $\tilde{Q}_i = p^{-1}(Q_i)$, and $p$ is a graded map,
  with $p_i \colon \tilde{Q}_i \to Q_i$ given by restriction.
  Every deck transformation $\phi \colon \tilde{Q} \isoto \tilde{Q}$
  is a graded map with $\phi_i \colon \tilde{Q}_i \isoto \tilde{Q}_i$.
  The deck transformation group $G = \Aut(p)$ is a graded group,
  with $G_i$ acting by covering transformations on $\tilde{Q}_i$,
  and this action turns $\tilde{Q}$ into a graded $G$-set.
\end{example}

The following definitions make the notions of this example explicit.
In the sequel we fix an index set $I$.
Whenever the context determines $I$ without ambiguity, 
the term ``graded'' will be understood to mean ``$I$-graded'', 
that is, graded with respect to our fixed set $I$.

\begin{definition}[graded quandles]
  A \emph{graded quandle} is a quandle $Q = \bigsqcup_{i \in I} Q_i$
  partitioned into subsets $(Q_i)_{i \in I}$ 
  such that $Q_i \ast Q_j = Q_i$ for all $i,j \in I$.
  This is equivalent to saying that each $Q_i$ is a union of connected components.
  A grading is equivalent to a quandle homomorphism $\gr \colon Q \to I$
  from $Q$ to the trivial quandle $I$ with fibres $Q_i = \gr^{-1}(i)$.

  A \emph{homomorphism} $\phi \colon Q \to Q'$ of graded quandles 
  is a quandle homomorphism such that $\phi(Q_i) \subset Q'_i$ for all $i \in I$,
  or equivalently $\gr = \gr' \circ \phi$.  Obviously, $I$-graded quandles 
  and their homomorphisms form a category, denoted $\Quandles_I$.
\end{definition}

\begin{definition}[graded groups]
  A \emph{graded group} is a group $G = \prod_{i \in I} G_i$ 
  together with the collection of groups $(G_i)_{i \in I}$ 
  that constitute the composition of $G$ as a product.
  A \emph{homomorphism} of graded groups $f \colon G \to H$ 
  is a product $f = \prod_{i \in I} f_i$ of homomorphisms 
  $f_i \colon G_i \to H_i$.  Obviously, $I$-graded groups and 
  their homomorphisms form a category, denoted $\Groups_I$.
  A graded subgroup of $G = \prod_{i \in I} G_i$ is a product
  $H = \prod_{i \in I} H_i$ of subgroups $H_i \subset G_i$.
\end{definition}

\begin{definition}[graded $G$-sets]
  A \emph{graded set} is a disjoint union $X = \bigsqcup_{i \in I} X_i$
  together with the partition $(X_i)_{i \in I}$.
  A \emph{graded map} $\phi \colon X \to Y$ between graded sets
  is a map satisfying $\phi(X_i) \subset Y_i$ for all $i \in I$.
  Graded sets and maps form a category, denoted $\category{Sets}_I$.

  A graded (right) action of a graded group $G$ on a graded set $X$
  is a collection of (right) actions $\alpha_i \colon X_i \times G_i \to X_i$,
  denoted by $\alpha(x,g) = x^g$.  This defines an action
  of $G$ on $X$ via the canonical projections $G \onto G_i$.
  A \emph{graded $G$-set} is a pair $(X,\alpha)$ consisting of 
  a graded set $X$ and a graded action $\alpha$ of $G$ on $X$.
  A \emph{morphism} $\phi \colon (X,\alpha) \to (Y,\beta)$ 
  between graded $G$-sets 
  is a graded map $\phi \colon X \to Y$ satisfying 
  $\phi(x^g) = \phi(x)^g$ for all $x \in X$ and $g \in G$.  
  Graded $G$-sets and their morphisms form a category, 
  denoted by $\Actions_I(G)$.
\end{definition}

\begin{remark}
  If the index set $I = \{*\}$ consists of one single element,
  then all gradings are trivial, and the categories of graded
  quandles, groups, and sets coincide with the usual (non-graded) notions.
\end{remark}

\begin{remark}
  As Mac\,Lane \cite[\textsection VI.2]{MacLane:1995} points out,
  it is often most convenient to consider a graded object $M$ as a collection 
  of objects $(M_i)_{i\in I}$; this is usually called an \emph{external grading}.
  Depending on the context and the category in which we are working,
  this can be reinterpreted as an \emph{internally graded} object,
  say $\prod_{i \in I} M_i$ or $\bigsqcup_{i \in I} M_i$ or $\oplus_{i \in I} M_i$ etc.
  
  For graded sets we use $\bigsqcup_{i \in I} X_i$, whereas for graded 
  groups the appropriate structure turns out to be $\prod_{i \in I} G_i$.
  As we have already mentioned, for quandles the situation 
  is special, because the decomposition $Q = \bigsqcup_{i \in I} Q_i$ 
  is not simply a disjoint union of quandles $Q_i$:
  in general we have to encode a non-trivial action
  $Q_i \times Q_j \to Q_i$, $(a,b) \mapsto a \ast b$.
\end{remark}

\subsection{Graded extensions}

\begin{definition}
  A graded quandle $Q$ is \emph{connected} (in the graded sense)
  if each set $Q_i$ is a connected component of $Q$.  Likewise,
  a graded covering $p \colon \tilde{Q} \to Q$ is said to be \emph{connected} 
  if each set $\tilde{Q}_i = p^{-1}(Q_i)$ is a connected component of $\tilde{Q}$.
  The covering $p$ is said to be \emph{galois} if, moreover, $\Aut(p)$ 
  acts transitively on the $i$th fibre $p^{-1}(q_i)$ for each $i \in I$.
\end{definition}

\begin{remark}
  Every galois covering $p \colon \tilde{Q} \to Q$ comes with
  the natural action $\Lambda \curvearrowright \tilde{Q}$ of the 
  graded deck transformation group $\Lambda = \Aut(p)$ satisfying 
  the following two axioms:
  \begin{enumerate} 
  \item[(E1)]
    $(\lambda\tilde{x}) \ast \tilde{y} = \lambda (\tilde{x}\ast\tilde{y})$
    and $\tilde{x} \ast (\lambda\tilde{y}) = \tilde{x}\ast\tilde{y}$
    for all $\tilde{x},\tilde{y}\in\tilde{Q}$ and $\lambda\in\Lambda$.
  \item[(E2)]
    $\Lambda_i$ acts freely and transitively on each fibre $p^{-1}(x)$ with $x \in Q_i$.
  \end{enumerate}
  Axiom (E2) then says that $\tilde{Q}_i \to Q_i$ is a 
  \emph{principal} $\Lambda_i$-covering, in the sense 
  that each fibre is a principal $\Lambda_i$-set.
  Notice, however, that we have to consider these
  actions individually over each component $Q_i$; the groups 
  $\Lambda_i$ act independently and may vary for different $i \in I$.
\end{remark}

\begin{definition}
  A \emph{graded extension} $E \colon \Lambda \curvearrowright \tilde{Q} \lto[p] Q$
  of a graded quandle $Q$ by a graded group $\Lambda$ consists of a surjective 
  quandle homomorphism $p \colon \tilde{Q} \to Q$ and a graded group action 
  $\Lambda \curvearrowright \tilde{Q}$ satisfying the axioms (E1) and (E2).
  They entail that $p$ is a quandle covering, and the action of $\Lambda$
  defines an injective homomorphism $\Lambda \to \Aut(p)$ of graded groups.
\end{definition}

\subsection{Universal coverings}

As before we will have to choose base points in order to obtain uniqueness 
properties.  To this end we equip each component with its own base point.  

\begin{definition}[pointed quandles] \label{def:MultiPointedQuandles}
  A \emph{pointed quandle} $(Q,q)$ is a graded quandle 
  $Q = \bigsqcup_{i \in I} Q_i$ with a base point $q_i \in Q_i$ for each $i \in I$.
  In other words, if the partition is seen as a quandle homomorphism
  $\gr \colon Q \to I$, then the choice of base points is 
  a section $q \colon I \to Q$, $\gr \circ q = \id_I$.
  We call $(Q,q)$ \emph{well-pointed} if $q$ specifies one base point 
  in each component, that is, the induced map $\pi_0 \circ q \colon I \to \pi_0(Q)$
  is a bijection between $I$ and the set of connected components of $Q$.

  A homomorphism $\phi \colon (Q,q) \to (Q',q')$ between pointed quandles 
  is a quandle homomorphism $\phi \colon Q \to Q'$ such that $\phi \circ q = q'$.  
  Obviously, $I$-pointed quandles and their homomorphisms 
  form a category, denoted $\Quandles_I^*$.
\end{definition}

\begin{lemma} \label{lem:MultiCoveringQuandle}
  Let $(Q,q)$ be a well-pointed quandle 
  with connected components $(Q_i,q_i)_{i \in I}$.
  Let $\Adj(Q)^\circ$ be the kernel of the 
  group homomorphism $\varepsilon \colon \Adj(Q) \to \Z$ 
  with $\varepsilon(\adj(Q)) = 1$.  For each $i \in I$ we define 
  \[
  \tilde{Q}_i := \bigl\{\; (a,g) \in Q_i \times \Adj(Q)^\circ \;|\; a = q_i^g \;\bigr\}, 
  \quad \tilde{q}_i := (q_i,1) .
  \]
  The disjoint union $\tilde{Q} = \bigsqcup_{i \in I} \tilde{Q}_i$
  becomes a graded quandle with the operations
  \begin{align*}
    (a,g) \ast (b,h) & := \bigl(\; a \ast b,\; g \cdot \adj(a)^{-1} \cdot \adj(b) \;\bigr),  
    \\
    (a,g) \tsa (b,h) & := \bigl(\; a \tsa b,\; g \cdot \adj(a) \cdot \adj(b)^{-1} \;\bigr).
  \end{align*}
  The quandle $\tilde{Q}$ comes with a natural augmentation
  $\tilde{Q} \lto[\rho] \Adj(Q) \lto[\alpha] \Inn(\tilde{Q})$,
  where $\rho(b,h) = \adj(b)$ and $\alpha$ is defined by the action
  \[
  \tilde{Q}_i \times \Adj(Q) \to \tilde{Q}_i \quad\text{with}\quad
  (a,g)^h := \bigl(\; a^h,\; \adj(q_i)^{-\varepsilon(h)} \cdot g h \;\bigr) .
  \]
  The subgroup $\Adj(Q)^\circ$ acts freely and transitively 
  on each $\tilde{Q}_i$.  As a consequence, the connected 
  components of $\tilde{Q}$ are the sets $\tilde{Q}_i$, 
  and so $\tilde{Q}$ is connected in the graded sense.  
  
  The canonical projection $p \colon (\tilde{Q},\tilde{q}) \to (Q,q)$ 
  given by $p(a,g) = a$ is a surjective quandle homomorphism, 
  and equivariant with respect to the action of $\Adj(Q)$.
  \qed
\end{lemma}

\begin{theorem}
  Let $(Q,q)$ be a well-pointed quandle and let $(\tilde{Q},\tilde{q})$ 
  be defined as above. 
  Then the projection $p \colon (\tilde{Q},\tilde{q}) \to (Q,q)$ 
  is the universal quandle covering of $(Q,q)$.  
  \qed
\end{theorem}


The verification of this and the following results in the graded case 
are a straightforward transcription of our previous arguments 
for the non-graded case of connected quandles, and will be omitted.

\subsection{Fundamental group and Galois correspondence} \label{sec:GradedGaloisCorrespondence}

\begin{definition}
  We call $\pi_1(Q,q_i) = \{ g \in \Adj(Q)^\circ \mid q_i^g = q_i \}$ 
  the fundamental group of the quandle $Q$ based at $q_i \in Q$. 
  For a pointed graded quandle $(Q,q)$ we define the graded fundamental 
  group to be the product $\pi_1(Q,q) := \prod_{i \in I} \pi_1(Q,q_i)$. 
\end{definition}

\begin{proposition} \label{prop:CanonicalIdentification2}
  For the universal covering $p \colon (\tilde{Q},\tilde{q}) \to (Q,q)$ as above, 
  we obtain a canonical isomorphism $\phi \colon \pi_1(Q,q) \isoto \Aut(p)$
  of graded groups from the graded left action 
  $\pi_1(Q,q_i) \times \tilde{Q}_i \to \tilde{Q}_i$
  defined by $h \cdot (a,g) = (a,hg)$.
  \qed
\end{proposition}

\begin{proposition}[functoriality]
  Every homomorphism $f \colon (Q,q) \to (Q',q')$ of pointed quandles 
  induces a homomorphism $f_* \colon \pi_1(Q,q) \to \pi_1(Q',q')$ 
  of graded fundamental groups.  We thus obtain a functor 
  $\pi_1 \colon \Quandles_I^* \to \Groups_I$ from the category 
  of $I$-pointed quandles to the category of $I$-graded groups.
  \qed
\end{proposition}

\begin{proposition}[lifting criterion] \label{prop:MultiLiftingCriterion}
  Let $p \colon (\tilde{Q},\tilde{q}) \to (Q,q)$ be a quandle covering
  and let $f \colon (X,x) \to (Q,q)$ be a quandle homomorphism from
  a well-pointed quandle $(X,x)$ to the base quandle $(Q,q)$.  Then 
  there exists a lifting $\tilde{f} \colon (X,x) \to (\tilde{Q},\tilde{q})$,
  $p \circ \tilde{f} = f$, if and only if 
  $f_* \pi_1(X,x) \subset p_* \pi_1(\tilde{Q},\tilde{q})$. 
  In this case the lifting $\tilde{f}$ is unique.
  \qed
\end{proposition}

\begin{theorem}[Galois correspondence for well-pointed coverings]
  Let $(Q,q)$ be a well-pointed quandle indexed by some set $I$.
  The canonical functors $\Coverings_I(Q,q) \to \Subgroups_I(\pi_1(Q,q))$ 
  and $\Subgroups_I(\pi_1(Q,q)) \to \Coverings_I(Q,q)$ establish
  a natural equivalence between the category of well-pointed coverings 
  of $(Q,q)$ and the category of graded subgroups of $\pi_1(Q,q)$.
  \qed
\end{theorem}


\begin{theorem}[Galois correspondence for general coverings]
  Let $(Q,q)$ be a well-pointed quandle indexed by some set $I$.
  The canonical functors $\Coverings(Q) \to \Actions_I(\pi_1(Q,q))$ 
  and $\Actions_I(\pi_1(Q,q)) \to \Coverings(Q)$ establish
  a natural equivalence between the category of coverings of $(Q,q)$ 
  and the category of graded actions of $\pi_1(Q,q)$.
  \qed
\end{theorem}

\begin{theorem} \label{thm:ExtensionNonConnectedCorrespondence}
  Let $(Q,q)$ be a well-pointed quandle indexed by some set $I$,
  and let $\Lambda$ be a graded group.  There exists a natural 
  bijection $\Ext(Q,\Lambda) \cong \Hom(\pi_1(Q,q),\Lambda)$.
  If $\Lambda$ is a graded abelian group, or more generally a graded 
  module over some ring $R$, then both objects carry natural $R$-module 
  structures and the natural bijection is a graded $R$-module isomorphism.
  \qed
\end{theorem}

\begin{example} \label{exm:QmnCovering}
  The covering theory of non-connected quandles allows us to complete 
  the discussion of the quandle $Q_{m,n} = \Z_m \sqcup \Z_n$ 
  begun in Example \ref{exm:DoubleCycle}.  
  From Proposition \ref{prop:HeisenbergGroup} we deduce that 
  \[
  \Adj(Q_{m,n})^\circ = \left\{ 
    \left(\begin{smallmatrix} 1 & -s & t \\ 0 & 1 & +s \\ 0 & 0 & 1
      \end{smallmatrix}\right) \mathrel{\big|} s \in \Z, t \in \Z_\ell \right\}
  \subset H/\gen{z^\ell} .
  \]
  The shown matrix acts as $a \mapsto a+s$ on $a \in \Z_m$, 
  and as $b \mapsto b-s$ on $b \in \Z_n$, which entails 
  $\pi_1(Q,a) = m\Z \times \Z_\ell$ and $\pi_1(Q,b) = n\Z \times \Z_\ell$.
  The universal covering $p \colon \tilde{Q} \onto Q_{m,n}$
  can be constructed as in Lemma \ref{lem:MultiCoveringQuandle}.
  After some calculation this leads to $Q_\ell = A \sqcup B$, 
  where $A$ and $B$ are copies of $\Z \times \Z_\ell$ 
  with $\ell = \gcd(m,n)$, and the quandle structure 
  \[
  (a,a') \ast (b,b') = \begin{cases}
    ( a, a'+b-a ) & \text{if } (a,a'),(b,b') \in A \text{ or if } (a,a'),(b,b') \in B , \\
    ( a+1, a'-b ) & \text{otherwise.}
  \end{cases}
  \]

  The quandle $Q_\ell$ has two connected components, 
  $A$ and $B$, so it is connected in the graded sense.
  The projection $p \colon Q_\ell \onto Q_{m,n}$ 
  is defined by $A \to \Z_m$, $(a,a') \mapsto a \bmod{m}$,
  and $B \to \Z_n$, $(b,b') \mapsto b \bmod{n}$.
  This is the universal covering of $Q_{m,n}$,
  and any other covering that is connected in the graded sense
  is obtained by quotienting out some graded subgroup 
  of $\Aut(p) \cong \pi_1(Q,a) \times \pi_1(Q,b)$.

  Notice that in the special case $\ell=1$ we obtain 
  the obvious covering $Q_{0,0} \onto Q_{m,n}$, 
  but even in this toy example the general case would be 
  difficult to analyze without the classification theorem.
\end{example}

\begin{Notes}
  \newcommand{\heis}[3]{\left(\begin{smallmatrix} 
        1 & #1 & #3 \\ 0 & 1 & #2 \\ 0 & 0 & 1 \end{smallmatrix}\right)}
  Here are the details: according to
  Proposition \ref{prop:HeisenbergGroup} we have 
  \[
  \adj(a) = \heis{1}{0}{a} \text{ for } a \in \Z_m, \text{ and }
  \adj(b) = \heis{0}{1}{-b} \text{ for } b \in \Z_n .
  \]
  
  As base points we choose $x =0 \in \Z_m$ and $y = 0 \in \Z_n$ respectively.

  Elements of the group $G = \Adj(Q_{m,n})^\circ$ 
  can be conveniently parametrized by 
  \[
  g(s,t) = \heis{-s}{+s}{t-s^2} \text{ with } s \in \Z, t \in \Z_\ell .
  \] 
  
  The quandle operations on $G$ over $\Z_m$ can then be calculated as:
  \begin{align*}
    g(s,t) \ast a & = \adj(x)^{-1} g(s,t) \adj(a) 
    = \heis{-1}{0}{0} \heis{-s}{+s}{t-s^2} \heis{1}{0}{a}
    \\ 
    & = \heis{-s}{s}{t-s^2+a-s} = g(s,t+a-s) 
    \\
    g(s,t) \ast b & = \adj(x)^{-1} g(s,t) \adj(b) 
    = \heis{-1}{0}{0} \heis{-s}{+s}{t-s^2} \heis{0}{1}{-b}
    \\
    & = \heis{-s-1}{s+1}{t-(s+1)^2-b} = g(s+1,t-b)
  \end{align*}

  We can also parametrize elements of the group $G$ by 
  \[
  h(s,t) = \heis{+s}{-s}{-t} \text{ with } s \in \Z, t \in \Z_\ell .
  \] 
  
  For the quandle operations on $G$ over $\Z_n$ we then find:
  \begin{align*}
    h(s,t) \ast a & = \adj(y)^{-1} h(s,t) \adj(a) 
    = \heis{0}{-1}{0} \heis{+s}{-s}{-t} \heis{1}{0}{a}
    \\
    & = \heis{s+1}{-s-1}{-t+a} = h(s+1,t-a) 
    \\
    h(s,t) \ast b & = \adj(y)^{-1} h(s,t) \adj(b) 
    = \heis{0}{-1}{0} \heis{+s}{-s}{-t} \heis{0}{1}{-b}
    \\
    & = \heis{+s}{-s}{-t-b+s} = h(s,t+b-s)
  \end{align*}
  
  This justifies the above formulae for $Q_\ell$.
\end{Notes}

\subsection{Application to link quandles} \label{sub:LinkQuandles1}


Given an $n$-component link $K = K_1 \sqcup \dots \sqcup K_n \subset \S^3$, 
we choose a base point $q_K^i \in Q_K$ for each link component $K_i$.
The adjoint group $\Adj(Q_K)$ is isomorphic to the fundamental group 
$\pi_K = \pi_1(\S^3 \minus K)$, and each element $q_K^i$ maps 
to a meridian $m_K^i = \adj(q_K^i) \in \pi_K$.  
We denote by $\ell_K^i \in \pi_K$ the corresponding longitude.

The universal covering $p \colon \tilde{Q}_K \onto Q_K$ can 
formally be constructed as in Lemma \ref{lem:MultiCoveringQuandle}.
Its geometric interpretation has been studied in \cite{Eisermann:2003}
in terms of quandle homology $H_2(Q_K)$ and orientation classes 
$[K_i] \in H_2(Q_K)$.  We are now in position to go one step 
further and determine the fundamental group: 

\begin{theorem} \label{thm:LinkQuandleFundamentalGroup}
  Over each component $Q_K^i \subset Q_K$ the automorphism group 
  of the universal covering $p \colon \tilde{Q}_K \onto Q_K$
  is given by $\Aut(p)_i = \pi_1(Q_K,q_K^i)= \gen{\ell_K^i}$.  
  For the graded fundamental group this means that 
  $\pi_1(Q_K,q_K) = \Aut(p) = \prod_{i=1}^n \gen{\ell_K^i}$.
\end{theorem}

\begin{proof}
  Fixing a link component $K_i$, we can construct a long link 
  $L \subset \R^3$ by opening $K_i$ while leaving all other components closed.
  This is the same as removing from the pair $(\S^3,K)$ a point on $K_i$
  so as to obtain the pair $(\R^3,L)$.  In particular, the correspondence 
  $(\S^3,K) \leftrightarrow (\R^3,L)$ is well-defined when we pass to isotopy classes.  
  The associated quandle $Q_L$ has two distinguished elements $q_L$ and $q_L^*$, 
  corresponding to the beginning and the end of the open component, respectively.
  The natural quandle homomorphism $p_i \colon Q_L \onto Q_K$ is the 
  quotient obtained by identifying $q_L$ and $q_L^*$, both being mapped 
  to $q_K^i = p_i(q_L) = p_i(q_L^*)$.

  While $p_i \colon Q_L \onto Q_K$ is in general not 
  an isomorphism between the quandles $Q_L$ and $Q_K$, the induced map 
  $\Adj(p_i) \colon \Adj(Q_L) \onto \Adj(Q_K)$ is always an isomorphism 
  between the adjoint groups $\Adj(Q_L) = \pi_L = \pi_1(\R^3 \minus L)$
  and $\Adj(Q_K) = \pi_K = \pi_1(\S^3 \minus K)$.  In particular, 
  this implies that $p_i \colon (Q_L,q_L) \onto (Q_K,q_K)$ is 
  a quandle covering, and an isomorphism over all components except $Q_K^i$.

  Let $\hat{p}_i \colon (\hat{Q}_K,\hat{q}_K^i) \onto (Q_K,q_K^i)$ 
  be the covering that is universal over $Q_K^i$ and an isomorphism
  over all other components.  Then one can construct an isomorphism
  $(Q_L,q_L) \isoto (\hat{Q}_K,\hat{q}_K^i)$ of quandle coverings 
  over $(Q_K,q_K^i)$ as in \cite[Theorem 30]{Eisermann:2003}.
  In particular, we obtain a canonical group isomorphism
  $\Aut(p_i) \cong \pi_1(Q_K,q_k^i)$ as in 
  Proposition \ref{prop:CanonicalIdentification2}.

  The longitude $\ell_K^i$ satisfies $(q_K^i)^{\ell_K^i} = q_K^i$,
  so $\ell_K^i \in \pi_1(Q_K,q_k^i)$.  Moreover, $(q_L)^{\ell_K^i} = q_L^*$, 
  so the quotient of $Q_L$ by the subgroup $\gen{\ell_K^i} \subset \Aut(p_i)$ 
  yields $\gen{\ell_K^i} \backslash Q_L = Q_K$.  
  To see this, notice that we have a canonical projection 
  $\gen{\ell_K^i} \backslash Q_L \onto Q_K$ 
  as a quotient of the covering $Q_L \onto Q_K$.
  Inversely, we have a canonical map $Q_K \onto \gen{\ell_K^i} \backslash Q_L$
  by the universal property of the quotient $Q_K = Q_L/(q_L=q_L^*)$.
  We conclude that $\Aut(p)_i = \Aut(p_i) = \gen{\ell_K^i}$. 
\end{proof}

A link component $K_i \subset K$ is called \emph{trivial}, if there exists 
an embedded disk $D \subset \S^3$ with $K_i = K \cap D = \partial D$.
Using the Loop Theorem of Papakyriakopoulos \cite{Papakyriakopoulos:1957}
we conclude:

\begin{corollary} \label{cor:TrivialComponents}
  For a link $K \subset \S^3$ the following assertions are equivalent:
  \begin{enumerate}
  \item
    The link component $K_i \subset K$ is trivial.
  \item
    The fundamental group $\pi_1(Q_K,q_K^i)$ 
    is trivial.
  \item
    The longitude $\ell_K^i \in \pi_K$ 
    is trivial.
  \end{enumerate}
  Conversely, if the link component $K_i$ is non-trivial,
  then the fundamental group $\pi_1(Q_K,q_K^i)$ of the quandle $Q_K$ 
  based at $q_K^i$ is freely generated by the longitude $\ell_K^i$.
\end{corollary}

\begin{proof}
  The implications $(1) \Rightarrow (2) \Leftrightarrow (3)$ 
  follow from $\pi_1(Q_K,q_K^i) = \gen{\ell_K^i}$,
  established in the previous theorem, while $(3) \Rightarrow (1)$ 
  is a consequence of the Loop Theorem \cite{Papakyriakopoulos:1957}.
  If $K$ is non-trivial, then $\ell_K^i$ is of infinite order,
  and thus freely generates $\pi_1(Q_K,q_K^i)$.
\end{proof}

\section{Fundamental groupoid of a quandle} \label{sec:FundamentalGroupoid}

As in the case of topological spaces, the choice of a base point $q \in Q$
in the definition of $\pi_1(Q,q)$ focuses on one connected component 
and neglects the others.  If we do not want to fix base points, 
then the fundamental groupoid is the appropriate tool.  
(See Spanier \cite[\textsection 1.7]{Spanier:1981},
Brown \cite[chap.\ 9]{Brown:1988}, and May \cite[chap.\ 3]{May:1999}).  
We shall expound this idea in the present section because it explains 
the striking similarity between quandles and topological spaces.

\subsection{Groupoids}

We recall that a \emph{groupoid} is a small category 
in which each morphism is an isomorphism.
In geometric language one considers its objects as ``points'' 
$a,b,\dots$ and it morphisms $a \to b$ as ``paths'' 
(or, more frequently, equivalence classes of paths).

\begin{example}
  The classical example is the fundamental groupoid $\Pi(X)$ 
  of a topological space $X$: this is the category whose objects 
  are the points $x \in X$ and whose morphisms $x \to y$ 
  are the homotopy classes of paths from $x$ to $y$. 
  There exists a morphism $x \to y$ if and only if 
  $x$ and $y$ belong to the same path-component.
  The group of automorphisms of an object $x$ is exactly 
  the fundamental group $\pi_1(X,x)$ of $X$ based at $x$.
\end{example}

\begin{example}
  Consider a set $Q$ with a group action $Q \times G \to Q$,
  denoted by $(a,g) \mapsto a^g$.  We can then define the groupoid
  \[
  \Pi(Q,G) := \{\; (a,g,b) \in Q \times G \times Q \;|\; a^g = b \;\} .
  \]
  Here the objects are given by elements $a \in Q$, 
  and the morphisms from $a$ to $b$ are the triples $(a,g,b) \in \Pi(Q,G)$.
  Their composition is defined by $(a,g,b) \circ (b,h,c) = (a,gh,c)$.
  There exists a morphism $a \to b$ if and only if 
  $a$ and $b$ belong to the same $G$-orbit.
  The group of automorphisms of an object $a$
  is exactly the stabilizer of $a$ in $G$.
\end{example}

\begin{definition}
  For a quandle $Q$ we call $\Pi(Q,\Adj(Q)^\circ)$
  the \emph{fundamental groupoid} of $Q$.
\end{definition}

\begin{remark}[connected components]
  Already Joyce noticed some analogy between 
  quandles and topological spaces when he introduced
  the terminology ``connected component'' of $Q$ to signify 
  an orbit with respect to the inner automorphism group $\Inn(Q)$.  
  (This was probably motivated by the example of symmetric spaces, 
  where both notions of connectedness coincide, 
  see Remark \ref{rem:ConnectedSymmetricSpaces}.)
  This turned out to be a very fortunate and intuitive wording, 
  and connectedness arguments have played a crucial r\^ole
  for all subsequent investigations of quandles.
  The connected components of the quandle $Q$ are precisely 
  those of the groupoid $\Pi(Q,\Adj(Q)^\circ)$,
  see Remark \ref{rem:CommutatorOrbits}.  
\end{remark}

\begin{remark}[fundamental group]
  According to the previous remark one can partition 
  a quandle $Q$ into the set $\pi_0(Q)$ of connected components,
  and with a little bit of na\"ivet\'e one could wonder 
  what the fundamental group $\pi_1(Q,q)$ should be. 
  In the above groupoid we recover the fundamental group
  $\pi_1(Q,q) = \{ g \in \Adj(Q)^\circ \mid q^g = q \}$ based at $q \in Q$
  as the group of automorphisms of the object $q$ in the category $\Pi(Q,\Adj(Q)^\circ)$.
  For base points $q,q'$ in the same component of $Q$, 
  these groups are isomorphic by a conjugation in $\Pi(Q,\Adj(Q)^\circ)$.
  As usual this isomorphism is not unique, unless $\pi_1(Q,q)$ is abelian.
\end{remark}

\begin{remark}[coverings]
  There exists an extensive literature on groupoids 
  because they generalize and simplify recurring arguments
  in seemingly different situations, notably in diverse Galois theories, 
  just as in our setting of coverings and fundamental groups of quandles.
  The universal covering quandle $(\tilde{Q},\tilde{q})$ 
  constructed in Lemmas \ref{lem:UniversalCovering} and
  \ref{lem:MultiCoveringQuandle} reappears here as 
  the set of paths based at $q$ (with arbitrary endpoint).
  This is exactly the path fibration used to construct 
  the universal covering of a topological space,
  or more generally of a groupoid.  We refer to the 
  excellent introduction of May \cite[chap.\ 3]{May:1999}.
\end{remark}

In conclusion, the ``generic part'' of quandle covering theory 
can be recast in the general language of groupoid coverings.
The initial problem, however, is to construct the appropriate groupoid.
Several groupoid structures are imaginable, 
and one cannot easily guess the appropriate one: a priori 
one can choose many groups acting on $Q$, for example 
$\Adj(Q)$, $\Aut(Q)$, $\Inn(Q)$, or $\Inn(Q)^\circ$, but only 
the choice $\Adj(Q)^\circ$ yields the groupoid that is dual to quandle coverings.  
The difficulty is thus resolved by first analyzing coverings, 
which seem to be the more natural notion.

It should also be noted that the unifying concept of groupoids 
does not cover the whole theory of quandle coverings.
Besides its ``generic'' aspects, the latter also has its distinctive 
``non-standard'' features.  These have been pointed out 
in \secref{sub:TournantsDangereux} and merit special attention.
This is why we have preferred to present all constructions in detail.

\subsection{Combinatorial homotopy} 

For future reference, let us give another derivation how 
the group $\Adj(Q)^\circ$ and the associated groupoid 
$\Pi(Q,\Adj(Q)^\circ)$ appear naturally --- as the groupoid 
of combinatorial paths modulo combinatorial homotopy.

\begin{definition}
  Let $Q$ be a quandle. Consider the graph $\Gamma$ with vertices $q \in Q$ 
  and edges $a \lto[b] c$ for each triple $a,b,c \in Q$ with $a \ast b = c$.
  A \emph{combinatorial path} from $q$ to $q'$ in $\Gamma$ is 
  a sequence of vertices $q=a_0,a_1,\dots,a_{n-1},a_n=q' \in Q$ 
  and arrows $b_1^{\varepsilon_1},\dots,b_n^{\varepsilon_n}$
  with $b_i \in Q$ and $\varepsilon_i \in \{ \pm1 \}$ for all $i$,
  such that $a_{i-1} \ast b_i = a_i$ for $\varepsilon_i = +1$
  and $a_{i-1} \tsa b_i = a_i$ for $\varepsilon_i = -1$.
  The sign $\varepsilon_i$ is just a convenient way 
  to denote the orientation of the $i$th arrow:
  \[
  (a \lto[b^+] a \ast b) = (a \lto[b] a \ast b)
  \quad\text{and}\quad
  (a \lto[b^-] a \tsa b) = (a \lot[b] a \tsa b) . 
  \]

  Let $P(Q)$ be the category having as objects the elements $q \in Q$
  and as morphisms from $q$ to $q'$ the set of combinatorial paths 
  from $q$ to $q'$.  Composition is given by juxtaposition:
  \[
  (a_0 \to \cdots \to a_m) \circ (a_m \to \cdots \to a_n) 
  = (a_0 \to \cdots \to a_m \to \cdots \to a_n) .
  \]

  Two combinatorial paths are \emph{homotopic} if they can be transformed
  one into the other by a sequence of the following 
  local moves and their inverses:
  \begin{enumerate}
  \item[(H1)]
    $a \lto[a] a$ is replaced by $a$, or $a \lot[a] a$ is replaced by $a$.
  \item[(H2)]
    $a \lto[b] a \ast b \lot[b] a$ is replaced by $a$, or $a \lot[b] a \tsa b \lto[b] a$ is replaced by $a$.
  \item[(H3)]
    $a \lto[b] a \ast b \lto[c] (a \ast b) \ast c$ is replaced by 
    $a \lto[c] a \ast c \lto[b \ast c] (a \ast c) \ast (b \ast c)$.
  \end{enumerate}
  
  We denote by $\Pi(Q)$ the quotient category having as objects 
  the elements $q \in Q$ and as morphisms from $q$ to $q'$ 
  the set of homotopy classes of combinatorial paths from $q$ to $q'$.
\end{definition}

\begin{figure}[hbtp]
  \centering
  \ifpdf\input{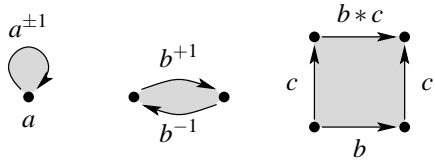}\else\input{qhomotopy.pstex_t}\fi
  \caption{Elementary homotopies for paths in $P(Q)$}
  \label{fig:ElementaryHomotopies}
\end{figure} 

\begin{proposition}
  The category $\Pi(Q)$ is a groupoid, 
  that is, every morphism is invertible.
  Moreover, there exists a natural isomorphism 
  $\Phi \colon \Pi(Q) \isoto \Pi(Q,\Adj(Q)^\circ)$,
  given by 
  \[
  [a_0 \lto[b_1^{\varepsilon_1}] \cdots \lto[b_n^{\varepsilon_n}] a_n] \mapsto (a_0, g, a_n)
  \quad\text{with}\quad 
  g = a_0^{-\sum_i \varepsilon_i} b_1^{\varepsilon_1} \cdots b_n^{\varepsilon_n} \in \Adj(Q)^\circ
  \]
\end{proposition}

\begin{proof}
  The homotopy relation (H2) above ensures that $\Pi(Q)$ is a groupoid.
  It is straightforward to verify that the map $\Phi$ is well-defined:
  a homotopy (H1) does not change the element $g \in \Adj(Q)^\circ$
  due to the normalization with $a_0^{-\sum_i \varepsilon_i}$.
  A homotopy (H2) translates to $b^{\pm} b^{\mp} = 1$.
  A homotopy (H3) translates to one of the defining relations
  $c \cdot (b \ast c) = b \cdot c$ of the adjoint group $\Adj(Q)$.
  By construction, $\Phi$ sends composition in $\Pi(Q)$ to 
  composition in $\Pi(Q,\Adj(Q)^\circ)$, so it is a functor.
  Obviously $\Phi$ is a bijection on objects $q \in Q$, and it 
  is easy to see that it is also a bijection on morphisms.
\end{proof}

\subsection{Classifying spaces} \label{sub:TopologicalRealization}

As usual, combinatorial paths and combinatorial homotopy 
can be realized by a suitable topological space $K$:
it suffices to take the graph $\Gamma$ as $1$-skeleton
and to glue a $2$-cell for each relation of type (H1) and (H3).
(Relation (H2) is automatic, since both $a \lto[b^{+}] a \ast b$ and
$a \lot[b^{-}] a \ast b$ are actually represented by the same edge.)  
This ensures that $\Pi(Q)$ is the edge-path groupoid of 
the resulting $2$-dimensional (cubical) complex $K$;
see Spanier \cite[\textsection 3.6]{Spanier:1981} for the simplicial case.

When we go back to the sources of quandle and rack cohomology,
we thus rediscover yet another approach to the fundamental group $\pi_1(Q,q)$ 
of a quandle $Q$, which is entirely topological and has the merit 
to open up the way to a full-fledged homotopy theory:
Fenn, Rourke, and Sanderson \cite{FennRourkeSanderson:1995}
constructed a classifying space $BX$ for a rack $X$,
which allowed them to define (co)homology and homotopy groups for each rack.
Their construction can be adapted to quandles $Q$, so that the
resulting classifying space $BQ$ is a topological model for 
quandle (co)homology $H_*(Q) = H_*(BQ)$ and $H^*(Q) = H^*(BQ)$. 
Our construction of $K$ corresponds precisely to the $2$-skeleton of $BQ$. 

The homotopy groups $\pi_n(BQ)$ have not yet played a r\^ole
in the study of quandles.  It turns out, however, that our algebraic
fundamental group $\pi_1(Q,q)$ coincides with the fundamental group
of the classifying space, $\pi_1(BQ,q)$, at least in the case of a connected quandle.
Starting from the algebraic notion of quandle covering, we thus recover 
and remotivate the topological construction of Fenn, Rourke, and Sanderson.  

\subsection{Theft or honest toil?}


In order to define the fundamental group of a quandle $Q$,
one could thus take its classifying space $BQ$ and set $\pi_1(Q,q) := \pi_1(BQ,q)$.
Does this mean that we could entirely replace the 
algebraic approach by its topological counterpart?  
Two arguments suggest that this is not so:

\begin{itemize}
\item
  Even with an independent \emph{topological} definition of $\pi_1(Q,q)$,
  one would still have to prove that the \emph{algebraic} covering 
  theory of quandles behaves the way it does, and in particular
  is governed by the fundamental group so defined,
  in order to establish and exploit their relationship.
\item
  Quandle coverings differ in some crucial details from
  topological coverings, which means that both theories
  cannot be equivalent in any superficial way.
  It is thus justified and illuminating
  to develop the algebraic theory independently.
\end{itemize}

In conclusion it appears that algebraic coverings are interesting in their own
right, and that the algebraic and the topological viewpoint are complementary.


\section{Extensions and cohomology} \label{sec:ExtensionCohomology}

Our goal in this final section is to establish a correspondence between 
quandle extensions $E \colon \Lambda \curvearrowright \tilde{Q} \to Q$ 
and elements of the second cohomology group $H^2(Q,\Lambda)$.
For abelian groups $\Lambda$ this is classical for group extensions 
(see for example Mac\,Lane \cite[\textsection{IV.4}]{MacLane:1995} or 
Brown \cite[\textsection{IV.3}]{Brown:1994}) and has previously been translated 
to quandle extensions. 
This correspondence has to be generalized in two directions 
in order to apply to our general setting:

\begin{itemize}
\item
  The usual formulation is most appealing for abelian 
  groups $\Lambda$, and has been independently developed 
  in \cite{CarterEtAl:2003} and \cite{Eisermann:2003}.  
  For general galois coverings and extension, however, 
  the coefficient group $\Lambda$ can be non-abelian.
\item
  For non-connected quandles the notion of extension 
  must be refined in the graded sense, because different
  components have to be treated individually.
  The corresponding cohomology theory $H^2(Q,\Lambda)$ deals 
  with a graded quandle $Q$ and a graded group $\Lambda$, 
  both indexed by some fixed set $I$.
\end{itemize}

For racks such a non-abelian cohomology theory has previously been proposed
by N.\,Andruskiewitsch and M.\,Gra\~na \cite[\textsection4]{AndruskiewitschGrana:2003}. 
In view of knot invariants, this has been adapted to a non-abelian quandle 
cohomology in \cite{CarterEtAl:2005}.  We will complete this approach 
by establishing a natural bijection between $\Ext(Q,\Lambda)$ and 
$H^2(Q,\Lambda)$ in the non-abelian graded setting, which specializes 
to the previous formulation in the abelian non-graded case.

\subsection{Non-abelian graded quandle cohomology} \label{sub:QuandleCohomology}

Let $Q = \bigsqcup_{i \in I} Q_i$ be a graded quandle and 
let $\Lambda = \prod_{i \in I} \Lambda_i$ be a graded group.
We do not assume $\Lambda$ to be abelian and will thus use multiplicative notation.

\begin{remark}
  The first cohomology $H^1(Q,\Lambda)$ consists of all graded maps 
  $g \colon Q \to \Lambda$ with $g(Q_i) \subset \Lambda_i$,
  such that $g(a) = g(a * b)$ for all $a,b \in Q$.
  These are the class functions, i.e.\ functions
  that are constant on each connected component of $Q$.
  Notice that the grading of $Q = \bigsqcup_{i \in I} Q_i$ turns $H^1$
  into a graded group, $H^1(Q,\Lambda) = \prod_{i \in I} H^1(Q,\Lambda)_i$.
  If $Q$ is graded connected, i.e.\ each $Q_i$ is a connected component of $Q$, 
  then $H^1(Q,\Lambda) = \prod_{i \in I} \Lambda_i = \Lambda$.
\end{remark}

In order to define $H^2(Q,\Lambda)$ we proceed as follows.

\begin{definition}
  The grading of the quandle $Q = \bigsqcup_{i \in I} Q_i$ induces a grading 
  of the product $Q \times Q = \bigsqcup_{i \in I} Q_i \times Q$.
  A \emph{$2$-cochain} is a graded map $f \colon Q \times Q \to \Lambda$
  with $f(Q_i \times Q) \subset \Lambda_i$,  such that $f(a,a) = 1$ 
  for all $a \in Q$.  We say that $f$ is a \emph{$2$-cocycle} if 
  \[ 
  f( a, b ) f( a \ast b, c ) = f(a,c) f( a \ast c, b \ast c )
  \quad\text{for all $a,b,c \in Q$.}
  \]
  We denote by $Z^2(Q,\Lambda)$ the set of $2$-cocycles.
  We say that two cocycles $f,f' \in Z^2(Q,\Lambda)$ are \emph{cohomologous} 
  if there exists a graded map $g \colon Q \to \Lambda$ 
  with $g(Q_i) \subset \Lambda_i$ such that 
  \[
  f(a,b) = g(a)^{-1} f'(a,b) g(a \ast b)  \quad\text{for all $a,b \in Q$.}
  \]
  This is an equivalence relation on $Z^2(Q,\Lambda)$, 
  and we denote by $H^2(Q,\Lambda)$ the quotient set.  
  Its elements are \emph{cohomology classes} $[f]$ 
  of $2$-cocycles $f \in Z^2(Q,\Lambda)$.
\end{definition}

\begin{remark}
  Notice that the set $C^2$ of $2$-cochains decomposes as $C^2 = \prod_{i \in I} C^2_i$
  where $C^2_i$ consists of maps $f_i \colon Q_i \times Q \to \Lambda_i$.
  Likewise, we obtain $Z^2 = \prod_{i \in I} Z^2_i$ and $H^2 = \prod_{i \in I} H^2_i$.
  
  In the case where $\Lambda$ is an abelian group, or more generally
  a module over some ring $R$, one can define in every degree $n \in \N$ 
  an $R$-module $C^n(Q,\Lambda)$ of quandle $n$-cochains with values in $\Lambda$,
  together with $R$-linear maps $\delta_n \colon C^n \to C^{n+1}$ satisfying 
  $\delta_{n} \delta_{n-1} = 0$.  Such a cochain complex allows us, as usual, 
  to define the submodule of $n$-cocycles $Z^n = \Ker(\delta_n)$ and 
  its submodule of $n$-coboundaries $B^n = \Im(\delta_{n-1})$, and 
  finally the cohomology $H^n = Z^n/B^n$ as their quotient module.
  This construction respects the $I$-grading, and so 
  cochains $C^n = \prod_{i \in I} C^n_i$, cocycles
  $Z^n = \prod_{i \in I} Z^n_i$, coboundaries $B^n = \prod_{i \in I} B^n_i$, 
  and finally cohomology $H^n = \prod_{i \in I} H^n_i$ are $I$-graded modules.

  In the non-abelian case we content ourselves with $H^1$ and $H^2$.
  Notice that $H^1$ can be given a group structure by point-wise multiplication.  
  For $H^2$ pointwise multiplication works if $\Lambda$ is abelian, 
  but it fails in the non-abelian case.  This means that the quotient 
  $H^2(Q,\Lambda)$ is in general only a set.  It has nonetheless 
  a canonical base point, namely the class $[1]$ of the trivial $2$-cocycle 
  $Q \times Q \to \{1\}$, which plays the r\^ole of the neutral element.
\end{remark}

\begin{remark}[functoriality in $Q$]
  Every graded quandle homomorphism $\phi \colon Q' \to Q$ induces 
  a natural graded map $\phi^* \colon H^2(Q,\Lambda) \to H^2(Q',\Lambda)$
  mapping the trivial class to the trivial class.
  More explicitly, $\phi^*$ sends $[f]$ to $[\phi^* f]$,
  where $f \in Z^2(Q,\Lambda)$ is mapped to $\phi^* f \in Z^2(Q,\Lambda)$ 
  defined by $(\phi^* f)(a',b') = f( \phi(a'), \phi(b')$.
\end{remark}

\begin{remark}[functoriality in $\Lambda$]
  Every graded group homomorphism $h \colon \Lambda \to \Lambda'$ induces 
  a natural graded map $h_* \colon H^2(Q,\Lambda) \to H^2(Q,\Lambda')$
  mapping the trivial class to the trivial class.
  More explicitly, $\phi_*$ sends $[f]$ to $[\phi f]$,
  defined by composing $f \colon Q \times Q \to \Lambda$ with the 
  group homomorphism $\phi \colon \Lambda \to \Lambda'$.
\end{remark}

\subsection{Classification of extensions} \label{sub:CentrExtClass}

It is a classical result of group cohomology that 
central extensions of a group $G$ with kernel $\Lambda$ 
are classified by the second cohomology group $H^2(G,\Lambda)$,
see for example Brown \cite[\textsection{IV.3}]{Brown:1994},
or Mac\,Lane \cite[\textsection{IV.4}]{MacLane:1995}.
We will now prove that an analogous theorem holds 
for quandles and their non-abelian graded extensions.

\begin{lemma} \label{lem:CentrExtCocycle}
  Let $E \colon \Lambda \curvearrowright \tilde{Q} \to Q$ 
  be a graded extension of a graded quandle $Q$ by a graded group $\Lambda$.
  Each set-theoretic section $s \colon Q \to \tilde{Q}$ 
  defines a unique graded map $f \colon Q \times Q \to \Lambda$ 
  such that $s(a) \ast s(b) = f(a,b) \cdot s(a\ast b)$.  
  This map $f$ is a quandle $2$-cocycle; it measures the 
  failure of the section $s$ to be a quandle homomorphism.
  Furthermore, if $s' \colon Q\to\tilde{Q}$ is another section,
  then the associated quandle $2$-cocycle $f'$ is homologous to $f$.
  In this way each extension $E$ determines a cohomology class 
  $\Phi(E) := [f] \in H^2(Q,\Lambda)$.
\end{lemma}

\begin{proof}
  Since the action of $\Lambda_i$ is free and transitive 
  on each fibre $p^{-1}(a)$ with $a \in Q_i$, 
  the above equation uniquely defines the map $f$. 
  Idempotency of $\tilde{Q}$ implies $f(a,a)=0$, and 
  self-distributivity implies the cocycle condition:
  \begin{alignat*}{3}
    & \big[ s(a) \ast s(b) \big] \ast s(c) && \quad = \quad
    f(a,b) f(a\ast b,c) & ~ & s\big[ (a\ast b)\ast c \big] \qquad\text{and}
    \\
    & \big[ s(a) \ast s(c) \big] \ast \big[ s(b) \ast s(c) ] && \quad = \quad
    f(a,c) f(a\ast c,b\ast c) & ~ & s\big[ (a\ast c)\ast (b\ast c) \big] .
  \end{alignat*}
  Since both terms are equal, we obtain
  $f(a,b) f(a\ast b,c) = f(a,c) f(a\ast c,b\ast c)$, 
  as desired, which means that $f$ is a $2$-cocycle.
  If $s'$ is another section, then there exists a graded map
  $g \colon Q \to \Lambda$ with $s'(a) = g(a) s(a)$.
  The defining relation $s'(a) \ast s'(b) = f'(a,b) s'(a\ast b)$
  thus becomes $g(a) s(a) \ast g(b) s(b) = f'(a,b) g(a\ast b) s(a\ast b)$.
  Comparing this to $s(a) \ast s(b) = f(a,b) s(a\ast b)$ we find
  that $f(a,b) = g(a)^{-1} f'(a,b) g(a \ast b)$, which means that
  $f$ and $f'$ are cohomologous.  In other words, the cohomology 
  class $[f]$ is independent of the chosen section $s$,
  and hence characteristic of the extension $E$.
\end{proof}

Conversely, we can associate with each quandle $2$-cohomology class
$[f] \in H^2(Q,\Lambda)$ an extension of $Q$ by $\Lambda$:


\begin{theorem} \label{thm:ExtensionCohomologyCorrespondence}
  Let $Q$ be a graded quandle and let $\Lambda$ be a graded group.
  For each extension $E \colon \Lambda \curvearrowright \tilde{Q} \to Q$
  let $\Phi(E)$ be the associated cohomology class in $H^2(Q,\Lambda)$.
  This map induces a natural bijection 
  $\Phi \colon \Ext(Q,\Lambda) \cong H^2(Q,\Lambda)$.
  If $\Lambda$ is an abelian group, or more generally 
  a module over some ring $R$, then $\Ext(Q,\Lambda)$ and 
  $H^2(Q,\Lambda)$ carry each a natural $R$-module structure,
  and $\Phi$ is an isomorphism of $R$-modules.
\end{theorem}

\begin{proof}
  We first note that $\Phi$ is well-defined 
  on equivalence classes of extensions.  If two extensions
  $E_1 \colon \Lambda \curvearrowright Q_1 \lto[p_1] Q$ and 
  $E_2 \colon \Lambda \curvearrowright Q_2 \lto[p_2] Q$ are 
  equivalent via a quandle isomorphism $\phi \colon Q_1 \to Q_2$,
  then every section $s_1 \colon Q \to Q_1$ induces a section
  $s_2 = \phi \circ s_1 \colon Q \to Q_2$, and by $\Lambda$-equivariance 
  the equation $s_1(a) \ast s_1(b) = f(a,b) \cdot s_1(a\ast b)$ 
  is translated to $s_2(a) \ast s_2(b) = f(a,b) \cdot s_2(a\ast b)$, 
  which means that $\Phi(E_1) = [f] = \Phi(E_2)$, as desired.

  To prove the theorem, we will construct an inverse map 
  $\Psi \colon H^2(Q,\Lambda) \to \Ext(Q,\Lambda)$ as follows.
  Given a quandle $2$-cocycle $f \colon Q \times Q \to \Lambda$, 
  we define the quandle $\tilde{Q} = \Lambda \times_f Q$ as the set 
  $\bigsqcup_{i \in I} \tilde{Q}_i$ with $\tilde{Q}_i = \Lambda_i \times Q_i$ 
  equipped with the binary operation
  \begin{align*}
    (u,a) \ast (v,b) & = \bigl(\; u f(a,b),\; a\ast b \;\bigr).
    \intertext{Idempotency is guaranteed by $f(a,a)=1$, 
      the inverse operation is given by}
    (u,a) \tsa (v,b) &= \bigl(\; u f(a\tsa b,b)^{-1},\; a\tsa b \;\bigr),
    \intertext{and self-distributivity follows from the cocycle condition:}
    \big[ (u,a) \ast (v,b) \big] \ast (w,c) 
    & = \bigl(\; u f(a,b), \; a \ast b \;\bigr) \ast (w,c) \\
    & = \bigl(\; u f(a,b) f( a \ast b, c ), \; (a\ast b) \ast c \;\bigr)
    \qquad\text{and} \\ 
    \big[ (u,a) \ast (w,c) \big] \ast \big[ (v,b) \ast (w,c) ] 
    & = \bigl(\; u f(a,c), \; a \ast c \;\bigr) \ast \bigl(\; v f(b,c), \; b \ast c \;\bigr) \\
    & = \bigl(\; u f(a,c) f(a\ast c,b\ast c), \; (a\ast c)\ast (b\ast c) \;\bigr) .
  \end{align*}
  
  The graded left action of $\Lambda$ on the quandle $\tilde{Q} = \Lambda \times_f Q$ 
  is defined by $\lambda \cdot (u,a) = ( \lambda u, a )$ 
  for all $(u,a) \in \tilde{Q}_i$ and $\lambda \in \Lambda_i$.
  It is straightforward to verify that we thus obtain a graded extension 
  $\Lambda \curvearrowright \Lambda\times_f Q \lto[p] Q$ with projection
  $p(u,a) = a$.

  Suppose that $f,f' \in Z^2(Q,\Lambda)$ are cohomologous, that is, there exists 
  $g \colon Q \to \Lambda$ such that $f'(a,b) = g(a)^{-1} f(a,b) g(a \ast b)$.  
  Then the corresponding extensions are equivalent via the isomorphism 
  $\phi \colon \Lambda \times_{f} Q \to \Lambda\times_{f'} Q$
  defined by $\phi(u,a) = (u g(a),a)$.  Hence we have constructed 
  a well-defined map $\Psi \colon H^2(Q,\Lambda) \to \Ext(Q,\Lambda)$.

  To see that $\Phi\Psi = \id$, let $f \in Z^2(Q,\Lambda)$ and consider 
  the section $s \colon Q \to \Lambda \times_f Q$ with $s(a) = (1,a)$.
  The corresponding $2$-cocycle is $f$, hence $\Phi\Psi = \id$.

  It remains to show that $\Psi\Phi=\id$.  Given an extension 
  $E \colon \Lambda \curvearrowright \tilde{Q} \to Q$, 
  we choose a section $s \colon Q \to \tilde{Q}$ and consider 
  the corresponding $2$-cocycle $f \in Z^2(Q,\Lambda)$.
  The map $\phi \colon \Lambda\times_f Q \to \tilde{Q}$ 
  given by $\phi(u,a) = u \cdot s(a)$ is then an equivalence 
  of extensions, which proves $\Psi\Phi=\id$.

  Naturality and the module structure are easily verified.
\end{proof}

\subsection{The Hurewicz isomorphism} \label{sub:HurewiczIsomorphism}

On the one hand, the Galois correspondence establishes
a natural bijection between quandle extensions 
$E \colon \Lambda \curvearrowright \tilde{Q} \to Q$
and group homomorphisms $\pi_1(Q,q) \to \Lambda$,
see Theorems \ref{thm:ExtensionFundamentalCorrespondence}
and \ref{thm:ExtensionNonConnectedCorrespondence}. 
On the other hand, the preceding cohomology arguments show that 
the second cohomology group $H^2(Q,\Lambda)$ classifies extensions,
see Theorem \ref{thm:ExtensionCohomologyCorrespondence}.
We thus arrive at the following conclusion:

\begin{corollary}
  For every well-pointed quandle $(Q,q)$ and 
  every graded group $\Lambda$ we have natural graded bijections 
  \[
  H^2(Q,\Lambda) \cong \Ext(Q,\Lambda) \cong \Hom(\pi_1(Q,q),\Lambda) .
  \]
  If $\Lambda$ is an abelian group, or more generally a module over some ring $R$, 
  then these objects carry natural $R$-module structures and 
  the bijections are isomorphisms of $R$-modules.
  \qed
\end{corollary}

Finally, we want to prove that $H_2(Q) \cong \pi_1(Q,q)_\mathrm{ab}$.
This is somewhat delicate if $Q$ has infinitely many components:
then the graded group $\pi_1(Q,q)$ is an infinite \emph{product},
whereas $H_2(Q)$ is an infinite \emph{sum} of abelian groups.
The correct formulation is as follows:

\begin{theorem}[Hurewicz isomorphism for quandles] \label{thm:Hurewicz}
  Let $(Q,q)$ be a well-pointed quandle with components $(Q_i,q_i)_{i \in I}$
  and graded fundamental group $\pi_1(Q,q) = \prod_{i \in I} \pi_1(Q,q_i)$.
  Then there exists a natural graded isomorphism 
  $H_2(Q) \cong \bigoplus_{i \in I} \pi_1(Q,q_i)_\mathrm{ab}$.
\end{theorem}

\begin{proof}
  In \secref{sub:TopologicalRealization} we have 
  constructed a $2$-complex $K$ that realizes the fundamental 
  groupoid $\Pi(Q,\Adj(Q)^\circ)$ of a given quandle $Q$, 
  and thus the fundamental group $\pi_1(Q,q_i) \cong \pi_1(K,q_i)$ 
  based at some given point $q_i \in Q$.
  Notice that the connected components of $K$
  correspond to the connected components of $Q$.

  We deduce an isomorphism $H_1(K) \cong H_2(Q)$
  as follows.  The combinatorial chain group $C_1(K)$
  is the free abelian group with basis given by the edges
  of the graph $\Gamma$, which is the $1$-skeleton of $K$.
  On the chain level we can thus define $f \colon C_1(K) \to C_2(Q)$
  by mapping each edge $(a \lto[b] a \ast b) \in C_1(K)$ 
  to the $2$-chain $(a,b) \in C_2(Q)$.  (For the definition
  of quandle homology, see \cite{CarterEtAl:2001} or \cite{Eisermann:2003}).  
  It is readily  verified that this maps $1$-cycles to $2$-cycles and induces 
  the desired isomorphism $H_1(K) \cong H_2(Q)$ on homology.  We conclude that
  \[
  H_2(Q) \cong H_1(K) \cong \oplus_{i \in I} \pi_1(K,q_i)_\mathrm{ab} 
  \cong \oplus_{i \in I} \pi_1(Q,q_i)_\mathrm{ab} 
  \]
  by appealing to the classical Hurewicz Theorem,
  see Spanier \cite[Theorem 7.5.5]{Spanier:1981}.
\end{proof}

\subsection{Application to link quandles} \label{sub:LinkQuandles2}

Having the Hurewicz isomorphism at hand, we can apply it to complete 
our study of links $K \subset \S^3$ and their quandles $Q_K$.
In particular we obtain an explicit correspondence between
the longitude $\ell_K^i \in \pi_1(Q_K,q_K^i)$, 
as explained in \secref{sub:LinkQuandles1},
and the orientation class $[K_i] \in H_2(Q_K)$, 
as explained in \cite[\textsection6.2]{Eisermann:2003}.

\begin{corollary}
  For every choice of base points $q_K^i \in Q_K^i$, the natural 
  Hurewicz homomorphism $h \colon \pi_1(Q_K,q_K) \to H_2(Q_K)$
  is an isomorphism of graded groups, mapping each longitude 
  $\ell_K^i \in \pi_1(Q_K,q_K^i)$ to the orientation class $[K_i] \in H_2(Q_K)$.
\end{corollary}

\begin{proof}
  We know from Theorem \ref{thm:LinkQuandleFundamentalGroup}
  that $\pi_1(Q_K,q_K) = \prod_{i=1}^n \gen{\ell_K^i}$ 
  is abelian, and so $h$ is an isomorphism.
  The longitude $\ell_K^i$ can be read from a link diagram,
  as explained in \cite[Theorem 13]{Eisermann:2003}, as a word
  in the generators of $\pi_K = \Adj(Q_K)$, which corresponds to 
  a path in the complex associated to the link quandle $Q_K$.
  Likewise, the homology class $[K_i] \in H_2(Q_K)$ can be read from 
  the link diagram, as explained in \cite[\textsection6.2]{Eisermann:2003},
  which corresponds to a $1$-cocycle in the same complex.
  The construction of the group homomorphism $h$ in the proof of 
  Theorem \ref{thm:Hurewicz} shows that $h(\ell_K^i) = [K_i]$.
\end{proof}

Consider two oriented links $K = K_1 \sqcup \dots \sqcup K_n$
and $K' = K'_1 \sqcup \dots \sqcup K'_n$ in $\S^3$, and 
their respective link quandles $Q_K$ and $Q_{K'}$.  
We have a natural bijection $\pi_0(K) \isoto \pi_0(Q_K)$.
Every quandle isomorphism $\phi \colon Q_K \isoto Q_{K'}$ 
induces a bijection $\tau \colon \pi_0(Q_K) \isoto \pi_0(Q_{K'})$
as well as a graded isomorphism $\phi_* \colon H_2(Q_K) \isoto H_2(Q_{K'})$. 
We also know that for each $i$ the group $H_2(Q_K)_i = \gen{[K_i]}$ 
is either trivial or freely generated by $[K_i]$, and the same holds 
for its isomorphic image $H_2(Q_K)_{\tau i} = \gen{[K'_{\tau i}]}$.  
This means that $\phi_* [K_i] = \pm [K'_{\tau i}]$ for all $i$. 



\begin{theorem}
  Two oriented links $K = K_1 \sqcup \dots \sqcup K_n$
  and $K' = K'_1 \sqcup \dots \sqcup K'_n$ in $\S^3$ 
  are ambient isotopic respecting orientations and numbering 
  of components if and only if there exists a quandle 
  isomorphism $\phi \colon Q_K \isoto Q_{K'}$ such that 
  $\phi_* [K_i] = [K'_i]$ for all $i=1,\dots,n$.
\end{theorem}

\begin{proof}
  Obviously, if $K$ and $K'$ are ambient isotopic,
  then the quandles $Q_K$ and $Q_{K'}$ are isomorphic.
  Conversely, consider an isomorphism $\phi \colon Q_K \isoto Q_{K'}$ 
  such that $\phi_* [K_i] = [K'_i]$ for all $i=1,\dots,n$.
  According to the characterization of trivial components
  in Corollary \ref{cor:TrivialComponents}, we can assume
  that all components of $K$ and $K'$ are non-trivial.
  We number the components $Q_K^1,\dots,Q_K^n$ of $Q_K$ 
  such that $[K_i] \in H_2(Q_K)$ is supported by $Q_K^i$.
  We choose a base point $q_K^i \in Q_K^i$ for each $i=1,\dots,n$.
  In the adjoint group $\Adj(Q_K)$ this determines group elements 
  $m_K^i = \adj(q_K^i)$.  For each $i$ there are two generators 
  $(\ell_K^i)^\pm \in \pi_1(Q_K,q_K^i)$ of the fundamental group, and 
  we choose $\ell_K^i$ corresponding to the given class 
  $[K_i] \in H_2(Q_K)$ under the Hurewicz isomorphism.
  In this way we recover the link group $\pi_K = \Adj(Q_K)$ 
  together with the peripheral data $(m_K^i,\ell_K^i)$ for each 
  link component $K_i$.  The quandle isomorphism $\phi \colon Q_K \isoto Q_{K'}$
  thus induces a group isomorphism $\psi \colon \pi_K \isoto \pi_{K'}$ 
  respecting the peripheral data.  According to Waldhausen's result 
  \cite[Corollary 6.5]{Waldhausen:1968}, there exists an orientation 
  preserving homeomorphism $f \colon (\S^3,K) \isoto (\S^3,K')$ such that 
  $f_* = \psi$; for details see \cite[Theorem 6.1.7]{Kawauchi:1996}.
  Moreover, $f$ can be realized by an ambient isotopy.
\end{proof}


\bibliographystyle{plain}
\bibliography{qcovering}

\end{document}

%% file: qcircles.pstex_t
\begin{picture}(0,0)%
\includegraphics{qcircles.pstex}%
\end{picture}%
\setlength{\unitlength}{4144sp}%
\begingroup\makeatletter\ifx\SetFigFont\undefined%
\gdef\SetFigFont#1#2#3#4#5{%
  \reset@font\fontsize{#1}{#2pt}%
  \fontfamily{#3}\fontseries{#4}\fontshape{#5}%
  \selectfont}%
\fi\endgroup%
\begin{picture}(2733,2224)(1441,-2162)
\put(1441,-736){\makebox(0,0)[lb]{\smash{{\SetFigFont{10}{12.0}{\rmdefault}{\mddefault}{\updefault}{\color[rgb]{0,0,0}$\Z$}%
}}}}
\put(3241,-736){\makebox(0,0)[lb]{\smash{{\SetFigFont{10}{12.0}{\rmdefault}{\mddefault}{\updefault}{\color[rgb]{0,0,0}$\Z$}%
}}}}
\put(3241,-2086){\makebox(0,0)[lb]{\smash{{\SetFigFont{10}{12.0}{\rmdefault}{\mddefault}{\updefault}{\color[rgb]{0,0,0}$\Z_5$}%
}}}}
\put(1441,-2086){\makebox(0,0)[lb]{\smash{{\SetFigFont{10}{12.0}{\rmdefault}{\mddefault}{\updefault}{\color[rgb]{0,0,0}$\Z_6$}%
}}}}
\end{picture}%

%% file: qhomotopy.pstex_t
\begin{picture}(0,0)%
\includegraphics{qhomotopy.pstex}%
\end{picture}%
\setlength{\unitlength}{4144sp}%
\begingroup\makeatletter\ifx\SetFigFont\undefined%
\gdef\SetFigFont#1#2#3#4#5{%
  \reset@font\fontsize{#1}{#2pt}%
  \fontfamily{#3}\fontseries{#4}\fontshape{#5}%
  \selectfont}%
\fi\endgroup%
\begin{picture}(2624,958)(499,-281)
\put(631,479){\makebox(0,0)[b]{\smash{{\SetFigFont{10}{12.0}{\rmdefault}{\mddefault}{\updefault}{\color[rgb]{0,0,0}$a^{\pm 1}$}%
}}}}
\put(631,-61){\makebox(0,0)[b]{\smash{{\SetFigFont{10}{12.0}{\rmdefault}{\mddefault}{\updefault}{\color[rgb]{0,0,0}$a$}%
}}}}
\put(1531,299){\makebox(0,0)[b]{\smash{{\SetFigFont{10}{12.0}{\rmdefault}{\mddefault}{\updefault}{\color[rgb]{0,0,0}$b^{+1}$}%
}}}}
\put(1531,-151){\makebox(0,0)[b]{\smash{{\SetFigFont{10}{12.0}{\rmdefault}{\mddefault}{\updefault}{\color[rgb]{0,0,0}$b^{-1}$}%
}}}}
\put(2206,164){\makebox(0,0)[b]{\smash{{\SetFigFont{10}{12.0}{\rmdefault}{\mddefault}{\updefault}{\color[rgb]{0,0,0}$c$}%
}}}}
\put(2611,-241){\makebox(0,0)[b]{\smash{{\SetFigFont{10}{12.0}{\rmdefault}{\mddefault}{\updefault}{\color[rgb]{0,0,0}$b$}%
}}}}
\put(3016,164){\makebox(0,0)[b]{\smash{{\SetFigFont{10}{12.0}{\rmdefault}{\mddefault}{\updefault}{\color[rgb]{0,0,0}$c$}%
}}}}
\put(2611,569){\makebox(0,0)[b]{\smash{{\SetFigFont{10}{12.0}{\rmdefault}{\mddefault}{\updefault}{\color[rgb]{0,0,0}$b \ast c$}%
}}}}
\end{picture}%